\documentclass[english,reqno,11pt,empty]{amsart}
\usepackage{amsmath, amsthm, amsfonts}
\usepackage{hyperref}
\hypersetup{
	colorlinks=true,
			citecolor=blue!60!black,
	linkcolor=red!60!black,
		urlcolor=green!40!black,
		filecolor=yellow!50!black,
	breaklinks=true,
	pdfpagemode=UseNone,
	bookmarksopen=false,
}
\usepackage{tcolorbox}
\usepackage{amsmath,amsthm,amssymb}
\usepackage{amscd,indentfirst,epsfig}
\usepackage{latexsym}
\usepackage{times}
\usepackage{enumerate}
\usepackage{mathrsfs}
\usepackage{stmaryrd}
\usepackage{amsopn}
\usepackage{amsmath}
\usepackage{amssymb,dsfont,mathtools}
\usepackage{amsfonts,bm}
\usepackage{amsbsy,amsmath}
\usepackage{amscd}
\usepackage{xcolor}
\usepackage{mathtools}

\usepackage[mono=false]{libertine}
\usepackage[T1]{fontenc}
\usepackage{amsthm}
\usepackage[normalem]{ulem} 
\usepackage[cal=euler, scr=boondoxo]{}
\usepackage{microtype}
\usepackage{physics}
\usepackage{numprint}

 \usepackage[mediumspace,mediumqspace,Grey,squaren]{SIunits}

 \usepackage{array}

\usepackage{latexsym}
\usepackage{enumerate}
\usepackage{mathrsfs}
\usepackage{stmaryrd}
\usepackage{amsopn}
\usepackage{amsmath,amsthm,amssymb}
\usepackage{mathtools}
\usepackage{amsfonts}
\usepackage{soul}
\usepackage{amsbsy}
\usepackage{amscd,indentfirst,epsfig}
\usepackage{amsfonts,latexsym,verbatim,amsbsy}

\usepackage{pgf,tikz}
\usepackage{amsmath,amsthm,amssymb}
\usepackage{amscd,indentfirst,epsfig}
\usepackage{latexsym}
\usepackage{times}
\usepackage{enumerate}
\usepackage{mathrsfs}
\usepackage{stmaryrd}
\usepackage{amsopn}
\usepackage{amsmath}
\usepackage{amssymb}
\usepackage{amsfonts,bm,dsfont}

\usepackage{dsfont,mathtools}
\usepackage{enumerate}
\usepackage{mathrsfs,esint}
\usepackage{stmaryrd}

\usepackage{amsopn,bm}
\usepackage{amsmath}
\usepackage{amssymb}
\usepackage{amsfonts}
\usepackage{amsbsy}
\usepackage{amscd,indentfirst,epsfig}
\usepackage{amsfonts,amsmath,latexsym,amssymb,verbatim,amsbsy,textcomp}
\usepackage{amsthm}
\usepackage{dsfont}
\usepackage{colordvi}
\usepackage{pstricks}
\usepackage{csquotes}

\usepackage{natbib}
\def \leq {\leqslant}                                                                
\def \geq {\geqslant}                                                                
\def\ind#1{\lower5pt\hbox{$\scriptstyle #1$}}                                        
\def \d {\mathrm{d}}                                                                 
\def \ds {\displaystyle}                                                             
\def \R{\mathbb R}                                                                   
\def\S{\mathbb S}                                                                    
\def\Z{\mathbb Z}                                                                    
\def\N{\mathbb N}                                                                    
\def\Cs{\mathcal C}																	 
\def\A{\mathcal A}																	 
\def\D{\mathscr D}																	 
\def\M{\mathcal M}																	 
           \usepackage{empheq}     																	 
\def \T {\mathbb{T}}                                                                 
\numberwithin{equation}{section}                                                     
\newcommand{\vertiii}[1]{{\left\vert\kern-0.25ex\left\vert\kern-0.25ex\left\vert #1  
    \right\vert\kern-0.25ex\right\vert\kern-0.25ex\right\vert}}                      
\newcommand{\verti}[1]{{\left\vert\kern-0.25ex\left\vert\kern-0.25ex\left\vert #1    
    \right\vert\kern-0.25ex\right\vert\kern-0.25ex\right\vert}}						 
\def\e{\varepsilon}      															    
\def \m {\bm{\varpi}}															        
			\def \l {\ell}																		        
\def \ho {h^{0}} 																	    
\def \hu {h^{1}}																		
\def \huu {\Psi}																		
\def \Mo {\Delta_{0}}																	
\def \ra {\Big\rangle}																	
\def \la {\Big\langle}																	
\def \en {\theta}																	
\def \vE {\vartheta}	
\def \lae {\lambda_{\e}}													            
\def \g {\gamma} 
\def \re {{\rm e}}																		

\def\E{\mathcal{E}} 																	
\def \B{\mathcal{B}}                                         							
\def\W {\mathbb{W}} 
\def \H {\mathcal{H}} 

\def\vet{v_{\ast}}                                           							

\def \vb {v_\ast}                                           							
\def \Q{\mathcal{Q}}                                        							
\def \LL {\mathscr{L}_{\re}} 							 							    
		
\def \n {{n}}													 							
\def \LLe {\mathscr{L}_{\e,t}}						 							    
\def \G {\mathcal{G}}	
\def \LL {\mathscr{L}_{1}} 							 							    
 	
\def \g {\gamma}							 										
\def \It {\int_{\R^{3} \times \S^{2}}}                       							
\def \IR {\int_{\R^{3}}}                                       							

\newtheorem{theo}{Theorem}[section]                          							
\newtheorem{prop}[theo]{Proposition}                         							
\newtheorem{cor}[theo]{Corollary}                            							
\newtheorem{lem}[theo]{Lemma}                                							
\newtheorem{defi}[theo]{Definition}                          							
\newtheorem{nb}[theo]{Remark}                                							
\newtheorem{exe}[theo]{Example}                          							
 
\begin{document}

\title[Hydrodynamic limit of the viscoelastic Boltzmann equation for hard-spheres]{The hydrodynamic limit of viscoelastic granular gases}

\author{Ricardo {\sc Alonso}}
\address[R. Alonso]{College of Science and Engineering, Division of Sciences, HBKU, Education City, Doha, Qatar.} \email{ralonso@hbku.qatar.edu}

\author{Bertrand {\sc Lods}}

\address[B. Lods]{Universit\`{a} degli Studi di Torino \& Collegio Carlo Alberto, Department of Economics and Statistics, Corso Unione Sovietica, 218/bis, 10134 Torino, Italy.}\email{bertrand.lods@unito.it}

\author{Isabelle {\sc Tristani}}
\address[I. Tristani]{Université Côte d’Azur, CNRS, LJAD, Parc Valrose, F-06108 Nice, France} 
\email{isabelle.tristani@univ-cotedazur.fr}

\maketitle
\begin{abstract} 
We obtain the first rigorous derivation of an incompressible Navier-Stokes-Fourier system with self-consistent and time-dependent forcing terms from the inelastic hard-spheres Boltzmann equation associated to the relevant case of viscoelastic granular gases.  The model's inelasticity is measured by the so-called restitution coefficient which, for viscoelastic particles, depends on the relative velocities of particles.  Through a suitable self-similar change of variables, a balanced dynamic between energy inflow and outflow naturally emerges in the model which permits its analysis.  In contrast, such balanced dynamic does not emerge naturally in the constant restitution case and has to be imposed in \cite{ALT}.  The exact self-similar rescaling, which allows to capture nontrivial inelastic-hydrodynamic effects, presents itself explicitly in terms of the Knudsen number and the restitution coefficient.  The consequence of such scaling is a non-autonomous rescaled Boltzmann equation whose solutions converge, in a specific weak sense, towards the aforementioned hydrodynamic limit.  The incompressible Navier-Stokes-Fourier system obtained by this process appears to be new in this context.  As a byproduct of the analysis, we determine the exact dissipation rate of the granular temperature known as \emph{Haff's law}.
\end{abstract}

\vspace{.5cm}

\vspace{0.4cm}
 
 {
  \hypersetup{linkcolor=blue}
 }

\tableofcontents

\section{Introduction}

The study of granular flows, encompassing phenomena such as avalanches, sands, grain mixing, presents a rich and challenging area of research.  These flows, composed of a large number of discrete particles interacting through inelastic collisions, friction, heating, exhibit complex behaviors that go beyond the classical theory of elastic gases.  The non-equilibrium nature of dilute granular flows, characterized by energy dissipation through inelastic collisions, requires the development of sophisticated theoretical frameworks to accurately predict their macroscopic properties. 

In this paper, we address the delicate problem of deriving hydrodynamic equations for a system of \emph{viscoelastic hard spheres} where the inelastic collisions are modelled by a velocity-dependent restitution coefficient.  Moving beyond traditional approaches that rely on simplified constitutive relations or assume a constant restitution coefficients \cite{Gold03, Poschel}, our work provides a rigorous derivation of an incompressible Navier-Stokes-Fourier system directly from the Boltzmann equation.  More specifically, we employ a novel self-similar change of scale that demonstrates a controlled dissipation process involving both nearly elastic and highly inelastic collisions.   In such scale, a new Boltzmann equation, suited for capturing inelastic and hydrodynamic effects, emerges.  A rigorous proof of convergence for this model towards a modified incompressible Navier-Stokes-Fourier system is obtained.  The approach presented here provides new insights into the macroscopic quantities of such system. We refer the reader to \cite[Section 1.5]{ALT} for a detailed account of the multiscale nature of granular gases and a thorough review of the relevant literature on the subject.

\subsection{The Boltzmann equation for granular gases}
The inelastic Boltzmann equation, which provides a statistical description of identical hard spheres under binary inelastic collisions, is given by
\begin{equation}\label{Bol-}
\partial_{t}F(t,x,v) + v\cdot \nabla_{x} F(t,x,v)=\Q_{\re}(F,F)
\end{equation}
supplemented with initial condition $F(0,x,v)=F_{\mathrm{in}}(x,v)$, where $F(t,x,v)$ is the density of particles with position $x \in \T_{L}^{3}$ and velocity $v \in \R^{3}$ at time $t\geq0$. The quadratic collision operator~$\Q_\re(\cdot,\cdot)$ models the interactions of hard-spheres by inelastic  binary collisions, described later on. At this point, we point out that the case $\re=1$ corresponds to the classical Boltzmann equation for elastic interactions. 

The dynamic of the Boltzmann equation \eqref{Bol-} is assumed to take place, for simplicity, in the flat torus
\begin{equation*}\label{torus}
\T_{L}^{3}=\R^{3}\slash (2\pi L\,\Z)^{3}\end{equation*}
for some typical length-scale $L >0$. This corresponds to periodic boundary conditions
$$F(t,x+2\pi L\bm{e}_{i},v)=F(t,x,v)\,, \qquad \forall \, i=1,2,3\,.$$
To capture the hydrodynamic behavior of the gas, we write the above equation in \emph{nondimensional form} introducing the dimensionless Knudsen number
$$\e:=\dfrac{\text{ mean free path }}{\text{ spatial length-scale }}$$
which is assumed to be small.  We introduce a rescaling of time-space to capture the hydrodynamic limit and introduce the particle density
\begin{equation}\label{eq:Scaling}
F_{\e}(t,x,v)=F\left(\frac{t}{\e^{2}},\frac{x}{\e},v\right), \qquad \forall \,t \geq 0\,.
\end{equation}
In this case, we choose for simplicity 
$L=\varepsilon^{-1}$  which ensures that $F_{\e}$ is defined on $\R^{+}\times \T^{3}\times \R^{3}$ with $\T^{3}
:=\T_{1}^{3}$.  In the sequel, we assume for simplicity that the torus $\T^{3}$ is equipped with the normalized Lebesgue measure, i.e.~$|\T^{3}|= 1$.  It is well-known that, in the classical elastic case, this scaling leads to the incompressible Navier-Stokes system.   We note that other scalings are possible yielding different hydrodynamic models. Under such a scaling, the typical number of collisions per particle per time unit is~$\e^{-2}$, more specifically, $F_{\e}$ satisfies the rescaled Boltzmann equation
\begin{subequations}
\begin{equation}\label{Bol-e}
\begin{aligned}
&\e^{2}\partial_{t}F_{\e}(t,x,v) + \e\,v\cdot \nabla_{x} F_{\e}(t,x,v)=\Q_{\re}(F_{\e},F_{\e})\,, \qquad (x,v) \in \T^{3}\times\R^{3}\,, \\
\end{aligned}
\end{equation}
supplemented with the initial condition
\begin{equation}\label{eq:init}
F_{\e}(0,x,v)=F_{\mathrm{in}}^{\e}(x,v)=F_{\mathrm{in}}\left(\frac x\e, v\right)\,.
\end{equation}
\end{subequations}
As mentioned before $\Q_{\re}$ denotes the collision operator modelling the interactions of hard-spheres by  inelastic binary collisions.  The inelasticity is described by the restitution coefficient $\re\in [0,1]$ defined as the ratio between the magnitude of the relative velocity normal component, along the separation line between the colliding sphere's centers at contact, after and before the collision. Namely, if $v$ and $\vb$  denote the velocities of two colliding particles before collision, their respective velocities $v'$ and $\vb'$ after collision are such that
\begin{equation}\label{coef}
(u'\cdot \n)=-(u\cdot \n) \,\re(|u \cdot \n|)\,,
\end{equation}
where the vector $\n \in \mathbb{S}^2$  determines the impact direction, that is, $\n$ stands for the unit vector that points from the $v$-particle center to the $\vb$-particle center at the instant of impact.  Here above
$$u:=v-\vb\,,\qquad u':=v'-\vb'\,,$$
denote respectively the relative velocity before and after collision.  Assuming  the granular particles to be perfectly smooth hard-spheres of mass
$m=1$ and that the momentum is preserved during the collision process, that is 
\begin{equation}\label{eq:conse}v'+\vb'=v+\vb\,,\end{equation}
the velocities after collision $v'$ and $\vb'$ are given, in virtue of \eqref{coef} and~\eqref{eq:conse}, by
\begin{equation}
\label{transfpre}
  v'=v-\frac{1+\re(|u \cdot \n|)}{2}\,(u\cdot \n)\n\,,
\qquad \vb'=\vb+\frac{1+\re(|u \cdot \n|)}{2}\,(u\cdot \n)\n\,.
\end{equation}
Notice that the kinetic energy is dissipated since
\begin{equation}\label{eq:diss}
|v'|^{2}+|\vb'|^{2}-|v|^{2}-| \vb|^{2}=-\frac{1-\re^{2}\big(|u \cdot n|\big)}{2}\,\big(u \cdot n\big)^{2}\leq 0\,.
\end{equation}
We are interested in viscoelastic restitution coefficients.  Such coefficients model the dissipation associated to the elastic deformation during collisions, thus, small relative velocities are associated to quasi elastic collisions and large relative velocities to highly inelastic collisions.  In Definition~\ref{defiC}, we make this concept more rigorous.  In particular, we assume that there exist positive constants~$\overline{\g} > \g >0$ and $\mathfrak{a}_{0},\mathfrak{b}_{0} >0$ such that
\begin{equation}\label{eq:ab00}
\left|\re(r)-1-\mathfrak{a}_{0}r^{\g}\right| \leq \mathfrak{b}_{0}r^{\overline\g}\,, \qquad \forall \,r >0
\end{equation}
which in particular provides the behaviour of $\re(r)$ for $r\simeq 0$.   As a particular case of our theory, we cover the case of viscoelastic hard spheres, refer to~\cite{A} and~{\cite{ALCMP}} for an additional discussion of viscoelastic restitution coefficients.
The collision operator $\Q_{\re}(\cdot,\cdot)$ is defined in weak form as
\begin{equation}\label{eq:weakN}
\int_{\R^{3}}\Q_{\re}(g,f)(v)\,\psi(v)\,\d v= \int_{\R^{3}\times \R^{3}\times \S^{2}}f(v)g(\vb)B_{0}(v-\vb,n)\left[\psi(v') -\psi(v) \right]\,\d v\,\d \vb\,\d \n
\end{equation}
where $\d \n$ denotes the Lebesgue measure on $\S^{2}$.  The collision kernel is given by
\begin{equation} \label{def:B0}
B_{0}(u,\n)=|u|\,b_{0}(\widehat{u}\cdot \n)\,, \qquad u=v-\vb\,, \qquad \widehat{u}:=\frac{u}{|u|}
\end{equation}
with an angular kernel $b_{0}(\cdot)$ which is assumed to satisfy, in addition to suitable technical assumptions, the Grad angular cut-off assumption
\begin{equation} \label{eq:cutoff}
\int_{\S^{2}}|\widehat{u}\cdot \n|^{-1}b_{0}(\widehat{u}\cdot\n)\,\d \n < \infty\,.
\end{equation}
In Sections \ref{Sec21} and \ref{Sec22} we describe an alternative representation of the collision mechanism as well as the strong form of the collision operator $\Q_{\re}$.  Additionally, the main assumptions on the restitution coefficient $\re(\cdot)$ are clearly stated.  
 
Consequences of the collision mechanism can be seen not only at the microscopic level but also at the macroscopic scale. More precisely, let us introduce the \textit{global} density and bulk velocity respectively defined by 
$$
\bm{R}_{\e}(t) := \int_{\T^{3}\times \R^{3}}F_{\e}(t,x,v)\, \d v\, \d x 
\qquad \text{and} \qquad 
\bm{U}_{\e}(t) := \int_{\T^{3}\times\R^{3}}v F_{\e}(t,x,v)\, \d v\, \d x\,.
$$
Due to \eqref{eq:weakN} and \eqref{eq:conse} one observes that these quantities are preserved over time, namely,
$$
\dfrac{\d}{\d t}\bm{R}_{\e}(t)= \frac{\d}{\d t}\bm{U}_{\e}(t)=0\,.
$$
Consequently, there is no loss of generality in assuming that
$$\bm{R}_{\e}(t)=\bm{R}_{\e}(0)=1\,, \qquad \bm{U}_{\e}(t)=\bm{U}_{\e}(0)=0\,, \qquad \forall \, t \geq0\,.$$
The main contrast between elastic and inelastic gases is that the \emph{granular temperature}
\begin{equation} \label{def:Teps}
\bm{T}_{\e}(t):=\int_{\T^{3}\times\R^{3}}|v|^{2}F_{\e}(t,x,v)\, \d v\, \d x\,, \qquad \forall \,t \geq0\,, \qquad \forall \, \e >0\,,
\end{equation}
is constantly decreasing
\begin{equation}\label{eq:TeHaff}
\frac{\d}{\d t}\bm{T}_{\e}(t)= -\frac{1}{\e^{2}} \mathscr{D}_{\re}(F_{\e}(t),F_{\e}(t))\,,\qquad \forall \,t \geq0\,,
\end{equation}
where $\mathscr{D}_{\re}(\cdot,\cdot)$ denotes the normalised energy dissipation associated to $\Q_{\re}$, see \cite{CMS}.   The energy dissipation is defined for suitable $f,g$ by
\begin{equation}\label{eq:Dre}
\mathscr{D}_{\re}(f,g):=\int_{\T^{3}}\d x\int_{\R^{3} \times \R^{3}}f(x,\vb)g(x,v)\bm{\Psi}_{\re}(|v-\vb|^2)\, \d v\, \d \vb
\end{equation} 
with 
\begin{equation} \label{def:Psire}
\bm{\Psi}_{\re}(r):= \frac{r^{\frac{3}{2}}}{2} \int_0^1 \big(1-\re(\sqrt{r}z)^2\big)z^3 \, \d z\,, \qquad \forall \, r >0\,. 
\end{equation}
We refer to Section \ref{Sec21} for more details. In addition, it is possible to show that, for any fixed $\e >0$ (small enough), 
$$\bm{T}_{\e}(t)\xrightarrow[t \to \infty]{} 0$$
which expresses the \emph{total cooling of granular gases}.  Determining the exact dissipation rate of the granular temperature is an important question known as \emph{Haff's law}, see {\cite{haff}}, which we also address here in the context of inhomogeneous viscoelastic granular flows. 
\subsection{Related literature and strategy}\label{sec:liter}
The derivation of hydrodynamic limits from kinetic equations, particularly the elastic Boltzmann equation, has been a central theme in kinetic theory for decades.  For the elastic Boltzmann equation, where collisions conserve kinetic energy, significant progress has been made in establishing rigorous connections between the Boltzmann equation and the compressible and incompressible Navier-Stokes equations following the program developed in \cite{BaGoLe1,BaGoLe2} which found later a rigorous justification in the works \cite{golseSR,golseSR1}. We refer the reader to \cite{SR,golse} for an up-to-date description of the relevant results in the field as well as to the introduction of \cite{ALT} for a description of the main approaches to achieve a rigorous derivation of hydrodynamic limits, for example using a spectral approach or renormalized solutions.  

Extending these results to inelastic interactions presents significant challenges; notably, the existence of renormalized solutions for the inhomogeneous inelastic Boltzmann equation remains an open problem due to the absence of known entropy control in the large data regime. We just mention here that the only rigorous result in the diffusive regime is the one obtained in~\cite{ALT}, we also mention that a rigorous result has been obtained by~\cite{jabin} in the mono-kinetic regime (which considers the extreme case of infinite energy dissipation rate) leading to the pressureless Euler system (corresponding to
sticky particles) in a one-dimensional framework.

The inelastic Boltzmann equation has been the subject of extensive research due to its relevance in the study of granular gases \cite{Poschel, garzo}.  A simplified approach involves assuming a constant restitution coefficient. While this makes the analysis more tractable, it may not accurately capture the behavior of real granular materials in the time-asymptotic.  The case of viscoelastic hard spheres, where the restitution coefficient depends on the impact velocity \cite{PoSc}, presents a more realistic albeit more challenging scenario. 

This work builds upon various attempts to understand the behaviour of viscoelastic hard spheres.  Besides the derivation of the model in \cite{Poschel,PoSc} with first principles, a systematic study of viscoelastic granular gases has been performed by the authors in a series of results mainly concerned with spatially homogeneous situations \cite{A,ALhaff,ALCMP} or considering thermal bath heating \cite{Tr,CMS}.  In particular, the work \cite{Tr} provides, to the best of our knowledge, the only Cauchy theory for spatially inhomogeneous Boltzmann equation associated to viscoelastic particles in the presence of a thermal bath. This work extends these previous efforts by providing a rigorous derivation of hydrodynamic equations for viscoelastic hard spheres.  We adapt the roadmap developed in \cite{ALT} which addressed the case of constant restitution in a nearly elastic regime.  A crucial difference, however, lies in the fact that in the viscoelastic case a natural energy inflow-outflow balance emerges in the model after self-similar rescaling.  This is not the case for the constant restitution model where such balance has to be enforced by a direct modification of the restitution coefficient associated to a nearly elastic regime.

The proof's strategy is based upon the following steps:
\begin{enumerate}[a)]
    \item \textit{\textbf{Rescaling and self-similar variables:}} A suitable self-similar rescaling of the Boltzmann equation based on the Knudsen number and the restitution coefficient is introduced.  This scale allows to capture the inelastic effects along the hydrodynamic limit through a non-autonomous Boltzmann equation.
    \item \textit{\textbf{Linearization and hypocoercivity:}}  The non-autonomous Boltzmann equation is linearized around a Maxwellian distribution, and the resulting linearized operator is analyzed. Hypocoercivity techniques are then employed to investigate the long-time behavior of solutions.
    \item \textit{\textbf{A Priori Estimates and compactness:}} Key to our approach is the derivation of \emph{a priori} estimates that bound the solutions in suitable Sobolev spaces. These \emph{a priori} estimates are obtained using energy method applied to a specific splitting of the Boltzmann equation under study.  Such a splitting, already used in \cite{ALT}, is reminiscent to the one introduced in \cite{bmam} for the elastic Boltzmann equation.  The \emph{a priori} estimates and compactness arguments allow us to prove the existence of a hydrodynamic limit.
    \item\textit{\textbf{Identification of the hydrodynamic limit:}} Finally, the hydrodynamic limiting equations satisfied by the macroscopic quantities are identified.  They correspond to a modified incompressible Navier-Stokes-Fourier system with self-consistent forcing terms.
\end{enumerate}
\subsection{Self-similar change of variables}\label{sec:SSCH} In order to understand the free-cooling inelastic Boltzmann equation~\eqref{Bol-e}, we perform a self-similar change of variables which allows to introduce an intermediate asymptotic regime, see \cite{ALCMP}.   For $F_\e$ solving~\eqref{Bol-e}-\eqref{eq:init} we write 
\begin{equation}
\label{eq:Fescal}
F_{\e}(t,x,v)={V}_{\e}(t)^{3}f_{\e}\big(\tau_{\e}(t),x,V_{\e}(t)v\big)\,,
\end{equation}
where $V_{\e}(t) \geq0$ and $\tau_{\e}(t) \geq 0$ is a new time scaling.  Assuming that 
\begin{equation} \label{eq:taueps}
\tau_{\e}(t)=\int_{0}^{t}\frac{\d s}{V_{\e}(s)}\,, \end{equation}
it follows that $f_{\e}$ solves the equation
\begin{equation}\label{BE0-In}
\e^{2}\partial_{t} f_{\e}(t,x,v)+\e v \cdot \nabla_{x} f_{\e}(t,x,v) + \e^{2}\xi_{\e}(t)\,\mathrm{div}_v (v f_{\e}(t,x,v))=  \Q_{\re_{\e,t}}(f_{\e},f_{\e})\,,
\end{equation}
with initial condition 
$$f_{\e}(0,x,v)=F^{\e}_{\mathrm{in}}(x,v)\,.$$
The change of variable \eqref{eq:Fescal} induces a rescaled restitution coefficient, namely from $\re(\cdot)$ to $\re_{\e,t}(\cdot)$. This is an important difference with the case of constant restitution coefficient \cite{ALT} for which the restitution coefficient was independent to a self-similar change of variable.
Notice that the definitions of $\xi_\e(t)$ and $\re_{\e,t}(\cdot)$ in~\eqref{BE0-In} only depend on~$V_\e(t)$. The choice of the scaling functions $V_{\e}(t)$, and consequently $\tau_{\e}(t),\xi_{\e}(t)$ and~$\re_{\e,t}(\cdot)$ constitutes one of the milestones to derive a hydrodynamic limit.  The details leading to the choice of $V_{\e}(t)$ are given in Section \ref{sec:scaling} and are based upon explicit sharp quantitative estimates between $\Q_{\re_{\l}}$ and $\Q_{1}$ for $\l \simeq 0$ where
\begin{equation} \label{def:e_ell}\re_{\l}(r):=\re(\l r)\,, \qquad \forall \,\l, r >0\,.
\end{equation}
 
Indeed, the strategy is based on a comparison of \eqref{BE0-In} with its classical elastic counterpart, involving $\Q_{1}$ and no drift term, for which the analysis is well understood.  At the formal level we seek that
$$\Q_{\re_{\e,t}}(f,f) \simeq \Q_{1}(f,f) \quad \text{and} \quad  \e^{2}\xi_{\e}(t) \simeq 0 \quad \text{as} \quad \e \simeq 0$$
in a way that a quasi-conservative regime is predominant in order $1$ time-space variations.   We note that the drift term in \eqref{BE0-In} effectively introduces energy in the rescaled granular gas, thus, a balance between the input-output energy fluxes along the hydrodynamic limit is a pre-requisite for the correct self-similar scaling.
As time evolves, a quasi-conservative regime is reached.  This time scale depends on $\e$ since the number of collisions per unit time plays a central role in the dissipation process.
Let us write the distribution $f_\e$ as
\begin{equation}\label{ansatz}
f_{\e}(t,x,v)=\M(v)+ \e h_{\e}(t,x,v)
\end{equation}
with $\M$ the universal Maxwellian distribution with unit mass, zero momentum, and energy~ $3\en_\star$
\begin{equation} \label{eq:max}
\M(v) :=(2\pi\en_{\star})^{-\frac{3}{2}}\exp\left(-\frac{|v|^{2}}{2\en_{\star}}\right)\,, \quad \forall\,v \in \R^{3}\,,
\end{equation} 
where $\en_{\star} >0$ is such that
\begin{equation}\label{eq:FinEe}\lim_{\e\to0}\frac{1}{\e}\int_{\T^{3}\times \R^{3}}\left(f_{\e}(0,x,v)-\M(v)\right)|v|^{2}\,\d v\,\d x=\lim_{\e\to0}\int_{\T^{3}\times \R^{3}} h_\e(0,x,v) \,|v|^2 \d v \, \d x = 0\,.
\end{equation}
Under the \emph{ansatz} \eqref{ansatz}, noted in \cite{ALCMP} for the time-asymptotic analysis, the equation on the fluctuation $h_\e$ writes
\begin{equation}\label{eq:heps} 
\partial_t h_\e + \frac1\e v \cdot \nabla_x h_\e + \xi_{\e}(t) \mathrm{div}_v (v h_\e) =
\frac{1}{\e^2} \LLe h_\e + \frac1\e \Q_{\e,t} (h_\e,h_\e) + S_\e(t,v)
\end{equation}
 where $\Q_{\e,t}:=\Q_{\re_{\e,t}}$.  The short-hand notations for the linear operator $\LLe$ and the source term~$S_\e$ are given by
\begin{equation} \label{def:LetSe}
\LLe h := \Q_{\e,t}(\M,h) + \Q_{\e,t}(h,\M)\,, \qquad  S_\e(t,v) := \frac{1}{\e^3} \Q_{\e,t}(\M,\M) - \frac1\e \xi_{\e}(t)\mathrm{div}_v(v \M)\,.
\end{equation}

\subsection{Notations and definitions}\label{sec:nota}

We introduce here functional spaces used along the document. For any nonnegative weight function $m\::\:\R^{3}\to \R^{+}$, all norm weights considered will depend only on the velocity, define $L^q_v L^p_x(m)$, $1 \leq p,q \leq +\infty$, as the Lebesgue space associated to the norm 
$$
\| h \|_{L^q_v L^p_x(m)} = \| \| h(\cdot, v)\|_{L^p_x} \, m(v) \|_{L^q_v}\,.
$$
Consider also the standard higher-order Sobolev generalizations $\W^{\sigma, q}_v \W^{s,p}_x(m)$ for any $\sigma, s \in \N$ defined by the norm 
$$
\| h \|_{\W^{\sigma, q}_v \W^{s,p}_x(m)} = \sum_{\substack {0\leq s' \leq s, \, 0 \leq \sigma' \leq \sigma, \\ s'+ \sigma' \leq \max(s,\sigma)}} 
\| \| \nabla_x^{s'} \nabla_v^{\sigma'} h(\cdot, v)\|_{L^p_x} \, m(v) \|_{L^q_v}\,.
$$
This definition reduces to the classical weighted Sobolev space $\W^{s,p}_{x,v}(m)$ when $q=p$ and $\sigma=s$. For $m \equiv 1$, we simply write $L^q_vL^{p}_x$ and $\W^{\sigma,q}_v\W^{s,p}_x$. 
 
\smallskip
\noindent
For any $z \in \R^{3} \setminus \{0\}$, we denote by $\widehat{z}=\frac{z}{|z|}$ the associated unit vector.
For two tensors $A=(A_{i,j}), B=(B_{i,j}) \in \mathscr{M}_{d}(\R)$, we denote by $A:B$ the scalar product $(A,B)=\sum_{i,j}A_{i,j}B_{i,j} \in \R$ as the trace of the matrix product $AB$.  For a vector function $w=w(x) \in \R^{3}$, the tensor $(\partial_{x_{i}}w_{j})_{i,j}$ is denoted as $\nabla_{\!x} w$. We  also write 
$(\mathrm{Div}_{x}A)^{i}=\sum_{j}\partial_{x_{j}}A_{i,j}(x).$

\smallskip
Throughout the paper, for $\bm{A,B}> 0$ we indicate $\bm{A} \lesssim \bm{B}$ 
whenever there is a positive constant $C>0$ depending only on fixed quantities but never on the asymptotic parameters or variables $\ell>0$, $\e \in (0,1]$ and $t \geq 0$ {such that $\bm{A} \leq C \bm{B}$}.  Such constant $C >0$ could depend on $\re(\cdot)$, for example through $\g$ and $\overline{\g}$.  We stress that the restitution coefficient is fixed.

The particular weight
$$\m_{s}(v):=(1+ |v|^{2})^{\frac{s}{2}}\,, \qquad \forall \,v \in \R^{3}\,,\qquad \forall \, s \geq 0\,,$$
will be of importance along the document. In particular, we will use in the sequel the following functional spaces
\begin{subequations}\label{def:func}
\begin{equation}\label{eq:defE}
\E:=L^1_v\W_x^{m,2}( {\m_q})\,, \qquad m > \frac{3}{2}\,, \quad q \geq 3+\g\,,\end{equation}
where $\g >0$ appears in \eqref{eq:ab00} and define  
\begin{equation}\label{eq:defE1}
\E_1:=L^1_v\W_x^{m,2}( \m_{q+1})\,, \qquad m > \frac{3}{2}\,, \quad q \geq 3+\g\,.\end{equation}
Finally, we will also work in the spaces 
\begin{equation}\label{eq:defH}
\H:=L^2_v\W^{m,2}_x(\M^{-\frac{1}{2}})\,, \qquad m > \frac{3}{2}\,,
\end{equation}
and 
\begin{equation}\label{eq:defH1}
\H_1:=L^2_v\W^{m,2}_x(\langle v \rangle^{\frac12}\M^{-\frac{1}{2}})\,, \qquad m > \frac{3}{2}\,,
\end{equation}
which is the  natural functional space to study the \emph{elastic} Boltzmann equation.
\end{subequations}

\subsection{Main results}\label{sec:main} We present the main results regarding well-possedness, hydrodynamic limit, and time-asymptotic, for solutions $h_\e=h_\e(t,x,v)$ to the Boltzmann equation~\eqref{eq:heps}.
 with initial condition
$h_\e(0,x,v)=h_\e^{\rm in}(x,v)\,.$
We seek solutions $h_{\e} \in \E$ to \eqref{eq:heps} where $\E$ is defined in \eqref{eq:defE}.
Recall from \eqref{ansatz} that $h_\e=\e^{-1}\left(f_\e-\M\right)$ is the fluctuation around $\M$ of  solutions $f_\e$ to the rescaled Boltzmann equation \eqref{BE0-In}.  It was mentioned that the strategy is based on enforcing the elastic Boltzmann operator $\Q_1$ and its linearized counterpart $\LL$ in equation \eqref{eq:heps}. To do so, 
it is convenient to consider a splitting of the solution $h_\e$ to \eqref{eq:heps} in the form
$$
h_{\e}(t)=\ho_{\e}(t)+\hu_{\e}(t)$$where $\ho=\ho_{\e} \in \E$ and $\hu=\hu_{\e} \in \H$ are the solutions of a system of equations which, through their coupling, is equivalent to \eqref{eq:heps}.  The details are presented in Section \ref{sec:nonlinear}. 
\begin{theo}\label{theo:main}
There exists a pair $(\varepsilon_0,{\alpha_{0}})$ of positive constants that depend on the mass and energy of $F_{\mathrm{in}}^\e$, $m$, and $q$ introduced in~\eqref{eq:defE}, such that for $\e\in(0,\e_0)$ and if
$$\|h_{\mathrm{in}}^\e\|_{ \E} \leq {\alpha_{0}},$$
then, the equation \eqref{eq:heps} has a unique solution $h_{\e} \in\mathcal{C}\big([0,\infty); \E\big) \cap L^1_{\rm{loc}}\big((0,\infty);\E_1\big)$ 
satisfying 
\begin{equation}\label{eq:heZg}
\left\|h_\e(t)\right\|_{ \E}\leq \frac{C ({\alpha_{0}})}{\sqrt{1+t}}\,,
\qquad \forall \, t \geq 0\,, \qquad \forall \, \e \in (0,\e_0)\,,
\end{equation}
for some positive constant $C({ \alpha_{0}}) >0$ independent of $\e$. Moreover, the  solution splits as
\begin{equation}\label{eq:heps0+1}
h_{\e}(t)=\ho_{\e}(t)+\hu_{\e}(t)\end{equation}
and, for any $T >0$ there exists $\bm{h}\in L^{2}\left((0,T);\H\right)$ such that, up to extraction of a subsequence, it holds that
	\begin{equation}\begin{cases}\label{eq:mode-conv}
	\left\{\ho_{\e}\right\}_{\e} \text{converges to $0$ strongly  in } L^{1}((0,T)\,;{\E_1})\,, \\[0.2cm]
	\left\{\hu_{\e}\right\}_{\e} \text{converges to $\bm{h}$ weakly in } L^{2}\left((0,T)\,;\H\right)\,.
	\end{cases}\end{equation}
In particular, there exist 
	\begin{align*}
	\varrho \in L^{2}&\left((0,T)\,;\W^{m,2}_{x}(\T^{d})\right)\,, \qquad
	u \in L^{2}\left((0,T)\,;\left(\W^{m,2}_{x}(\T^{d})\right)^{d}\right)\,,  \\
	&\qquad\text{and} \qquad \vE \in L^{2}\left((0,T)\,;\W^{m,2}_x(\T^{d})\right)\,, 
	\end{align*}
such that
	\begin{equation}\label{def:bmh}
	\bm{h}(t,x,v)=\left(\varrho(t,x)+u(t,x)\cdot v + \frac{1}{2}\vE(t,x)(|v|^{2}-3\en_{\star})\right)\M(v)
	\end{equation}
where $\M$ is the Maxwellian distribution \eqref{eq:max} and $(\varrho,u,\vE)=(\varrho(t,x),u(t,x),\vE(t,x))$ are suitable solutions to the following incompressible Navier-Stokes-Fourier system with forcing term
\begin{equation}\label{eq:NSFint}
\begin{cases}
\partial_{t}u-\frac{\nu_{0}}{\en_{\star}}\,\Delta_{x}u + {\en_{\star}}u\cdot \nabla_{x}\,u+\nabla_{x}p= \xi(t)\,u\,,\\[6pt]
\partial_{t}\,\vE-\dfrac{\nu_{1}}{\en_{\star}^{2}}\,\Delta_{x}\vE + \en_{\star}\,u\cdot \nabla_{x}\vE
	=3\frac{1-\g}{2}\,\en_{\star}^{2}\,\xi(t)\,\vE\,,\\[8pt]
\mathrm{div}_{x}u=0\,, \qquad \varrho + \en_{\star}\,\vE = 0\,,
\end{cases}
\end{equation}
subject to initial conditions $(\varrho_{\mathrm{in}},u_{\mathrm{in}},\vE_{\mathrm{in}})$, entirely determined by the limiting behaviour of $h_{\rm in}^{\e}$ as~$†\e \to 0$. The viscosity $\nu_0 >0$ and heat conductivity $\nu_1 >0$ as well as the forcing coefficient $\xi(t)$ are explicit.
 \end{theo}
 \begin{nb} Theorem \ref{theo:main} states peculiar version of a more general result in which the decay is given by $(1+t)^{-1+ \delta}$ for $\delta$ arbitrarily small, see Theorem \ref{theo:main-cauchy1} in Section \ref{sec:hydro}.  
\end{nb}
\begin{nb}  Notice that the source term $S_\e(\cdot)$ defined in~\eqref{def:LetSe} is not belonging to $L^1((0,\infty),\E_1)$ and this explains why we can only have here
$$\int_0^T \|h_\e(t)\|_{\E_1}\, \d t \leq C_T\,, \qquad \forall \, T >0\,.$$
Notice however that the constant $C(\alpha_0)$ in the bound \eqref{eq:heZg} is not depending on such local $L^1$-bound.\end{nb}
The physical relevance of Theorem \ref{theo:main} is reflected in its effectiveness to provide an explicit cooling rate of the granular temperature, namely \emph{Haff's law} both for the global and local temperatures. 
\begin{theo}[\textit{\textbf{Haff's law -- local and global}}]\label{theo:haff}  Consider a solution $f_{\e}(t,x,v)=\M+\e h_\e(t,x,v)$ to \eqref{BE0-In} as constructed in Theorem \ref{theo:main} with $\e_{0}$ and $\alpha_0$ sufficiently small.  {We recall that the global temperature $\bm{T}_{\e}$ associated to $F_\e$ is defined in~\eqref{eq:Fescal}, we also define its local temperature as}
\begin{equation*}
{T}_{\e}(t,x):=\int_{\R^{3}}{F}_{\e}(t,x,v)|v|^{2}\,\d v \qquad \forall \,(x,t) \in \T^{3} \times \R^+\,.
\end{equation*}
Then, for all \emph{fixed}  {$\e \in (0,\e_0)$}, 
$$
\e^{\frac{4}{\g+1}} \left(\e^{-\frac2\g} + t\right)^{-\frac{2}{\g+1}} \lesssim T_\e(t,x) \lesssim \e^{\frac{4}{\g+1}} \left(\e^{-\frac2\g} + t\right)^{-\frac{2}{\g+1}} \qquad \forall \, \, (x,t) \in \T^3 \times \R^+\,.$$
In particular, 
$$
\e^{\frac{4}{\g+1}} \left(\e^{-\frac2\g} + t\right)^{-\frac{2}{\g+1}} \lesssim \bm{T}_\e(t) \lesssim \e^{\frac{4}{\g+1}} \left(\e^{-\frac2\g} + t\right)^{-\frac{2}{\g+1}} \qquad \forall \, t \in \R^+\,.
$$
\end{theo}

 \subsection{Outline of the paper}
Section \ref{sec:collL} states the properties of the collision operator for viscoelastic hard-spheres.  In particular, it shows the scaling properties of the collision operators~$\mathscr{L}_{\re}$ and $\Q_{\re}$ under the rescaling of the restitution coefficient $\re_\ell(\cdot)$ introduced in~\eqref{def:e_ell} 
and obtain new quantitative estimates about the differences $\Q_{\re_\l}-\Q_1$ and $\mathscr{L}_{\re_{\l}}-\LL$ in various norms. These results are important for deriving the sharp scaling function $V_\e(t)$  pertinent to the self-similar change of variable \eqref{eq:Fescal}. The derivation of this scaling function is done in Section \ref{sec:scaling}. The linearized operator associated to \eqref{eq:heps} is studied in Section \ref{sec:LIN} where it is recalled the splitting of $\mathscr{L}_{\re}$ as well as the hypocoercivity results as derived in \cite{ALT}, see also \cite{BCMT,CRT}, which allow the implementation of the energy method as described in the points c)--d) of Section \ref{sec:liter}.  The energy method and the \emph{a priori} estimates are described in Section \ref{sec:nonlinear} where a splitting of equation \eqref{eq:heps} is described.  In Section \ref{sec:hydro}  we explain how to obtain the existence of solution to the Boltzmann equation as well as the derivation of the hydrodynamic limit for the studied inelastic Boltzmann equation. This provides a complete proof of Theorem \ref{theo:main}.  Finally, Section \ref{sec:orig} gives a proof of Theorem \ref{theo:haff} and show the validity of the local and global Haff's law for viscoelastic hard spheres.  The paper ends with several appendices, namely, Appendix~\ref{app:hydro} which recalls tools useful for the derivation of hydrodynamic limit, see also \cite{ALT}, Appendix \ref{app:collL} which recalls important estimates for the collision operator $\Q_\re$ in the case of viscoelastic hard spheres. Finally, Appendix \ref{app:vis} shows that the true viscoelastic particles, as derived in \cite{PoSc}, meets the assumptions introduced in Section~\ref{Sec21} providing the main application to our results.

\subsection*{Acknowledgments.} B. L. gratefully acknowledges the financial support from the Italian Ministry of Education, University and Research (MIUR), Dipartimenti di Eccellenza grant 2022-2027, as well as the support from the de Castro Statistics Initiative, Collegio Carlo Alberto (Torino). I. T. was supported by the French government through the France 2030 investment plan managed by the ANR, as part of the Initiative of Excellence Universit\'e  C\^ote d'Azur under reference number ANR-15-IDEX-01 as well as by the MaDynOS ANR-24-CE40-3535-01.

\section{Summary of useful results about the collision operator}\label{sec:collL}

 \subsection{Collision mechanism in granular gases}\label{Sec21} The Boltzmann equation for granular gases is a well-accepted model that describes a system composed by a large number of granular particles which are assumed to be hard-spheres with equal mass undertaking inelastic collisions.  The collision mechanism is characterized by the \textit{coefficient of normal restitution} denoted by $\re \in(0,1]$.  For viscoelastic particles the coefficient of normal restitution depends solely on the impact velocity, that is, if $v$ and $\vb$  denote the velocities of two particles before collision, their respective velocities~$v'$ and~$\vb'$ after collision satisfy \eqref{coef} and \eqref{transfpre}. For the coefficient of normal restitution, we adopt the following definition, see \cite{A}.
\begin{defi}\label{defiC} A coefficient of normal restitution $\re \: :\: r \mapsto \re(r) \in (0,1]$ belongs to the class~$\mathcal{R}_0$ if it satisfies the following:
\begin{enumerate}
\item The mapping  $r \in \mathbb{R}^{+} \longmapsto \re(r) \in (0,1]$ is  {differentiable} and non-increasing.
\item The mapping $r\in\mathbb{R}^{+}\longmapsto \eta_{\re}(r):=r\,\re(r )$ is strictly increasing.
\item $\lim_{r \to \infty} \re(r)=\re_0 \in [0,1)$.
\end{enumerate}
Moreover, for a given $\gamma > 0$ we say that $\re(\cdot)$ belongs to the class $\mathcal{R}_\gamma$ if it belongs to $\mathcal{R}_0$ and there exist positive $\mathfrak{a}_0$, $\mathfrak{b}_0$, and $\overline{\g}>\frac{3}{2}\g$ such that for any $r\geq0$, 
\begin{equation}\label{gamma0bis}
|\re(r)-1+\mathfrak{a}_{0}r^{\g}| \leq \mathfrak{b}_{0} r^{\overline{\g}}\,.\end{equation}
In addition, we assume that there exist $\alpha >0$ and $m \in (0,1)$ such that for any $r \geq 0$, 
\begin{equation}\label{eq:Je-e}
\frac{\d}{\d r} \eta_{\re}(r) \gtrsim \re(r)^\alpha\,,  \qquad  \qquad \left|{\frac{\d}{\d r} \eta_{\re}(r)}-\re(r)\right| \lesssim |1-\re(r)|\,,
\end{equation}
and 
\begin{equation}\label{eq:Je-e1}
\limsup_{z\to\infty}z^{-m}\eta_{\re}^{-1}(z)\lesssim 1\,. \end{equation}
\end{defi} 

\begin{nb} Assumptions \eqref{eq:Je-e} and \eqref{eq:Je-e1} are most likely of technical nature.  Thus, we only emphasize the role of $\g$ in the definition of $\mathcal{R}_{\g}$.  We highlight below several consequences of the assumptions made on $\re(\cdot)$. 
\begin{itemize}
\item[--] Using point  $\textit{(2)}$ in Definition \ref{defiC}, the Jacobian of the transformation \eqref{transfpre} is given by
\begin{equation} \label{def:Je} 
J_{\re}:=J_{\re}\big(|u\cdot \n|\big)=\left|\dfrac{\partial (v',\vb')}{\partial (v,\vb)}\right|
=\re\big(|u\cdot \n|\big) + |u\cdot \n|\dfrac{\d \re}{\d r}(|u\cdot \n|)=\frac{\d}{\d r}\eta_{\re}(|u\cdot \n|)  >0\,. 
\end{equation}
In particular, the first part of the technical assumption \eqref{eq:Je-e} implies that 
\begin{equation} \label{eq:boundJac}
J_{\re}\big(|u\cdot \n|\big) \gtrsim\re\big(|u\cdot \n|\big)^{\alpha}\,, \qquad \forall \, u \in \R^{3}\,,\,\,n \in \S^{2}\,.
\end{equation}
\item[--]
 The inequality~\eqref{eq:Je-e1} implies that for  any $r \geq 0$, 
\begin{equation}\label{eq:infre}
\re(r) \gtrsim (1+r)^{1-\frac1m}\,.
\end{equation}
Indeed, using that $\theta_\re^{-1}$ is strictly increasing, it implies that 
$$
\limsup_{r \to \infty} r^{-m+1} \re(r)^{-m} \lesssim 1
\quad \text{i.e.} \quad 
\liminf_{r \to \infty} r^{1-m} \re(r)^m \gtrsim 1
$$
and this gives \eqref{eq:infre} recalling that $\re(0)=1$.
\item[--] For the rescaled restitution coefficient $\re_\ell(\cdot) := \re (\ell \cdot)$, for $\ell>0$, it follows that 
$$
\eta_{\re_{\l}}(r)=\frac{1}{\l}\eta_{\re}(\l r)
\quad \text{so that} \quad 
\eta_{\re_{\l}}^{-1}(z)=\frac{1}{\l}\eta_{\re}^{-1}(\l z)\,.
$$
Thus, \eqref{eq:Je-e1} implies that for $r$ sufficiently large
\begin{equation}\label{eq:Je-e2}
\eta_{\re_{\l}}^{-1}(r) \lesssim \l^{m-1}r^{m}\,, \quad \forall \, \ell \in (0,1]\,.
\end{equation}
\end{itemize}
\end{nb}
\begin{exe} \label{nb:visco} An important model is the viscoelastic hard spheres.  The coefficient of normal restitution is given by the expansion
\begin{equation}\label{visc}
\re(r)=1+ \sum_{k=1}^\infty (-1)^k \, a_k \, r^{\frac{k}{5}}\,, \qquad \forall \,r > 0\,,
\end{equation}
where $a_k \geq 0$ for any $k \in  \mathbb{N}$.  We refer to \cite{Poschel,PoSc} for the physical considerations leading to such expansion.  We show in  {Appendix \ref{app:vis}} that $\re(\cdot)$ belong to the class $\mathcal{R}_{\g}$ with $\g=\frac{1}{5}$, $\overline{\g}=2\g$, and $\re_{0}=0$.  \end{exe}
 
\noindent
Due to the previous fundamental example, any model with coefficient of normal restitution belonging to $\mathcal{R}_{\gamma}$ is referred to as \textit{generalized viscoelastic particles} model.  We now fix the restitution coefficient $\re(\cdot)\in\mathcal{R}_\gamma$.  In particular, $\gamma>0$ is fixed as well.

\subsection{Strong form and alternative representation of the Boltzmann operator} \label{Sec22}
We provide a strong form of the Boltzmann collision operator defined in weak form in \eqref{eq:weakN}.
Pre-collisional velocities $('v,'\vb)$, resulting in $(v,\vb)$ after collision, can be introduced through the relation
\begin{equation}\label{'v'vb}
v=\,'v-\frac{1+\re\big( | 'u \cdot n| \big)}{2} \big(\, 'u \cdot n \big) n\,, \quad \vb=\,'\vb+\frac{1+\re\big( |'u \cdot n| \big)}{2}\big(\, 'u \cdot n \big) n, \quad 'u:=\,'v-'\vb\,.
\end{equation} In particular, the energy relation in the collision mechanism can be written as
\begin{equation}\label{energ}
|v|^{2}+| \vb|^{2}=|'v|^{2}+|'\vb|^{2}-\frac{1-\re^{2}\big(|'u \cdot n|\big)}{2}\,\big(\,'u \cdot n\big)^{2}, \qquad
u\cdot n=-\re\big(|'u \cdot n|\big)\big(\,'u \cdot n\big).
\end{equation}
For a pair of distributions $f=f(v)$ and $g=g(v)$, the Boltzmann collision operator can be defined in strong form as the difference of two nonnegative operators, gain and loss operators respectively
\begin{equation*}
\Q_{\re}\big(f,g\big)=\Q_{\re}^{+}\big(f,g\big)-\Q_{\re}^{-}\big(f,g\big),
\end{equation*}
with
\begin{multline}\label{Boltstrong}
\Q_{\re}^{+}\big(f,g\big)(v)=\It \dfrac{B_0(u,\n)}{\re\big(|'u\cdot n|\big)J_{\re}\big(|'u\cdot n|\big)}f('v)g('\vb)\,\d\vb\,\d\n\,,\\
\text{and}\qquad\Q_{\re}^{-}\big(f,g\big)(v)=f(v)\It  B_0(u,\n)g(\vb)\,\d\vb\,\d\n\,,
\end{multline}
where $J_{\re}$ is defined in \eqref{def:Je} and we recall that the collision kernel $B_0(u,\n)$ is defined in~\eqref{def:B0} is of the form $B_0(u,\n)=|u|\,b_0\big(\widehat{u} \cdot \n\big)$
where $b_0(\cdot)$, the scattering kernel, is a nonnegative function. The main technical assumption on $b_{0}(\cdot)$ is derived from the following equivalent form of the Boltzmann collision operator based on the so-called $\sigma$-parametrization or post-collisional velocity parametrization.  More precisely, let $v$ and $\vb$ be the collisional particles' velocities, recall that~$u=v-v_{\ast}$ and $\widehat{u} = u/|u|$, performing in~\eqref{transfpre} the change of unknown
\begin{equation*}
\sigma=\widehat{u}-2 \,(\widehat{u}\cdot \n)\n \in \mathbb{S}^2 
\end{equation*}
provides an alternative parametrization of the unit sphere $\mathbb{S}^2$.  In this case, the impact velocity becomes
\begin{equation*}
|u\cdot\n|=|u| \,|\widehat{u} \cdot \n|=|u| \sqrt{\frac{1-\widehat{u} \cdot \sigma}{2}}\,.
\end{equation*}
Therefore,  the post-collisional velocities $(v',\vb')$ given in~\eqref{transfpre} can be written as
\begin{multline}\label{co:transf}
v'=v+\frac{1}{4}\left(1+\re\left(|u|\sqrt{\frac{1-\widehat{u}\cdot \sigma}{2}}\right)\right)\,(|u|\sigma-u)\,,\\
\text{and}\qquad v_{\ast}'=v_{\ast}-\frac{1}{4}\left(1+\re\left(|u|\sqrt{\frac{1-\widehat{u}\cdot \sigma}{2}}\right)\right)\,(|u|\sigma-u)\,.
\end{multline}
Relations \eqref{co:transf} lead to the weak formulation of the collision operator in the $\sigma$-representation
\begin{equation}\label{eq:weak}\begin{split}
 \int_{\R^{3}} \Q_{\re} (g,f)(v)\, \psi(v)\, \d v  
 &= \int_{\R^3 \times \R^3} f(v) g(v_\ast) |v-v_*| \int_{\S^2} (\psi(v') - \psi(v)) \, \d\sigma \, \d v_\ast \, \d v \\
 &=  \frac{1}{2}  \int_{\R^{3} \times \R^{3}} f(v)\,g(v_{\ast})\,|v-v_{\ast}|
\mathcal{A}_{\re}[\psi](v,v_{\ast})\, \d v_{\ast}\, \d v\,,
\end{split}
\end{equation}
where
\begin{equation}    \label{eq:psi} \mathcal{A}_{\re}[\psi](v,v_{\ast}) :=
    \int_{\S^{2}}(\psi(v')+\psi(v_{\ast}')-\psi(v)-\psi(v_{\ast}))b(\sigma \cdot \widehat{u})\, \d{\sigma}\,.
\end{equation}
Here, $\d\sigma$ denotes the Lebesgue measure on $\S^{2}$.  In the case of granular gases there is no loss of generality in assuming
\begin{equation}\label{unitb}
\int_{\S^{2}}b(\sigma \cdot \widehat{z})\, \d{\sigma}=1\,, \qquad \forall \, \widehat{z} \in \S^{2}\,.
\end{equation}
The mapping $b(\cdot)$ defined on $(-1,1)$ is identified with a function defined on the sphere $\S^{2}$.  Besides \eqref{unitb}, we additionally assume that 
\begin{equation}\label{eq:bW11}
b(\cdot) \in \W^{1,1}((-1,1))\end{equation}
and, due to the Sobolev embedding $\W^{1,1}((-1,1)) \hookrightarrow L^{\infty}((-1,1))$, we notice that $b(\cdot)$ is actually bounded.  Moreover, the angular kernels $b_0$ and $b$ respectively introduced in \eqref{def:B0} and \eqref{eq:psi} are related through 
\begin{equation}\label{eq:bb0}
	b(\widehat u \cdot \sigma) = |\widehat u \cdot n|^{-1} b_0(\widehat u \cdot n) \,. 
\end{equation}
For any fixed vector $\widehat{z}$, the angular kernel $b_0$ defines a measure on the sphere through the mapping $ \n \in \mathbb{S}^{2}\mapsto b_0\big( \widehat{z}\cdot n \big)\in[0,\infty]$.  Assumption \eqref{unitb} translates into
\begin{equation}\label{normalization}
\left\| \big(\widehat{z}\cdot n \big)^{-1}\;b_{0}(\widehat{z}\cdot n) \right\|_{L^{1}(\mathbb{S}^{2},\,\d \n)}=2\pi \left\| s^{-1}\,b_{0}(s) \right \|_{L^{1}((-1,1),\,\d s)}=1\,,
\end{equation}
which is exactly the assumption we made on $b_0$ in~\eqref{eq:cutoff}.
The identification \eqref{eq:bb0} implies that the mapping $s \in (-1,1) \mapsto s^{-1}b_{0}(s)$ belongs to $\W^{1,1}((-1,1))$ and, as such, is bounded, and so is~$b_{0}(\cdot)$. In particular, 
\begin{equation}\label{eq:bL2}
\int_{-1}^{1}s^{-2}b_{0}^{2}(s)\,\d s < \infty\,.
\end{equation}
Such assumption is required to obtain some estimates for $\Q_{\re}$ in $L^{2}$-spaces with Gaussian weights, see the proof of Lemma \ref{lem:QEMM}.  In this document, we consider the anisotropic hard spheres model for which the kinetic potential is equal to $|\cdot|$ in \eqref{def:B0} with angular kernel $b_0$ satisfying \eqref{normalization}-\eqref{eq:bL2}.  The hard spheres model is a particular case of this setting for which
\begin{equation*}
b_0(\widehat{u} \cdot \n)=\tfrac{1}{4\pi}|\widehat{u}\cdot \n|\,.
\end{equation*}
\subsection{Estimates on the difference between inelastic and elastic collisions}
We recall that, for~$\ell>0$, the rescaled restitution coefficient $\re_\ell(\cdot)$ is defined in~\eqref{def:e_ell}. 
Our main technical contribution in this section is a quantitative study of the differences~$\Q_{\re_{\ell}}-\Q_{1}$ and the linearized counterpart~$\mathscr{L}_{\re_{\ell}}-\LL$.  We significantly improve previous estimates obtained in \cite{ACG,MiMo3}   and extend them to the case of non constant restitution coefficient. 
We recall that from~ \eqref{co:transf}--\eqref{eq:weak}, we have
\begin{equation*}
\int_{\mathbb{R}^3}\Q^{+}_{e_\l}(g,f)(v) \psi(v) \, \d v= \int_{\mathbb{R}^3}\int_{\mathbb{R}^3}f(v)g(v - u)|u|\bigg(\int_{\mathbb{S}^{2}}\psi(v'_\l)b(\widehat{u}\cdot\sigma)\,\d\sigma\bigg) \, \d u \, \d v
\end{equation*} 
where
\begin{equation*}
v'_\l = v - \frac{\beta_\l}{2}\,\big( u - |u|\sigma \big)\quad\text{with}\quad \beta_\l:=\beta\Bigg( \l\,|u| \sqrt{\frac{1-\widehat{u}\cdot\sigma}{2}}\Bigg) := \frac{1}{2}\Bigg( 1 + \re\Bigg(\l \, |u| \sqrt{\frac{1-\widehat{u}\cdot\sigma}{2}}\Bigg)\Bigg)\,,
\end{equation*} 
and we have used that the impact velocity writes $|\widehat u \cdot n| = |u| \sqrt{\frac{1-\widehat u \cdot \sigma}{2}}$. 
\begin{prop}\label{prop:Qel}  
For any $q\geq0$, any $\l>0$, and suitable $f$ and $g$, we have that
\begin{equation*}
\big\| \Q^{+}_{\re_\l}(f,g) - \Q^{+}_1(f,g)  \|_{L^{2}_{v}(\m_q)} 
\lesssim \l^{\gamma} \, \| f\|_{L^{1}_{v}(\m_{q+3+\gamma})}\|g\|_{\W^{1,2}_{v}(\m_{q+3+\gamma})}
\end{equation*}
and
\begin{equation*}
\big\| \Q^{+}_{\re_\l}(g,f) - \Q^{+}_1(g,f)  \|_{L^{2}_{v}(\m_q)} 
\lesssim \l^{\gamma} \, \| f\|_{L^{2}_{v}(\m_{q+3+\gamma})}\|g\|_{\W^{1,2}_{v}(\m_{q+3+\gamma})}\,.
\end{equation*}
\end{prop}
\begin{proof}
Fix $\delta>0$ and express $b = b_{\delta} + b^{c}_{\delta}$ with $b_{\delta}$ supported in $[-1+\delta, 1 - \delta]$ and with $\| b'_{\delta} \|_{1} \lesssim  {\| b \|_{\infty}}+\| b' \|_{1}$ and $\| b^{c}_{\delta} \|_{1} {\lesssim \delta \|b\|_{\infty}}$.    Define ${B}_{\delta} := |u| \, b_{\delta}$ and ${B}^{c}_{\delta} := |u| \, b^{c}_{\delta}$ and consider $\psi \in L^2_v$ a test function, then using Proposition \ref{propB:CMS} it follows that
\begin{multline*}
\int_{\mathbb{R}^{3}}\Big( \Q^{+}_{{B}_{\delta},\re_\l}(f,g) - \Q^{+}_{{B}_{\delta},1}(f,g) \Big) \psi\, \d v \lesssim \l^{\gamma}\int_{\mathbb{R}^{3}} \Q^{+}_{{B}_{ {\delta},\gamma},1}(f,g)\, |\psi| \, \d v \\+ \l^{\gamma}\int^{1}_{0}\int_{\mathbb{R}^{3}} \Q^{+}_{\bar{B}_{ {\delta},\gamma},\tilde{\re}_{s}}(f,h)\, |\psi| \, \d v \,\d s
\end{multline*}
where  $B_{{\delta},\gamma} := B_{\delta} |\cdot|^{\gamma}$ and $\bar{B}_{{\delta},\gamma} := \max\{B_{\delta}, |\nabla_u B_{\delta}| \} |\cdot|^{2+\gamma}$.  Moreover, $\tilde{\re}_{s}$ belongs to the class~{$\mathcal{R}_0$} for $s\in [0,1]$ and $h= g+|\nabla g|$.  By the convolution inequalities in \cite[Theorem~1]{ACG}
\begin{equation}\label{Bgamma1}
\begin{aligned}
\int_{\mathbb{R}^{3}} \Q^{+}_{B_{{\delta}, \gamma},1}(f,g)\, |\psi| \, \d v 
\lesssim \|b_{{\delta}}\|_{1}\|f\|_{ L^{1}_{v}(\m_{1+\gamma}) } \| g \|_{ L^{2}_{v}(\m_{1+\gamma})} \|\psi\|_{L^{2}_{v}} \\
\lesssim \|b\|_{1}\|f\|_{ L^{1}_{v}(\m_{1+\gamma}) } \| g \|_{ L^{2}_{v}(\m_{1+\gamma})} \|\psi\|_{L^{2}_{v}}\,.
\end{aligned}
\end{equation}
Note that $|\nabla_u B_{\delta}| \leq b_{\delta}+  {|u|} |b'_{\delta}| $, therefore, using again \cite[Theorem 1]{ACG} 
\begin{multline}\label{Bgamma2}
\int^{1}_{0}\d s\int_{\mathbb{R}^{3}} \Q^{+}_{\bar{B}_{\gamma, {\delta}},\tilde{\re}_{s}}(f,h)\, |\psi| \, \d v\,\ds 
\\
\lesssim \big( {\|b\|_{\infty}} + \|b\|_1+ \|b'\|_{1}\big) \|f\|_{ L^{1}_{v}(\m_{3+\gamma}) } \| h \|_{ L^{2}_{v}(\m_{3+\gamma})} \|\psi\|_{L^{2}_{v}}\,.
\end{multline}
Meanwhile, for the portion corresponding to $b^{c}_{\delta}$ and the associated kernel $B^{c}_{\delta}$, one has that
\begin{multline}\label{Bcomp}
\int_{\mathbb{R}^{3}}\Big( \Q^{+}_{B^{c}_{\delta},\re_\l}(f,g) - \Q^{+}_{B^{c}_{\delta},1}(f,g) \Big) \psi\, \d v \lesssim \|b^{c}_{\delta}\|_{1}\|f\|_{L^{1}_{v}(\m_{1})} \| g \|_{ L^{2}_{v}(\m_{1})} \|\psi\|_{L^{2}_{v}}  \\
\lesssim \delta \, {\|b\|_\infty}  \|f\|_{ L^{1}_{v}(\m_{1}) } \| g \|_{ L^{2}_{v}(\m_{1})} \|\psi\|_{L^{2}_{v}}\,.
\end{multline}
Choosing $\delta  {= \ell^\gamma}$ and using that $\|h\|_{L^{2}_v(\m_{3+\gamma})} \lesssim \|g\|_{\W^{1,2}_v(\m_{3+\gamma})}$, it follows, by gathering the estimates \eqref{Bgamma1}, \eqref{Bgamma2} and \eqref{Bcomp}, that
\begin{equation*}
\int_{\mathbb{R}^{3}}\Big( \Q^{+}_{B,\re_\l}(f,g) - \Q^{+}_{B,1}(f,g) \Big) \psi\, \d v  \lesssim \l^{\gamma} \|f\|_{ L^{1}_{v}(\m_{3+\gamma})} \| g \|_{ \W^{1,2}_{v}(\m_{3+\gamma})} \|\psi\|_{L^{2}_{v}}\,,
\end{equation*}
since by assumption the quantities $\|b\|_\infty$, $\|b\|_1$ and $\|b'\|_1$ are finite from~\eqref{eq:bW11}. This proves the first inequality without weights, that is the case $q=0$.  Adding weights, the case $q>0$, follows in standard fashion using the test function $\tilde{\psi} = \langle v \rangle^{q}\psi$ and distributing the weight in the functions $f$ and $g$ by invoking the dissipation of energy $|v'_{\l}|^{2}\leq |v|^{2} + |v_{*}|^{2}$. 
The treatment of the second inequality is similar and left to the reader.
\end{proof}
Recalling that $\M$ has been defined in \eqref{eq:max}, the above result allows to estimate  the difference between the linear operators $\mathscr{L}_{\re_{\l}}$ and $\LL$ defined as 
\begin{equation} \label{def:LL}
\mathscr{L}_{\re_{\l}}(h):=\Q_{\re_{\l}}(h,\M)+\Q_{\re_{\l}}(\M,h)
\qquad \text{and} \qquad 
\LL h := \Q_1(h,\M) + \Q_1(\M,h)\,.
\end{equation}
\begin{cor}\label{lem:LeL1} 
For any $q \geq 0$, $\kappa >\frac{3}{2}$, and $\ell>0$, we have that
\begin{equation}\label{eq:LeL1kappa}
\left\|\mathscr{L}_{\re_{\l}}(h)-\LL(h)\right\|_{L^{1}_{v}(\m_{q})} 
\lesssim \l^{\g}\left\|h\right\|_{L^{2}_{v}(\m_{q+2\kappa+3+\g})}\,.
\end{equation}
\end{cor} 
\begin{proof} 
Notice that from the Cauchy-Schwarz inequality and Proposition \ref{prop:Qel}, 
\begin{align*}
\|\mathscr{L}_{\re_{\l}}(h)-\LL(h)\|_{L^{1}_{v}(\m_{q})} &\lesssim  \|\mathscr{L}_{\re_{\l}}(h)-\LL(h)\|_{L^{2}_{v}(\m_{q+\kappa})} \\
&\lesssim \l^{\g}\left(\|h\|_{L^{2}_{v}(\m_{q+\kappa+3+\gamma})}+\|h\|_{L^{1}_{v}(\m_{q+\kappa+3+\gamma})}\right)
\end{align*}
since $\kappa >\frac{3}{2}$. One deduces the final result using the Cauchy-Schwarz inequality. 
\end{proof}
The following result, established in \cite[Theorem 3.11]{CMS} for exponential weights is easily modified for polynomial weights.
\begin{prop}\label{prop:conti} There exists an explicit constant $\l_{\star} \in (0,1)$ such that for any $q \geq0$ and any~$\ell \in (0,\ell_\star)$,
\begin{equation}\label{L1diff}
\| \Q^+_{\re_{\l}}(f,g)-\Q^+_1(f,g) \|_{L^1_v(\m_{q})}
\lesssim  \l^{\frac{\gamma}{8+3\gamma}}\|f\|_{L^{1}_v(\m_{q+1})}\,\|g\|_{\W^{1,1}_v(\m_{q+1})}\,,
\end{equation}
and
\begin{equation*}
\| \Q^+_{\re_{\l}}(f,g)-\Q^+_1(f,g) \|_{L^1_v(\m_{q})}
\lesssim  \l^{\frac{\gamma}{8+3\gamma}}\|g\|_{L^{1}_v(\m_{q+1})}\,\|f\|_{\W^{1,1}_v(\m_{q+1})}\,.
\end{equation*}
\end{prop}
A consequence of this proposition is the following corollary which can be seen as an estimate on the difference between the linearized operator 
$\mathscr{L}_{\re_{\l}}$ and $\LL$ in the graph norm of $L^{1}_{v}(\m_{q})$.
\begin{cor}\label{cor:Line} 
For any $q \geq 0$ and any $\ell \in (0,\ell_\star)$,
$$
\left\|\mathscr{L}_{\re_{\l}}(h)-\LL(h)\right\|_{L^{1}_v(\m_{q})}\lesssim \l^{\frac{\g}{8+3\g}}\|h\|_{L^{1}_{v}(\m_{q+1})}\,.
$$
\end{cor}
We now establish a general property on the term $\Q_{\re_{\l}}(\M,\M)$ in the subsequent two lemmas. Notice that the entries are Maxwellian distributions, then the term $\Q_{\re_{\l}}(\M,\M)$ is actually equal to~$\Q_{\re_{\l}}(\M,\M)-\Q_{1}(\M,\M)$. We show that it scales with the optimal rate $\l^{\g}$ even in $L^{2}$-spaces with Gaussian weights. Observe that the choice $h=\M$ in Corollary \ref{lem:LeL1} gives already an optimal scaling in all $L^{1}_v(\m_{q})$ spaces. 
\begin{lem}\label{lem:QEMM}
For any $\l \in (0,1]$ we have that
$$
\|\Q_{\re_{\l}}(\M,\M)\|_{L^{2}_{v}(\M^{-\frac{1}{2}})} \lesssim \l^{\g}\,.
$$
\end{lem}
\begin{proof} Using the strong form \eqref{Boltstrong} of the collision operator leads to 
$$
\Q_{\re_{\l}}\big(\M,\M\big)(v)=\It \left(\dfrac{\M('v)\M('\vb)}{\re_{\l}\big(|'u\cdot n|\big)J_{\re_{\l}}\big(|'u\cdot n|\big)}-\M(v)\M(\vb)\right)B_0(u,\n)\,\d\vb\,\d\n\,,
$$
where $'v$ and $'\vb$ corresponds to the pre-collisional velocities.  Also, recall that $\re_{\l}(r)=\re(\l r)$. Using \eqref{energ} gives that
$$
\Q_{\re_{\l}}\big(\M,\M\big)(v)=\It  \M(v)\M(\vb) \mathcal{G}_{\l}(v,\vb,n)B_0(u,\n)\,\d\vb\,\d\n\,,$$
where
$$\mathcal{G}_{\l}(v,\vb,n):=\dfrac{1}{\re_{\l}\big(|'u\cdot n|\big)J_{\re_{\l}}\big(|'u\cdot n|\big)}\exp\left(-\frac{1-\re_{\l}^{2}(|'\!u\cdot n|)}{4\en_{\star}}\left('u\!\cdot n\right)^{2}\right)-1\,.$$
We split $\mathcal{G}_{\l}(v,\vb,n)=\mathcal{G}_{\l}^{(1)}(v,\vb,n)-\mathcal{G}_{\l}^{(2)}(v,\vb,n)$ with
$$\mathcal{G}_{\l}^{(1)}(v,\vb,n):=\exp\left(-\frac{1-\re_{\l}^{2}(|'\!u\cdot n|)}{4\en_{\star}}\left('u\!\cdot n\right)^{2}\right)\left(\dfrac{1}{\re_{\l}\big(|'u\cdot n|\big)J_{\re_{\l}}\big(|'u\cdot n|\big)}-1\right)\,,$$
and
$$\mathcal{G}_{\l}^{(2)}(v,\vb,n):=\exp\left(-\frac{1-\re_{\l}^{2}(|'\!u\cdot n|)}{4\en_{\star}}\left('u\!\cdot n\right)^{2}\right)-1\,.$$
From \eqref{eq:boundJac} we have that
$$\left|\frac{1}{\re_{\l}(|'u\cdot\n|)J_{\re_{\l}}(|'u\cdot n|)}-1\right|  {\lesssim}\, \re_{\l}(|'\!u\cdot n|)^{-1-\alpha}\left|1-\re_{\l}(|'u\cdot n|)J_{\re_{\l}}(|'u\cdot n|)\right|\,.$$
Recall that from \eqref{eq:infre}, we have that $\re_\ell(z) \gtrsim (1+\ell z)^{1-1/m}$, so that if $\l \leq 1$,
$$\left|\frac{1}{\re_{\l}(|'u\cdot\n|)J_{\re_{\l}}(|'u\cdot n|)}-1\right| \lesssim \left(1+|'u\cdot n|\right)^{(1+\alpha){\big(\frac{1}{m}-1\big)}}\left|1-\re(\l|'u\cdot n|)J_{\re}(\l|'u\cdot n|)\right|$$
where the fact that $J_{\re_{\l}}(z)=J_{\re}(\l z)$ was used.   Now, we deduce from  \eqref{eq:Je-e} and \eqref{def:Je} that
$$|1-\re(\l z)J_{\re}(\l z)| \leq |1-\re^{2}(\l z)| + \re(\l z)| {J_{\re}}(\l z)-\re(\l z)| \lesssim |1-\re(\l z)|$$
which together with \eqref{gamma0bis} yields
\begin{align*}
&\left|\mathcal{G}_{\l}^{(1)}(v,\vb,\n)\right| \\
&\qquad \lesssim
\big(\l^{\g}|'\!u\cdot n|^{\g} +\l^{\overline \g}|'u\cdot n|^{\overline \g} \big) (1+|'u\cdot n|)^{ {(1+\alpha)\big(\frac1m-1\big)}}\exp\left(-\frac{1-\re_{\l}^{2}(|'\!u\cdot n|)}{4\en_{\star}}\left('u\!\cdot n\right)^{2}\right)\,.
\end{align*}
The estimate for $\mathcal{G}_{\l}^{(2)}(v,\vb,\n)$ is simpler.  Using the elementary inequality $1-e^{-z} \leq z$ for all~$z\geq0$ and using \eqref{gamma0bis} it holds that
$$|\mathcal{G}_{\l}^{(2)}(v,\vb,n)| \lesssim ('u\cdot n)^{2}\left(1-\re^{2}_{\l}(|'u\cdot n|)\right)\lesssim \l^{\g}|'u\cdot n|^{\g+2}+\l^{\overline{\g}}|'u\cdot n|^{2+\overline{\g}}\,.$$
This proves that there exist $p$ and $\overline{p} >0$ such that
$$|\mathcal{G}_{\l}(v,\vb,\n)| \lesssim \left(\l^{\g}|'u\cdot n|^{p}+\l^{\overline{\g}}|'u\cdot n|^{\overline{p}}\right).$$
Thanks to Cauchy-Schwarz inequality it follows that
\begin{multline*}
\|\Q_{\re_{\l}}(\M,\M)\|_{L^{2}_v(\M^{-\frac{1}{2}})}^{2}=\int_{\R^{3}}\M(v)\left(\int_{\R^{3}\times \S^{2}}\M(\vb)\mathcal{G}_{\l}(v,\vb,n)B_{0}(u,\n)\,\d n\,\d\vb\right)^{2}\d v\\
\leq \int_{\R^{3}}\M(v)\left(\int_{\R^{3}\times \S^{2}}\M(\vet)B_{0}(u,n)^{2}\,\d n\,\d\vet\right)\left(\int_{\R^{3}\times \S^{2}}\M(\vb)|\mathcal{G}_{\l}(v,\vb,\n)|^{2}\,\d n\,\d \vet\right)\d v\,.
\end{multline*}
Since $b_{0}$ satisfies \eqref{eq:bL2}, it holds that
$$ \int_{\R^{3}\times \S^{2}}\M(\vet)B_{0}(u,n)^{2}\,\d n\,\d\vet  \lesssim \langle v\rangle^{2}\,,$$
and thus,
\begin{multline*}\|\Q_{\re_{\l}}(\M,\M)\|_{L^{2}_v(\M^{-\frac{1}{2}})}^{2} \lesssim \int_{\R^{3}\times\R^{3}\times \S^{2}}\langle v\rangle^{2}\M(v)\M(\vb)|\mathcal{G}_{\l}(v,\vb,\n)|^{2}\,\d n \,\d \vet\,\d v\\
\lesssim \l^{2\g}\int_{\R^{3}\times\R^{3}\times \S^{2}}\langle v\rangle^{2}\M(v)\M(\vb)\left(|'u\cdot n|^{2p}+\l^{2(\overline{\g}-\g)}|'u\cdot n|^{2\overline{p}}\right)\,\d\n\, \d \vet \, \d v\,.
\end{multline*}
Recalling from \eqref{'v'vb}, $|'u\cdot n|= \eta_{\re_{\l}}^{-1}(|u\cdot n|)$.  And using \eqref{eq:Je-e2}, we deduce that
\begin{multline*}
\|\Q_{\re_{\l}}(\M,\M)\|_{L^{2}(\M^{-\frac{1}{2}})}^{2} \\
\lesssim \l^{2\g}\int_{\R^{3}\times\R^{3} }\langle v\rangle^{2}\M(v)\M(\vb)
\left(1+\l^{2p+m-1}|u|^{2p+m}+\l^{2(\overline{\g}-\g+\overline{p})+m-1}|u|^{2\overline{p}+m}\right)\,\d \vet\,\d v\,, \end{multline*}
which gives the result since $2p+m-1 >0$ and $2\overline{p}+m-1 >0$  because $m \in (0,1)$ and~$\l \leq 1$.\end{proof} 
Let us introduce the Maxwellian with normalized mass, vanishing momentum, and unit temperature 
\begin{equation} \label{def:M}
M(v) := (2\pi)^{-\frac32} \exp \left(- \frac{|v|^2}{2} \right)\,, \quad \forall \, v \in \R^3\,,
\end{equation} 
and  the $3+s$-order moments of $M$ for $s>0$
\begin{equation} \label{def:Kgamma}
K_s := \int_{\R^3} M(v) |v|^{3+s} \, \d v\,.
\end{equation}
We define the positive constant $\mathfrak{a}_1$ as
\begin{equation} \label{def:a1}
\mathfrak{a}_1 :=  \frac{\mathfrak{a}_{0}K_{\g}}{3(4+\g)}(2\en_{\star})^{\frac{1+\g}{2}}\,,
\end{equation}
where $\mathfrak{a}_0$ has been introduced in \eqref{gamma0bis} and $\en_\star$ is the temperature of $\M$ defined in \eqref{eq:max}.
Considering the flux energy balance between $\Q_{\re_{\l}}(\M,\M)$ and the drift term of order $\mathfrak{a}_1$, the smallness flux rate can be improved. 
\begin{lem}\label{lem:Ener}  Consider $\M$ as defined in \eqref{eq:max} and $\mathfrak{a}_1$ in \eqref{def:a1}.  It holds that
\begin{equation*}
\left|\int_{\R^{3}}|v|^{2}\left[\Q_{\re_{\l}}(\M,\M)-\mathfrak{a}_1\,\l^{\g}\mathrm{div}_v \left(v \M(v)\right)\right]\d v\right|
\lesssim \max\left(\l^{\overline{\gamma}},\l^{2\g},\l^{\overline{\g}+\g}\right)\,. 
\end{equation*}
\end{lem}
\begin{proof}  Set 
$$D_{\l}:=-\int_{\R^{3}}\Q_{\re_{\l}}(\M,\M)\,|v|^{2}\,\d v=\int_{\R^{3}}\int_{\R^{3}}\M(v)\M(\vet)\bm{\Psi}_{e_{\l}}(|v-\vet|^{2})\,\d v\,\d\vet\,.$$
Then,
\begin{equation*} 
\int_{\R^{3}}|v|^{2}\left[\Q_{\re_{\l}}(\M,\M)-\mathfrak{a}_1\,\l^{\g}\mathrm{div}_v \left(v \M(v)\right)\right]\d v= 
6\mathfrak{a}_1\l^{\g}\en_{\star} -D_{\l}\,.
\end{equation*}
Compute the ratio
$$
\tau_{\l}:=\dfrac{D_{\l}}{6\mathfrak{a}_1 {\ell}^{\g}\en_{\star}}=\dfrac{1}{6\mathfrak{a}_1\l^{3+\g}\en_{\star}}\displaystyle\int_{\R^{3}\times\R^{3}}\M(v)\M(\vet)\bm{\Psi}_{\re}(\l^{2}|v-\vet|^{2})\,\d v\, \d\vet
$$
where we recall the definition~\eqref{def:Psire} and the fact that $\bm{\Psi}_{\re_{\l}}(r^{2})=\l^{-3}\bm{\Psi}_{\re}(\l^{2}r^{2})$.  Noticing that the Maxwellian defined in~\eqref{def:M} satisfies $\M(v)=\en_{\star}^{-\frac{3}{2}}M(\en_{\star}^{-\frac{1}{2}}v)$, it holds that
$$\int_{\R^{3} \times \R^{3}}\M(v)\M(\vet)\bm{\Psi}_{\re}(\l^{2}|v-\vet|^{2})\,\d v\,\d\vet=\int_{\R^{3} \times \R^{3}}M(v)M(\vet)\bm{\Psi}_{\re}(\l^{2}\en_{\star}|v-\vet|^{2})\,\d v\,\d\vet\,.
$$
We then use the unitary change of variable $y=\frac{v+\vet}{\sqrt{2}}$ and $y_*=\frac{v-\vet}{\sqrt{2}}$ combined with the identity~$M(v)M(\vet)=M(y)M(y_*)$ to deduce that
$$\int_{\R^{3} \times \R^{3}}\M(v)\M(\vet)\bm{\Psi}_{\re}(\l^{2}|v-\vet|^{2})\,\d v\,\d\vet=\int_{\R^{3}}M(y_*)\bm{\Psi}_{\re}(2\l^{2}\en_{\star}|y_*|^{2})\,\d y_*\,.$$
Moreover,
\begin{multline*}
\frac{1}{\l^{3+\g}}\bm{\Psi}_{\re}(2\l^{2}\en_{\star}|v|^{2})-\frac{\mathfrak{a}_{0}(2\en_{\star})^{\frac{\g+3}{2}}}{4+\g}|v|^{3+\g}\\
=\dfrac{(2\en_{\star})^{\frac{\g+3}{2}}|v|^{3+\gamma}}{2}\int_0^1 \left(\frac{1-\re^2(\l\sqrt{2\en_{\star}}\,|v|\,z)}{(\l\sqrt{2\en_{\star}}|v|\,z)^\gamma} -2\mathfrak{a}_{0}\right)z^{3+\gamma}\,\d z\,.
\end{multline*}
From \eqref{gamma0bis}, we deduce that
\begin{equation}\label{eq:PsiEL}
\left|\frac{1}{\l^{3+\g}}\bm{\Psi}_{\re}(2\l^{2}\en_{\star}|v|^{2})-\frac{\mathfrak{a}_{0}(2\en_{\star})^{\frac{\g+3}{2}}}{4+\g}|v|^{3+\g}\right| \lesssim \l^{\overline{\gamma}-\gamma} |v|^{3+\overline{\gamma}}+\l^\gamma\,|v|^{3+2\gamma}+  \l^{\overline{\gamma}} |v|^{3+\gamma+\overline{\gamma}}\end{equation}
holds for any $\l \geq0$ and $v \in\R^{3}$.  Thus,
$$\left|\tau_{\l}-\frac{1}{6\mathfrak{a}_1\en_{\star}}\frac{\mathfrak{a}_{0}(2\en_{\star})^{\frac{\g+3}{2}}}{4+\g}K_\gamma\right| \lesssim \max\left(\l^{\overline{\gamma}-\gamma},\l^{\g},\l^{\overline{\g}}\right)\,,$$
where $K_\gamma$ has been defined in~\eqref{def:Kgamma}.  Using the definition of $\mathfrak{a}_1$ in~\eqref{def:a1}, we obtain that 
$$\left|\tau_{\l}-1\right| \lesssim \max\left(\l^{\overline{\gamma}-\gamma},\l^{\g},\l^{\overline{\g}}\right).$$
This choice of $\mathfrak{a}_1$ is unique.  Indeed, once $\g,\en_{\star}$ have been fixed, there is a unique $\mathfrak{a}_1 >0$ such that $\lim_{\l\to0}\tau_{\l}=1$. This translates into 
$$
\left|D_{\l}-6\mathfrak{a}_1\en_{\star}\l^{\g}\right| \lesssim \max\left(\l^{\overline{\gamma}},\l^{2\g},\l^{\overline{\g}+\g}\right)\,,
$$
which gives the result.
\end{proof}
\subsection{Derivation of the sharp scaling}\label{sec:scaling}
We explicit the sharp self-similar scaling to observe the hydrodynamic behaviour of viscoelastic granular gases.  Recall from the introduction that, starting from a solution $F_{\e}$ to the free-cooling inelastic Boltzmann equation~\eqref{Bol-e}, we look for a self-similar change of variable as in~\eqref{eq:Fescal} 
with $V_{\e}(t) \geq0$ and
$\tau_{\e}(\cdot)\::\:[0,\infty) \to [0,\infty)$
which is a strictly increasing mapping with $\tau_{\e}(0)=0$ and $\lim_{t\to\infty}\tau_{\e}(t)=\infty$ for any $\e >0$.   The family of mappings $(f_{\e}(t))_{\e}$ is such that
$$f_{\e}(t,x,v)=\M(v)+\e\,h_{\e}(t,x,v)\,,$$
where $(h_{\e}(t))_{\e}$ is aimed to converge, as $\e \to 0$ to a profile  which depends on $(t,x)$ only through macroscopic the quantities $(\varrho(t,x),u(t,x),\vE(t,x)) \in \R^{+}\times\R^{3}\times \R^{3}$, namely, the hydrodynamic approximation.

We adopt a perturbative  approach using well-know results and properties of the elastic Boltzmann equation. This is possible because the Gaussian equilibrium $\M$ as well as the linearized and quadratic Boltzmann collision operators for elastic interactions emerge naturally as the gas evolves in the aforementioned scale, as noted in \cite{ALCMP}.  Thus, it is natural to assume that when $\e \simeq 0$, solutions $f_{\e}(t,x,v)$ solve a quasi conservative Boltzmann equation in a precise quantitative sense to be determined.
Consequently, it is central to the argument to determine the functions $V_{\e}(t)$ and $\tau_{\e}(t)$ to precisely balance the inflow-outflow of energy produced by the drift term and inelastic collisions along the hydrodynamic limit to capture the dynamics in an order one time-space scale. 

\color{black}

We observe that if $\tau_\e$ is such that~\eqref{eq:taueps} holds, 
the equation on the fluctuation $h_\e$ then writes~\eqref{eq:heps}
with
\begin{equation} \label{def:eepst}
\xi_{\e}(t):=\dot{V}_{\e}(s_{\e}(t))\,, \qquad 
\re_{\e,t}(r):=\re\left(\l_{\e}(t)\,r\right)\,,  \qquad \l_{\e}(t):=\frac{1}{V_{\e}(s_{\e}(t))}
\end{equation}
where $\dot{V}_{\e}(\cdot)$ denoting the derivative of the mapping $V_{\e}(\cdot)$ and $s_{\e}(\cdot)$ is the inverse function of the mapping~$t \mapsto \tau_{\e}(t)$.
The condition \eqref{eq:taueps} reads $\dot{\tau_{\e}}(t)=\frac{1}{V_{\e}(t)}$ and is assumed to get a constant coefficient in front of the $\partial_{t}f_{\e}$ term.  Because $\re(\cdot)$ is a non constant restitution coefficient, it is transformed into $\re_{\e,t}(\cdot)$  under the self-similar change of variable \eqref{eq:Fescal}.  More precisely, one can prove that, see \cite{CMS}, given a distribution $F=F(v)$ one has that
$$G_{\l}(v)=\l^3 F(\l v) \implies \l^2 \Q_\re(F,F)(\l v)=\Q_{\re_{\l}}(G_{\l},G_{\l})(v)\,, \qquad \forall \,(v,\ell) \in \R^{3} \times (0,\infty)\,,$$
where $\re_\l$ is defined in \eqref{def:e_ell}.
Recall from \eqref{eq:heps} that $\Q_{\e,t}=\Q_{\re_{\e,t}}$ while the linear operator $\LLe$ and the source term $S_\e$ are defined in~\eqref{def:LetSe}.
We point out that the condition  \eqref{gamma0bis} gives that for any $r > 0$, $t >0$, $\e \in (0,1]$,
\begin{equation}\label{eq:gamEet}
|\re_{\e,t}(r) - 1 + \mathfrak{a}_0 \l_{\e}(t)^\gamma r^\gamma|  \leq \mathfrak{b}_0 \l_{\e}(t)^{\bar \gamma} r^{\bar \gamma}\,. 
\end{equation}
The specific choice of the scaling function $V_{\e}(t)$ emerges by the following heuristics.
\begin{enumerate}[1)]
\item  {We want the linearized operator $\LLe$ to have finite kinetic energy for any $t >0$ as $\e \to 0$. The infinite energy case is of a vastly different nature and is out of reach of the method and tools used in this paper.} Notice also that, in the elastic case, the kinetic energy of the linearized term converges to $0$ (it is actually zero for any $\e >0$).  {Our goal being to get a non-trivial limit in the inelastic case, namely one that keeps track of the dissipative nature of granular gas in the macroscopic limit, we design our self-similar change of variables so that the kinetic energy of the linearized term is of order $1$ in the limit $\e \to 0$.} One can check without difficulty that
$$\e^{-2}\int_{\R^{3}}|v|^{2}\LLe h_{\e}(t,x,v)\,\d v=-2\e^{-2}\int_{\R^{3}} h_{\e}(t,x,v)\M(\vet)\bm{\Psi}_{\e,t}(|v-\vet|^{2})\,\d v\,\d\vet$$
where $\bm{\Psi}_{\e,t} := \bm{\Psi}_{\re_{\e,t}}$ and we recall that $\bm{\Psi}_\re$ is defined in~\eqref{def:Psire}. Assuming for a while that~$\l_\e(t)$ is converging to $0$ as $\e \to 0$, one can deduce from \eqref{eq:gamEet} that
$$\frac{1}{\e^{2}}\bm{\Psi}_{\e,t}(r^{2}) \simeq\frac{{1}}{\e^{2}}\frac{\mathfrak{a}_{0}}{4+\g}\l_{\e}(t)^{\g}r^{3+\g} \quad \text{as} \quad \e \simeq 0\,.$$
Here, we expect $(h_{\e})_{\e}$ to converge towards some ``hydrodynamic limit'' $\bm{h}(t,x,v)$. As such,
\begin{equation*}
\e^{-2}\int_{\R^{3}}|v|^{2}\LLe h_{\e}(t,x,v)\,\d v \simeq -\frac{2}{\e^{2}}\frac{\mathfrak{a}_{0}}{4+\g}\l_{\e}(t)^{\g}\int_{\R^{3}}\bm{h}(t,x,v)\M(\vet)|v-\vet|^{3+\g}\,\d v\,\d\vet\end{equation*}
 {In view of the above discussion}, it is then natural to expect here $\e^{-2}\l_\e(t)^\g$ to be of order one (w.r.t. $\e$), i.e. 
\begin{equation}\label{eq:eLeG}
\e^{-2}\l_\e(t)^\g \simeq  \nu(t) \quad \text{as} \quad \e \simeq 0 \end{equation}
for some suitable function $\nu(t)$ which may depend on $t$ but no longer depends on $\e$. 
 \item A second point relates to the connexion between $V_{\e}(t)$ and the Haff's law.  From the scaling~\eqref{eq:Fescal} and the definition of $\bm{T}_{\e}(t)$ in~\eqref{def:Teps}, one sees that
$$\bm{T}_{\e}(t)=V_{\e}(t)^{-2}\int_{\R^{3}}f_{\e}(\tau_{\e}(t),x,w)|w|^{2}\,\d w\,,$$
where, with \eqref{ansatz}, one expects that
$$\int_{\R^{3}}f_{\e}(\tau_{\e}(t),x,w)|w|^{2}\d w \simeq \int_{\R^{3}}\M(w)|w|^{2}\d w=3\en_{\star} \quad \text{for} \quad  \e \simeq 0\,.$$
This means that
\begin{equation}\label{eq:thermalspeed}\bm{T}_{\e}(t) \simeq V_{\e}(t)^{-2}\,\end{equation}
The equivalence \eqref{eq:thermalspeed} express the physical meaning of the scaling function $V_\e(t)$ seen as the square root of the inverse of the gas'\emph{thermal velocity}, as governed by Haff's law \cite[Chapter 7, Eq. (7.2)]{Poschel}.  
A rigorous proof of such law is given in \cite{SIAM} (in the spatially homogeneous case) showing that
$$\bm{T}_{1}(t) \simeq (1+t)^{-\frac{2}{\g+1}},$$
suggesting the \emph{ansatz}
\begin{equation}\label{eq:temp}V_{\e}(t)=\left(\mathfrak{b}_{\e} + \mathfrak{a}_{\e}t\right)^{\frac{1}{\g+1}}\end{equation}
for some positive parameters $\mathfrak{b}_{\e} >0, \mathfrak{a}_{\e} >0$ to be determined.
\end{enumerate}
Other technical considerations determine $V_\e(t)$ as well. Indeed, it is convenient to write \eqref{eq:heps} as a perturbation of the classical elastic Boltzmann equation introducing the elastic operators $\Q_{1}$ and $\mathscr{L}_{1}$ as much as possible to find the equivalent formulation
\begin{multline}\label{eq:hL1}
\partial_{t} h_{\e} + \frac1\e v \cdot \nabla_x h_\e 
=
\frac{1}{\e^2} \mathscr{L}_{1} h_\e+ \frac{1}{\e}\Q_{1}(h_{\e},h_{\e}) \\
+\underset{\textrm{new error term } \widetilde{S}_{\e}(t,x,v)}{\underbrace{\frac{1}{\e^{2}}\left(\LLe h_\e - \mathscr{L}_{1} h_{\e}\right) + \frac1\e \left(\Q_{\e,t} (h_\e,h_\e)-\Q_{1}(h_{\e},h_{\e}) \right)-\xi_{\e}(t)\mathrm{div}_{v}(vh_{\e})+ S_\e(t,v)}}
\end{multline}
where the last line is understood as a source term whose first two terms are aimed to be small. 
\begin{itemize}
\item[--] It is natural to assume that the drift term and $\Q_{\e,t}-\Q_{1}$ to be of the same magnitude, that is, 
$$\e^{2}\xi_{\e}(t) \simeq \|\Q_{\e,t}(\cdot,\cdot)-\Q_{1}(\cdot,\cdot)\|$$
for a suitable $L^{1}$-norm $\|\cdot\|$, see \cite{ALhaff,ALCMP} for the homogeneous case.
Moreover, on the basis of the result of Proposition \ref{prop:Qel}, 
one is led to the scaling
\begin{equation}\label{eq:LEXI}
\e^{2}\xi_{\e}(t) \simeq \l_{\e}(t)^{\g}.
\end{equation}
Using \eqref{eq:eLeG} amounts to assume $\xi_\e(t)$ is independent of $\e$, at least for small values of $\e$.
\item[--] The source term $S_{\e}$ defined in~\eqref{def:LetSe} acts as a controlled and balanced term along the hydrodynamic limit. The problem is that $\|S_{\e}\|$ has a penalisation term of order $\e^{-3}$. However, our analysis only needs the smallness of $\bm{\pi}_{0}S_{\e}$, where $\bm{\pi}_{0}$ denotes the projection on the kernel of $\mathscr{L}_{1}$ (see~\eqref{def:pi0} for the explicit definition),  which is the relevant term for the macroscopic effects of $S_{\e}$.   Additionally, due to the conservation properties of the collision operator and the symmetry properties of the Maxwellian state $\M$ it follows that 
\begin{equation} \label{eq:conservSeps}
\int_{\R^3} S_\e(t,v) \left(\begin{array}{c}1 \\v \end{array}\right) \, \d v = \left(\begin{array}{c}0 \\0 \\\end{array}\right) 
\qquad \forall \,t \geq 0\,, \, \forall\,\e>0\,.
\end{equation}
Therefore, the balanced and controlled nature of $\bm{\pi}_{0}S_{\e}$ only relates to its kinetic energy flux $$\int_{\R^{3}}|v|^{2}S_{\e}(t,v)\,\d v \simeq 0 \quad \text{as} \quad \e \to 0\,.$$
Using Lemma \ref{lem:Ener} together with \eqref{eq:LEXI}, this is possible only if
\begin{equation}\label{eq:xieA1}
\e^{2}\xi_{\e}(t)=\mathfrak{a}_{1}\l_{\e}(t)^{\g}
\end{equation}
where $\mathfrak{a}_{1}$ is given by \eqref{def:a1}. 
\end{itemize}
Let us proceed as follows: insert the ansatz \eqref{eq:temp} in the definition \eqref{eq:taueps} of $\tau_\e(\cdot)$, one computes that
\begin{equation} \label{eq:taueps}
\tau_{\e}(t)=\frac{1+\g}{\g\mathfrak{a}_{\e}}\left[\left(\mathfrak{b}_{\e}+\mathfrak{a}_{\e}t\right)^{\frac{\g}{\g+1}}-\mathfrak{b}_{\e}^{\frac{\g}{\g+1}}\right]\,, 
\end{equation}
with inverse function
$$s_{\e}(t)=\frac{1}{\mathfrak{a}_{\e}}\left[\left(\mathfrak{b}_{\e}^{\frac{\g}{\g+1}}+\frac{\g\mathfrak{a}_{\e}}{1+\g}\,t\right)^{\frac{1+\g}{\g}}-\mathfrak{b}_{\e}\right]\,.$$
Then, from \eqref{def:eepst}, it follows that
$$\xi_{\e}(t)=\dot{V}_{\e}(s_{\e}(t))=\frac{\mathfrak{a}_{\e}}{1+\g}\left(\mathfrak{b}_{\e}^{\frac{\g}{\g+1}}+\frac{\g\mathfrak{a}_{\e}}{1+\g}\,t\right)^{-1}$$
and
$$\l_{\e}(t)=\frac{1}{V_{\e}(s_{\e}(t))}=\left(\mathfrak{b}_{\e}^{\frac{\g}{\g+1}}+\frac{\g\mathfrak{a}_{\e}}{1+\g}\,t\right)^{-\frac{1}{\g}}\,.$$
On the one hand, \eqref{eq:xieA1} implies that
$$\mathfrak{a}_{\e}=(1+\g)\e^{-2}\mathfrak{a}_{1}\,.$$
On the other hand, one observes that 
$\xi_{\e}(0)=\mathfrak{a}_{1}\e^{-2}\mathfrak{b}_{\e}^{-\frac{\g}{\g+1}}$
and, since we look for $\xi_\e(t)$ to be independent of $\e$, one assumes
$$\mathfrak{b}_{\e}=\e^{-\frac{2(1+\g)}{\g}}\mathfrak{b}_{1}$$
for some arbitrary $\mathfrak{b}_{1} >0$.  Such a choice leads to
\begin{equation}
\label{eq:Ve}
V_{\e}(t)=\left(\e^{-\frac{2(1+\g)}{\g}}\mathfrak{b}_{1}+(1+\g)\mathfrak{a}_{1}\e^{-2}t\right)^{\frac{1}{\g+1}}=\e^{-\frac{2}{\g}}\left(\mathfrak{b}_{1}+(1+\g)\mathfrak{a}_{1}\e^{\frac{2}{\g}}t\right)^{\frac{1}{\g+1}}, \end{equation}
and
\begin{equation}\label{def:xi}
\l_{\e}(t)=\e^{\frac{2}{\g}}\left(\mathfrak{b}_{1}^{\frac{\g}{\g+1}}+ \g\mathfrak{a}_{1}\,t\right)^{-\frac{1}{\g}}, \qquad \xi_{\e}(t)=\xi(t)=\mathfrak{a}_{1}\left(\mathfrak{b}_{1}^{\frac{\g}{\g+1}}+ \g\mathfrak{a}_{1}\,t\right)^{-1}\,.
\end{equation} Consequently,
\begin{equation} \label{def:zt}
\l_{\e}(t)=\e^{\frac{2}{\g}}z(t) \quad \text{with} \quad  z(t):=\left(\mathfrak{b}_{1}^{\frac{\g}{\g+1}}+ \g\mathfrak{a}_{1}\,t\right)^{-\frac{1}{\g}}\end{equation}
where $z(\cdot)$ is also independent of $\e$. 

\begin{nb}
 
The constant $\mathfrak{a}_1$ is chosen satisfying~\eqref{def:a1} so that the conclusion of Lemma~\ref{lem:Ener} holds. We shall also need an upper bound on $\mathfrak{b}_{1}$ with respect to some constants given by the problem that will allow us to use a perturbative approach around the elastic case. More precisely, we suppose that~$\mathfrak{b}_1$ is small enough so that the condition on $\lambda_\star$ given in Proposition~\ref{prop:Psi} is satisfied. We do not enter into details here into the link between $\lambda_\star$ and $\mathfrak{b}_1$. However, we mention that we shall see in the subsequent analysis that the condition on $\lambda_\star$ can actually be reformulated in terms of a condition on the supremum of $z(\cdot)$ and thus in terms of an upper bound on $\mathfrak{b}_1$. 
\end{nb}

 \color{black}
  
\section{The complete inelastic linearized operator}\label{sec:LIN}  Recalling that the collisional elastic operator $\LL$ has been defined in~\eqref{def:LL}, the complete linearized operator for the \emph{elastic operator} is defined for any $\e \in (0,1]$ by
\begin{equation} \label{def:Ge}
\mathcal{G}_{\e}h:=\e^{-2} \LL(h) - \e^{-1}v \cdot \nabla_{x} h\,, 
\end{equation}
while the rescaled inelastic (time-depending) linearized operator is given by
\begin{equation} \label{def:Get}
\mathcal{G}_{\e,t}h:=\e^{-2} \LLe(h) - \e^{-1} v \cdot \nabla_{x}h - \xi(t) \mathrm{div}_v( v   h)\,, 
\end{equation}
where we recall that $\LLe$ and $\xi(t)$ have been defined respectively in~\eqref{def:LetSe} and~\eqref{def:xi}. 
We begin with recalling several properties of the \emph{elastic linearized operator}.
\subsection{Decomposition of $\LL$}\label{sec:hypo} Let us now recall the following decomposition of $\LL$  introduced in {\cite{GMM,Tr}} {(see also~\cite{AGT} for proofs adapted to the case of  a general collision kernel $b$)}. For any $\delta\in (0,1)$, we consider the cutoff function $0\leq\Theta_\delta = \Theta_\delta(v,\vb	, \sigma) \in\Cs^{\infty}(\R^{3}\times \R^{3}\times \S^{2})$, assumed to be bounded by $1$, which equals $1$ on 
$$J_{\delta}:=\left\{(v,\vb,\sigma)\in \R^{3}\times \R^{3}\times \S^{2}\,\Big|\,|v|\leq \delta^{-1}\,,\, 2\delta \leq |v-\vb|\leq \delta^{-1}\,,\,
 |\!\cos\theta| \leq 1-2\delta \right\},$$
and whose support is included in $J_{\frac{\delta}{2}}$ (where $\cos \theta=\langle \frac{v-\vb}{|v-\vb|},\sigma\rangle$). We then set 
\begin{equation*}
\begin{split}
\mathscr{L}_{1}^{S,\delta}h(v)& =\int_{\R^3\times \S^{2}} 
\big[\M(\vb')h(v') +\M(v')h(\vb')-\M(v)h(\vb)\big]\\
&\hspace{6cm}\times|v-\vb|\,{ b(\cos \theta)}\,\Theta_\delta (v,\vb,\sigma) \, \d\vb\, \d\sigma\,, \\ 
\mathscr{L}_{1}^{{R,\delta}}h(v)&=  \int_{\R^{3}\times \S^{2}}
\big[\M(\vb')h(v') +\M(v')h(\vb')-\M(v)h(\vb)\big] \\
&\hspace{6cm}\times|v-\vb|\,{ b(\cos \theta)}\,(1-\Theta_\delta(v,\vb,\sigma) ) \, \d\vb\, \d\sigma\,,
\end{split}\end{equation*}
where of course here, $v',\vb'$ are denoting the \emph{elastic post-collisional velocities}. 
It holds then that~$\mathscr{L}_{1}= \mathscr{L}_{1}^{{S,\delta}}+ \mathscr{L}_{1}^{{R,\delta}}-\Sigma_\M$
where using~\eqref{unitb}, $\Sigma_{\M}$ denotes the mapping 
\begin{equation}\label{eq:SigmaM}
\Sigma_{\M}(v):=\int_{\R^{3}}\M(\vb)|v-\vb	|\, \d\vb\,, \qquad \forall \, v\in \R^{3}\,.\end{equation}
 {Recall that there exist $\sigma_0>0$ and $\sigma_1>0$ such that 
\begin{equation} \label{eq:collfreq}
\sigma_0 \, \m_1(v) \leq \Sigma_{\M}(v) \leq \sigma_1 \, \m_1{(v)}\,, \qquad \forall \, v \in \R^{3}\,. 
\end{equation}}
Introduce
\begin{equation*} {A}^{(\delta)} :=  \mathscr{L}_{1}^{S,\delta}\qquad \text{ and } \qquad
{B}_{1}^{(\delta)}  := \mathscr{L}_{1}^{{R,\delta}}-\Sigma_\M
\end{equation*}
so that $\mathscr{L}_{1}={A}^{(\delta)}+  {B}_{1}^{(\delta)}$.  Let us now recall the known dissipativity results for the elastic Boltzmann operator in $L^{1}_{v}L^{2}_{x}$-based Sobolev spaces, see {\cite[Lemmas~4.12,~4.14 \& Lemma~4.16]{GMM}}:

\begin{lem}\phantomsection\label{lem:dissip}For any $k \in \N$ and $\delta >0,$ there are two positive constants $C_{k,\delta} >0$ and $R_{\delta} >0$ such that
$\mathrm{supp}\left({A}^{(\delta)}h\right)\subset B(0,R_{\delta})$
and
\begin{equation}\label{eq:Adelta}
 {\| {A}^{(\delta)}h\|_{\W^{k,2}_{v}(\R^{3})}} \leq C_{k,\delta}\|h\|_{L^{1}_{v}(\m_{1})}\,, \qquad \forall \,  h \in L^{1}_v(\m_{1})\,.\end{equation}
Moreover, the following holds:

for any $q > 2$ and any $\delta \in (0,1)$ it holds
{\begin{multline}\label{eq:B1delta}
\int_{\R^d} \|h(\cdot,v)\|_{L^{2}_{x}}^{-1} \left(\int_{\T^3} \left(B_{1}^{(\delta)}h(x,v)\right) h(x,v) \, \d x\right) \m_q(v)\, \d v \\
\leq \left(\Lambda_{q}(\delta)-1\right)\|h\|_{L^{1}_{v}L^{2}_{x}(\m_{q}\Sigma_\M)}
\end{multline}}
where $\Lambda_{q}\::\:(0,1) \to \R^{+}$ is some explicit function such that $\lim_{\delta\to0}\Lambda_{q}(\delta)=\frac{4}{q+2}.$
\end{lem}\phantomsection

\subsection{Splitting of $\mathcal{G}_{\e,t}$}

With the previous decomposition of $\LL$, we can split $\LLe$ as
\begin{equation} \label{eq:splitLalpha}
\LLe={B}_{\e,t}^{(\delta)} + {A}^{(\delta)} \qquad \text{where} \qquad {B}_{\e,t}^{(\delta)}:={B}_{1}^{(\delta)}+ \LLe-\LL 
\end{equation} 
which depend naturally on $\delta,\e$ and $t$. We insist here on the fact that these linearized operators are depending on $t$. Such a splitting yields the decomposition of $\mathcal{G}_{\e,t}$ as
$$\mathcal{G}_{\e,t}=\mathcal{A}_{\e}^{(\delta)}+ \mathcal{B}_{\e,t}^{(\delta)}$$
where
$$\mathcal{A}_{\e}^{(\delta)}h:=\e^{-2} A^{(\delta)}h\,,$$
and
$$\mathcal{B}_{\e,t}^{(\delta)}h:=\e^{-2} {B}_{\e,t}^{(\delta)} -\e^{-1} v\cdot \nabla_{x}h\,-\xi(t)\mathrm{div}_v (v h)\,.$$
One has the following properties of $\mathcal{B}_{\e,t}^{(\delta)}$ in~$L^1_vL^2_x$-based spaces. 
\begin{prop}\phantomsection\label{prop:dissip} 
For any $k\geq 0$ and $q >2$, there exist $\e_0>0$, ${\delta}_{q}^{\dagger} >0$ and ${\nu}_{q} >0$ such that  
$$\mathcal{B}_{\e,t}^{(\delta)} + \e^{-2}{\nu}_{q} \;\, \text{ is dissipative in $L^1_v \W^{k,2}_{x}(\m_{q})$} 
$$
for any  $\delta \in (0,{\delta}_{q}^{\dagger})$, any $t \geq 0$ and any $\e \in (0,\e_0)$.  
\end{prop}
\begin{nb} \label{rem:dissip}
Let us be more precise on the estimate of dissipativity we obtain in $L^1_vL^2_x(\m_q)$ for further use: for any $\e \in (0,\e_0)$, any $t \geq 0$ and any $\delta \in (0,{\delta}^{\dagger}_{q})$, we have
$$
\int_{\R^{3}} \|h(\cdot,v)\|_{L^2_x}^{-1}\left(\int_{\T^3}  \mathcal{B}_{\e,t}^{(\delta)}(h)(x,v)  h(x,v) \, \d x\right) \, \m_{q}(v)\, \d v  
\leq - \e^{-2} \nu_q \|h\|_{L^1_vL^2_x(\m_{q+1})}\,. 
$$
\end{nb}
\begin{proof}
Notice that derivatives with respect to the $x$-variable commute with the operator $\mathcal{B}_{\e,t}^{(\delta)}$ and this allows to prove the result, without loss of generality, in the special case $k=0$. The proof follows the similar one in \cite[Proposition 2.8]{ALT} that we recall here. From the definition of $\B_{\e,t}^{(\delta)}$, one has 
$$\mathcal{B}_{\e,t}^{(\delta)}(h)=\sum_{i=0}^{3}C_{i}(h)$$ 
with 
$$C_{0}(h):=\e^{-2}B_{1}^{(\delta)}h\,, \qquad  C_{1}(h):=-\e^{-1} v \cdot \nabla_{x}h\,,\quad C_{2}(h):=\e^{-2}\left(\LLe-\LL\right)h\,,$$
and  $C_{3}(h):=-\xi(t)\mathrm{div}_{v}(vh)\,,$ and correspondingly (with obvious notations),
{\begin{equation*}
\int_{\R^{3}} \|h(\cdot,v)\|_{L^2_x}^{-1}\left(\int_{\T^3}  \mathcal{B}_{\e,t}^{(\delta)}(h)(x,v)  h(x,v) \, \d x\right) \, \m_{q}(v)\, \d v  
 =:\sum_{i=0}^{3}I_{i}(h)\,.
 \end{equation*}} 
First, since $C_{1}$ is a divergence operator with respect to $x,$ one has $I_1(h)=0$. According to~\eqref{eq:B1delta}, by taking $\delta$ small enough so that $\Lambda_q(\delta)<1$ (which is possible since $q>2$), we have
$$I_{0}(h) \leq \e^{-2} \sigma_0\left(\Lambda_{q}(\delta)-1\right)\|h\|_{L^{1}_{v}L^{2}_{x}(\m_{q+1})}\,.$$
Then, it follows from Cauchy-Schwarz inequality that
$$
I_2(h) \leq 
\e^{-2} \int_{\R^{3}} \|(\LLe-\LL) h(\cdot,v)\|_{L^2_x}  \m_q(v) \, \d v =\e^{-2}\left\|\LLe h-\LL h\right\|_{L^{1}_{v}L^{2}_{x}(\m_{q})}\,.
$$
Now, recalling that $\LLe h=\mathscr{L}_{\re_{\l}}$ with $\l=\l_{\e}(t)$ defined in~\eqref{def:xi}, one deduces from Corollary~\ref{cor:Line} that
$$\left\|\LLe h-\LL h\right\|_{L^{1}_{v}L^{2}_{x}(\m_{q})} \lesssim \l_{\e}(t)^{\frac{\g}{8+3\g}}\|h\|_{L^{1}_{v}L^{2}_{x}(\m_{q+1})}$$
provided $\l_{\e}(t) < \l_{\star}.$ We recall that $\l_{\e}(t)=\e^{\frac{2}{\g}}z(t)$ with $z(t)$ defined in~\eqref{def:xi}. It satisfies $z(t) \leq z(0) = \mathfrak{b}_1^{-1/(\g+1)}$ for any $t \geq 0$, we thus see that by choosing $\e_{\star}:=\l_\star^{\frac{\g}{2}} \mathfrak{b}_1^{\frac{\g}{\g+1}}$, we will have that~$\l_{\e}(t) < \l_{\star}$ for any $\e \in (0,\e_{\star})$ and any $t\geq0.$ Then, we deduce that there exists $c_0>0$ such that for any $\e \in (0,\e_\star)$ and any~$t \geq 0$, we have 
$$\left\|\LLe h-\LL h\right\|_{L^{1}_{v}L^{2}_{x}(\m_{q})}\leq c_0 \, \e^{\frac{2}{8+3\g}}\|h\|_{L^{1}_{v}L^{2}_{x}(\m_{q+1})}\,.$$
Finally, direct computations as in \cite[Proposition 2.8]{ALT} yield
$$I_{3}(h) \leq q\, \xi(t) \|h\|_{L^{1}_{v}L^{2}_{x}(\m_{q+1})}\,.$$
Using the definition of $\xi(t)$ in~\eqref{def:xi}, we notice that  $\xi(t) \lesssim 1$, which implies that there exists a constant $c_1>0$ such that for any $t \geq 0$,
\begin{equation} \label{eq:I3}
I_{3}(h) \leq c_1 \|h\|_{L^{1}_{v}L^{2}_{x}(\m_{q+1})}\,.
\end{equation}
Gathering the previous estimates, one obtains 
\begin{equation}\label{eq:Bal}
\begin{split}
&\int_{\R^{3}} \|h(\cdot,v)\|_{L^2_x}^{-1} \left(\int_{\T^3}  \mathcal{B}_{\e,t}^{(\delta)}(h)(x,v)  h(x,v) \, \d x\right) \, \m_{q}(v)\, \d v \\
&\qquad \leq \e^{-2}\,\left(\sigma_0\,(\Lambda_{q}(\delta)-1) + c_{0}\,\e^{\frac{2}{8+3\g}}  +c_{1} \,\e^2\right)\|h\|_{L^{1}_{v}L^2_{x}(\m_{q+1})}\,.
\end{split}
\end{equation}
Recalling that $\lim_{\delta\to 0}(\Lambda_{q}(\delta)-1)=-\frac{q-2}{q+2} <0$, we can pick $\delta_{q}^{\dagger}$ and $\e_0 \in (0,\e_\star)$ small enough  so that
  $$
  - \nu_{q}:= \sup \left\{c_{0}\,\e^{\frac{2}{8+3\g}}  + \sigma_0\,(\Lambda_{q}(\delta)-1) +c_{1}\,\e^2\,;\,\e \in (0,\e_0),\, \delta \in (0,\delta_{q}^{\dagger})  \right\} <0 $$
and get the result. 
\end{proof}

For the rest of the paper, we fix $\delta \in (0,\delta^\dagger_q)$ (with {$q$ fixed in~\eqref{eq:defE}) and $\delta^\dagger_q$ defined in Proposition~\ref{prop:dissip}) so that the conclusions of Proposition~\ref{prop:dissip} and Remark~\ref{rem:dissip} are satisfied for~
$$\B_{\e,t} := \B_{\e,t}^{(\delta)}\,, \qquad \forall \, \e \in (0,\e_0)\,, \qquad \forall \,t \geq0$$
in the functional space~$\E=L^{1}_{v}(\W^{m,2}_{x}(\m_{q}))$ introduced in~\eqref{eq:defE}.

\begin{nb}\label{rem:dissip2} The above result has the next consequence on semilinear evolution equations governed by the operator~$\B_{\e,t}^{(\delta)}$. Given $h_{0} \in \E$ and $g=g(t) \in \E$, one can compute the evolution of $\|h(t)\|_{\E}$ of any solution to
$$\partial_{t}h(t)=\B_{\e,t} h(t) + g(t), \qquad h(t=0)=h_{0}\,.$$
Indeed, recalling that 
$$\|h(t)\|_{\E}=\sum_{0 \leq s \leq m}\left\|\|\nabla_{x}^{s}h(t,\cdot,v)\|_{L^{2}_{x}}\m_{q}(v)\right\|_{L^{1}_{v}}$$
one sees that
\begin{multline*}
\frac{\d}{\d t}\|h(t)\|_{\E}=\sum_{s=0}^m\int_{\T^{3}\times\R^{3}}\left\|\nabla_{x}^{s} h(t,\cdot,v)\right\|_{L^{2}_{x}}^{-1}\nabla_{x}^{s}(\partial_{t}h(t,x,v)) \cdot \nabla_{x}^{s}h(t,x,v)\m_{q}(v)\,\d v\,\d x\\
=\sum_{s=0}^m\int_{\T^{3}\times\R^{3}}\left\|\nabla_{x}^{s} h(t,\cdot,v)\right\|_{L^{2}_{x}}^{-1}\left[\B_{\e,t}\left(\nabla_{x}^{s}h(t,x,v)\right)+\nabla_{x}^{s}g(t,x,v)\right] \cdot \nabla_{x}^{s}h(t,x,v)\m_{q}(v)\,\d v\,\d x
\end{multline*}
where we used that $\B_{\e,t}$ and $\nabla_{x}^{s}$ commute. Therefore, according to Remark~\ref{rem:dissip}, for any $\e \in (0,\e_0)$ and any $t\geq 0$,
\begin{equation}\label{eq:remdissip}
\dfrac{\d}{\d t}\|h(t)\|_{\E} \leq -\e^{-2}\nu_{q} \|h(t)\|_{L^{1}_{v}(\W^{m,2}_{x}(\m_{q+1}))}+\|g(t)\|_{\E}\,.\end{equation}
\end{nb}

 \subsection{Hypocoercive norm for the elastic problem} \label{subsec:hypo}

We recall here some hypocoercivity results for the elastic linearized operator $\G_{\e}$ defined in~\eqref{def:Ge} as derived in our previous contribution \cite{ALT} and based upon results from \cite{BCMT,CRT}. We point out that these hypocoercivity results are involving spaces which \emph{do  not involve any derivative in the $v$-variable} and apply namely to the spaces $\H$ and $\H_1$ introduced in~\eqref{eq:defH}-\eqref{eq:defH1}.

Before stating our result, we recall the expression for the spectral projection $\bm{\pi}_{0}$ onto the kernel~$\mathrm{Ker}(\LL)$ of the linearized collision operator $\LL$ \emph{seen as an operator acting in velocity only} on the space $L^{2}_{v}(\M^{-\frac{1}{2}})$: 
\begin{equation}\label{def:pi0}
\bm{\pi}_{0}(g) :=\sum_{i=1}^{5}\left(\int_{ \R^{d} }g\,\varphi_{i}\,\d v \right)\,\varphi_{i}\,\M\,,
\end{equation}
where 
\begin{equation} \label{def:Psii}
\varphi_{1}(v):=1\,,\quad \varphi_{i}(v):=\frac{v_{i-1}}{\sqrt{\en_{\star}}} \quad  \quad (i=2,\ldots,4)\,, \quad \text{and} \quad \varphi_{5}(v):=\frac{1}{\en_{\star}\sqrt{6}}(|v|^{2}-3\en_{\star})\,,
\end{equation}
where
$$\int_{\R^{3}}\varphi_{i}(v)\varphi_{j}(v)\M(v)\, \d v=\delta_{i,j}\,, \qquad \forall \, i,j=1,\ldots,5\,.$$
 Note also (see for example~\cite{briant}) that the spectral projection $\mathbf{P}_{0}$ onto the kernel of $\G_{\e}$ defined in~\eqref{def:Ge} is given by 
 \begin{equation}\label{eq:P0}
\mathbf{P}_{0}(g) :=\sum_{i=1}^{5}\left(\int_{\T^{3}\times\R^{3}}g\,\varphi_{i}\,\d v\, \d x\right)\,\varphi_{i}\,\M\end{equation}
and so the difference with respect to the spectral projection $\bm{\pi}_0$ is that an additional spatial integration is performed.

\begin{prop}\label{prop:hypoco} On the space $\H$ defined in \eqref{eq:defH}, there exists a norm $\vertiii{\cdot}_{\H}$ with associated inner product $\langle\!\langle\cdot,\cdot\rangle\!\rangle_{\H} $ equivalent to the standard norm $\|\cdot\|_{\H}$ for which 
 there exist $\mathrm{a}_{1}>0$ and~$\mathrm{a}_{2} >0$ such that for any $\e \in (0,1]$,
\begin{equation}\label{eq:hypoco}
\langle\! \langle\, \G_{\e}h,h\,  \rangle\!\rangle_{\H}  \leq -\frac{\mathrm{a}_{1}}{\e^{2}}\left\|\left(\mathbf{Id}-\bm{\pi}_{0}\right)h\right\|_{\H_{1}}^{2} -\mathrm{a}_{1}\|h\|_{\H_{1}}^{2}-\mathrm{a}_{2}\vertiii{h}_{\H}^{2}
\end{equation}
holds true for any $h=(\mathbf{Id-P}_{0})h \in \D(\G_{\e})\subset \H$ where $\bm{\pi}_{0}$ (resp. $\mathbf{P}_0$) is defined in \eqref{def:pi0} (resp.~\eqref{eq:P0}). 
Moreover, there exists $\mathrm{a}_{3} >0, C >0$ (independent of $\e$) such that
\begin{equation}
\label{lem:hypodrift}
-\langle\! \langle \, \mathrm{div}_v  \left(v h\right),h\,\rangle\!\rangle_{\H} \leq \mathrm{a}_{3}\langle\!\langle \m_{2}h,h\rangle\!\rangle_{\H}+C\e\|h\|_{\H}^{2}\end{equation}
for any $h \in L^2_v \W^{m,2}_x(\langle \cdot \rangle \M^{-\frac12})$.  
\end{prop}
{
\begin{nb}
Remark that the equivalent norm $\vertiii{\cdot}_{\H}$ actually depends on $\e$  but we do not mention this dependency in our notation because this norm is equivalent to the usual one~$\|\cdot\|_\H$ \emph{uniformly} in~$\e \in (0,1]$ since  there exists $C_{\H} >0$ independent of $\e \in (0,1]$ such that
\begin{equation}\label{eq:CH-CH}
C_{\H}\|h\|_{\H} \leq \vertiii{h}_{\H} \leq C_{\H}^{-1}\|h\|_{\H}\,\,, \qquad \forall \, h \in \H\,.\end{equation}
\end{nb}}
We have the following result which is an easy adaptation of \cite[Lemma 3.4]{ALT} in the spirit of \cite[Theorem 4.7]{bmam}. 
\begin{lem}\label{theo:briant2.4}
The inner product $\langle\!\langle\cdot,\cdot\rangle\!\rangle_{\H}$ associated to the norm $\vertiii{\cdot}_{\H}$ on $\H$ constructed in Proposition \ref{prop:hypoco} is such that for any $h_{1},h_{2} \in \H_1$, $g \in \H_{1}{\cap \mathrm{Range}(\mathbf{Id-P}_{0})}$ (where we recall that $\H_1$ is defined in~\eqref{eq:defH1}), for any $\e \in (0,1]$ and any $t \geq 0$,
\begin{multline}\label{eq:briant}
\langle\!\langle \left(\mathbf{Id}-\bm{\pi}_{0}\right)\Q_{\e,t}(h_{1},h_{2}),g\rangle\!\rangle_{\H}
\lesssim \left(\|h_{1}\|_{\H_{1}}\|h_{2}\|_{\H}+\|h_{1}\|_{\H}\|h_{2}\|_{\H_{1}}\right)\,\left\|\left(\mathbf{Id}-\bm{\pi}_{0}\right)g\right\|_{\H_{1}}\\
+ { \e\|h_{1}\|_{\H}\|h_{2}\|_{\H}\|\bm{\pi}_{0}g\|_{\H}}\,.
\end{multline}
\end{lem}
\begin{proof} The proof is exactly the same as that of \cite[Lemma 3.4]{ALT} since we check easily that, for any restitution coefficient $\re(\cdot)$ of class $\mathcal{R}_{\g}$ and any $\l \lesssim 1$, the estimates 
\begin{multline*}
\left\langle \left(\mathbf{Id}-\bm{\pi}_{0}\right)\Q_{\re_{\l}}(h_{1},h_{2}), \bm{\pi}_{0}g\right\rangle_{\H}=0\,, \quad
 \left\langle \left(\mathbf{Id}-\bm{\pi}_{0}\right)\Q_{\re_{\l}}(h_{1},h_{2}), g^{\perp}\right\rangle_{\H}=\left\langle  \Q_{\re_{\l}}(h_{1},h_{2}), g^{\perp}\right\rangle_{\H}\end{multline*}
are valid with $g^{\perp}=(\mathbf{Id}-\bm{\pi}_{0})g$, $g \in \mathrm{Range}(\mathbf{Id-P}_{0})$. Then, introducing as in \cite[Lemma 3.4]{ALT}, the space
$$\H_{-1}:=L^{2}_{v}\W^{m,2}_{x}(\M^{-\frac{1}{2}}\langle \cdot\rangle^{-1})\,,$$ 
one has
$$\left\|\Q_{\re_{\l}}(h_{1},h_{2})\right\|_{\H_{-1}}
\lesssim \|h_{1}\|_{\H}\|h_{2}\|_{\H}$$
and one can estimate $\left\langle \Q_{\re_{\l}}(h_{1},h_{2}), g^{\perp}\right\rangle_{\H}$ using Lemma \ref{lem:trilinear} to conclude as in \cite[Lemma 3.4]{ALT} where we observed that all constants are independent of $\l \lesssim 1$ thanks to Lemma~\ref{lem:Lam} and Remark~\ref{nb:Lambdae}. Finally, we remark that $\l_\e(t) \lesssim 1$ uniformly in $\e \in (0,1]$ and $t \geq 0$, which yields the conclusion.
\end{proof}
\section{Nonlinear analysis}\label{sec:nonlinear}
 
We now apply the results obtained so far to the study of Eq.~\eqref{eq:heps}.  {Let us recall that the functional spaces $\E$ and $\E_1$ in which we are going to carry out our analysis are defined in~\eqref{eq:defE}-\eqref{eq:defE1} and that 
the Hilbert spaces $\H$ and $\H_1$ are defined in~\eqref{eq:defH}. 
We also recall here that the Maxwellian~$\M$ introduced in~\eqref{eq:max} is the steady state of $\mathscr{L}_{1}$ defined in~\eqref{def:LL} while $\H$ is a Hilbert space on which the elastic Boltzmann equation is well-understood (see Subsection~\ref{subsec:hypo}).
{Motivated by the estimate given in Lemma \ref{lem:LeL1}, we also define the following functional space
\begin{equation}\label{def:E1E2}
 \E_{2}:=L^2_{v}\W^{m,2}_{x}(\m_{q+{2\kappa}+\g+3})\,,  \qquad \kappa > \frac{3}{2},
\end{equation}
which is such that
$\E_{2} \hookrightarrow L^{1}_{v}\W^{m,2}_{x}(\m_{q+\g+3})$.}

\medskip
\noindent
Notice here that we imposed {the condition $m>\frac{3}{2}$  in order to ensure that the continuity of the embedding  $\W^{m,2}_{x}(\T^{3}) \hookrightarrow L^{\infty}_{x}(\T^{3})$, which allows us to treat nonlinear terms thanks to the underlying Banach algebra structure.  Notice also that our analysis is based on the fact that $\mathcal{A}_{\e}$ has regularization properties in velocity, namely $\mathcal{A}_\e \in \mathscr{B}(\E,\H)$ and on the following continuous embeddings: 
$$\H \hookrightarrow \E_{2} \hookrightarrow \E_1 \hookrightarrow \E.$$ Since~$\mathcal{A}_{\e}$ has no regularisation effect on the spatial variable, we are forced to have the same number of spatial derivates in the spaces $\E$ and $\H$. Taking $q \geq3+\g$ in~\eqref{eq:defE} allows us to control the dissipation of kinetic energy $\displaystyle \mathscr{D}_\re(f,f)$ defined in~\eqref{eq:Dre}.

We here give an estimate  {(that can be obtained straightforwardly from~\eqref{gamma0bis}) and that will be useful in the subsequent section}
on $\bm{\Psi}_{\e,t} = \bm{\Psi}_{\re_{\e,t}}$, see~\eqref{def:Psire} for the definition.
\begin{lem} \label{lem:Psiet}
For any $r \geq 0$, $\e \in (0,1]$ and $t \geq 0$, we have:
\begin{equation}\label{eq:PsiE}
\bm{\Psi}_{\e,t}(r^{2}) \lesssim \l_{\e}^{\g}(t)r^{3+\g} = \e^{2}z^{\g}(t)r^{3+\g}
\end{equation}
where we recall that $\l_\e(\cdot)$ and $z(\cdot)$ are defined in~\eqref{def:xi}-\eqref{def:zt}.
More precisely,  
\begin{equation}\label{eq:PsiEE}
\left|\bm{\Psi}_{\e,t}(r^{2})-\frac{\mathfrak{a}_{0}}{4+\g}\l_{\e}(t)^{\g}\,r^{3+\g}\right| \lesssim \max\left(\l_{\e}(t)^{2\g}r^{3+2\g},\l_{\e}(t)^{\g+\overline{\g}}r^{3+\g+\overline{\g}},\l_{\e}(t)^{\overline{\g}}r^{3+\overline{\g}}\right)
\end{equation}
holds for any $\e >0,t \geq0,$ and $r >0.$
 \end{lem}
\begin{proof} The proof of \eqref{eq:PsiEE} is an easy adaptation of that of \eqref{eq:PsiEL}.\end{proof}

\subsection{Splitting of the equation} 
We first recall some notations:
\begin{itemize}
\item  the linear operator $\LLe$ and the source term~$S_\e$ are defined in~\eqref{def:LetSe},
\item  the linear operators $\G_\e$ and $\G_{\e,t}$ are respectively defined in~\eqref{def:Ge} and~\eqref{def:Get},
\item the projection $\bm{\pi}_0$ (resp. $\mathbf{P}_0$) onto the kernel of $\LL$ (resp. $\G_\e$) is defined in~\eqref{def:pi0} (resp.~\eqref{eq:P0}),
\item the anti-drift coefficient $\xi(\cdot)$ is defined in~\eqref{def:xi},
\item the collision operator $\Q_{\e,t}$ is defined  as $\Q_{\e,t} = \Q_{\re_{\e,t}}$ where $\re_{\e,t} = \re(\l_\e(t) \cdot)$ with $\l_\e(t)$ defined in~\eqref{def:xi}.
\end{itemize}
We adapt the approach of {\cite{bmam}-\cite{ALT}} and decompose the solution~$h=h_{\e}$ of~\eqref{eq:heps}  into
$$h(t,x,v)=\ho(t,x,v)+\hu(t,x,v)$$
where $\ho=\ho_{\e} \in \E$ and $\hu=\hu_{\e} \in \H$ are the solutions to the following system of equations (in order to lighten the notations, in this whole section, we shall omit the dependence on $\e$ for $h$, $\ho$ and $\hu$): 
{\begin{equation}\label{eq:h0}
	\hspace{-.4cm}\left\{\begin{array}{ccl}
	\partial_{t} h^0\!\!\!\!&=&\!\!\!\mathcal{B}_{\e,t} h^0
	+ \, \e^{-2}{(\LLe-\LL) h^1 + \xi(t) \mathrm{a}_{3}\m_{2}(v)\left(\mathbf{Id} -\mathbf{P}_0\right)h^1}
	\\ [10pt]
	&& \!\!\!
	+ \,  \e^{-1}\Big[\Q_{\e,t}(\ho,\ho)+\Q_{\e,t}(\ho,\hu)+\Q_{\e,t}(\hu,\ho)+ \bm{\pi}_{0}\Q_{\e,t}(\hu,\hu)\Big] 
	+ \mathbf{P}_{0}S_\e
	 \,,\\ [10pt]
	h^0(0,x,v)\!\!\!\!&=&\!\!\!h_{\rm in}^\e(x,v) \in \mathcal E
	\end{array}\right.
	\end{equation}
where $\mathrm{a}_{3}$ is a positive constant defined in Lemma~\ref{lem:hypodrift}, the projectors $\bm{\pi}_0$, $\mathbf{P}_0$ in~\eqref{def:pi0}-\eqref{eq:P0}
and 
\begin{equation}\label{eq:h1}
\left\{
\begin{array}{ccl}
\partial_{t} \hu&=& {\G_{\e}\hu} 
+ \xi(t)\left[-\mathrm{div}_v\left(v \hu\right)-\mathrm{a}_{3}\m_{2}\left(\mathbf{Id-P}_{0}\right)\hu\right]  \\[10pt]
&&\phantom{++++++++}+ \e^{-1}\left(\mathbf{Id}-\bm{\pi}_{0}\right)\Q_{\e,t}(\hu,\hu) + \A_\e \ho + \left(\mathbf{Id-P}_{0}\right)S_{\e}
\,,\\[10pt]
\hu(0,x,v)&=&0 \in \H\,.
\end{array}\right.
\end{equation}
Since
$$\G_{\e,t}\hu=\G_{\e}\hu-\xi(t)\mathrm{div}_v \left(v \hu\right) +\e^{-2}\Big[\LLe\hu-\LL\hu\Big]$$
while
$$\G_{\e,t}\ho=\B_{\e,t}\ho+\A_{\e}\ho\,,$$
one checks easily that  $h =\ho +\hu $ satisfies  indeed
\begin{equation}\label{BEh}
\partial_{t}h=\G_{\e,t}h+\e^{-1}\Q_{\e,t}(h,h)\,, \qquad h(0,x,v)=h^{\e}_{\mathrm{\mathrm{in}}}(x,v) \in \mathcal{E}\end{equation}
which is exactly \eqref{eq:heps} complemented with the initial condition $h(0)=h^{\e}_{\mathrm{\mathrm{in}}}$.} We recall that comments about the splitting of~\eqref{eq:heps} into two equations~\eqref{eq:h0}-\eqref{eq:h1} are given in the Introduction. 

Before starting the analysis of equations~\eqref{eq:h0} and~\eqref{eq:h1}, we recall that $h$ satisfies  
\begin{equation}\label{eq:massmomH}
\int_{\T^{3}\times\R^{3}}h(t,x,v) \left(\begin{array}{c}1\\v\end{array}\right)\, \d v\, \d x=\left(\begin{array}{c}0 \\0\end{array}\right)\end{equation}
and, according to \eqref{eq:FinEe}
$$\int_{\T^{3}\times\R^{3}}h(0,x,v)|v|^{2}\,\d v\,\d x
\longrightarrow 0 \quad \text{as} \quad \e \to0\,$$
where we recall $h=h_\e$.
Recalling the definition of $\mathbf{P}_{0}$ in~\eqref{eq:P0}, since the part of the projection related to the dissipation of energy will play a particular role in our analysis, we define
\begin{equation}\label{eq:PP0}
\mathbb{P}_{0}h:=\sum_{i=1}^{4}\left(\int_{\T^{3}\times\R^{3}}h\,\varphi_{i}\,\d v\, \d x\right)\,\varphi_{i}\,\M\,, \quad \Pi_{0}h:=\left(\int_{\T^{3}\times\R^{3}}h\varphi_{5}\,\d v\, \d x\right)\,\varphi_{5}\,\M\,,
\end{equation}
so that $\mathbb{P}_0 + \Pi_0 = \mathbf{P}_0$.
 Recall also that the eigenfunctions $\varphi_{j}$ are such that
$$\int_{\R^{3}}\varphi_{i}(v)\varphi_{j}(v)\M(v)\, \d v=\delta_{i,j}\,, \qquad \forall \, i,j=1,\ldots,5\,,$$
which in particular implies that, in the Hilbert space $\H$ one has $\mathbf{Id}-\mathbf{P}_{0}=\mathbf{P}_{0}^{\perp}$. 

{
The rest of the section is dedicated to the proof of {\em a priori estimates} on $\ho$ and $\hu$. To this end, during the rest of the section, we assume that {$\ho\in \E$, {$\hu \in \H$}} are solutions to~\eqref{eq:h0}-\eqref{eq:h1} and that there exists $\Mo \leq 1$ such that
\begin{equation} \label{def:Mo}
\sup_{t\geq0}\big(\|\ho (t)\|_{\E}  +  \|\hu(t)\|_{{{\H}}} \big) \leq \Mo\,.
\end{equation}
Mention also that the multiplicative constants involved in the forthcoming estimates of this section may depend on $\Mo \leq 1$. We will only mention it when necessary.}
We also prove a technical Lemma that will be used repeatedly in the rest of the analysis:
\begin{lem}\label{lem:IntAlg}
Given $a,b,\alpha >0$, there are  $C >0$ and $\mu_{*} >0$ both depending only on $a,b,\alpha$ such that 
\begin{equation} \label{eq:intzeta}
\int_{0}^{t}\exp\left(-\mu(t-s)\right)(b+as)^{-\alpha}\,\d s \leq \frac{C}{\mu}(b+at)^{-\alpha}\,,
\end{equation}
holds for any $\mu >\mu_{*}$ and any $t \geq0.$\end{lem}
\begin{proof} Given $t >0$, we introduce
$$I_{1}(t):=\int_{0}^{\frac{t}{2}}\exp\left(-\mu(t-s)\right)(b+as)^{-\alpha}\,\d s\,, \qquad I_{2}(t):=\int_{\frac{t}{2}}^{t}\exp\left(-\mu(t-s)\right)(b+as)^{-\alpha}\,\d s\,.$$
First it holds
$$I_{1}(t) \leq  \exp\left(-\frac{\mu}{2}t\right)\int_{0}^{\frac{t}{2}}(b+as)^{-\alpha}\,\d s \leq b^{-\alpha}\frac{t}{2}\,\exp\left(-\frac{\mu}{2}t\right) \leq \frac{2b^{-\alpha}}{e\mu}\,\exp\left(-\frac{\mu}{4}t\right)$$
where we used the fact that $s\exp(-s) \leq e^{-1}$ for all $s  \geq 0$. To obtain the wanted bound, we remark that the mapping $t \geq 0 \mapsto (b+at)^{\alpha}\,\exp\left(-\frac{\mu}{4}t\right)$ is non-increasing if $\mu > \mu_{*}=\frac{4\alpha a}{b}$. This implies that for any $t \geq 0$, $(b+at)^{\alpha}\,\exp\left(-\frac{\mu}{4}t\right) \leq b^\alpha$ so that
$$
I_{1}(t) \leq \frac{2}{e\mu}(b+at)^{-\alpha}\,, \qquad \forall \,t  \geq 0\,, \quad \forall \, \mu > \mu_{*}\,.$$
Now, since the mapping $s \mapsto (b+as)^{-\alpha}$ is non-increasing, it holds
$$I_{2}(t) \leq \left(b+a\frac{t}{2}\right)^{-\alpha}\int_{\frac{t}{2}}^{t}\exp\left(-\mu(t-s)\right)\d s \leq \frac{1}{\mu}\left(b+a\frac{t}{2}\right)^{-\alpha} \leq \frac{2^{\alpha}}{\mu}(b+at)^{-\alpha}\,.$$
Combining the estimates for $I_{1}(t)$ and $I_{2}(t)$ gives the result.\end{proof}

\subsection{Estimating $\ho$}  
For the part of the solution $\ho(t)$ in $\E$, we have the following estimate which is based upon the dissipativity of $\B_{\e,t}$ as established in Proposition \ref{prop:dissip} (see also Remark~\ref{rem:dissip}):
\begin{prop}\label{prop:h0}
Let $\mu_0 \in (0,\nu_{q})$ (see Proposition~\ref{prop:dissip} and Remark~\ref{rem:dissip} for the definition of $\nu_q$). Then, there exists an explicit~$\e_1 \in (0,\e_0)$ (where $\e_0$ has been defined in Proposition~\ref{prop:dissip}) such that for any~$\e \in (0,\e_1)$ and any $t \geq 0$, 
\begin{equation}\label{eq:estimah00}
\|\ho(t)\|_{\E} \lesssim \|h_{\rm in}\|_{\E}\,\exp\left(-\frac{\mu_0}{\e^2}t\right)+ \e^{2}z^\gamma(t) \,,
\end{equation}
where we recall that $z(\cdot)$ is defined in~\eqref{def:zt}. {Moreover, for any non-increasing and differentiable function $\eta : \R^+ \to \R^+$ such that $\eta(0) >0$ and $\ds\int_0^\infty \eta(\tau) z^\gamma(t) \, \d \tau < \infty$, we have 
	\begin{equation}\label{eq:estimah01}
	\frac{1}{\e^2} \int_0^\infty \eta(\tau) \|h^0(\tau)\|_{\E_1} \, \d \tau \lesssim 1 + \|h_{\rm in}\|_\E\,. 
	\end{equation}
In particular, we can deduce that $h^0 \in L^1_{\rm loc}((0,\infty);\E_1)$ with $\|h^0\|_{L^1((0,T);\E_1)} \lesssim \e^2$ for any $T>0$.}
\end{prop}
 \begin{proof}  We look at \eqref{eq:h0} as a semilinear equation governed by $\B_{\e,t}$ and, using Remark~\ref{rem:dissip} (in particular \eqref{eq:remdissip}), we deduce that
 \begin{multline*}
\dfrac{\d}{\d t}\|\ho(t)\|_{\E} \leq -\frac{\nu_q}{\e^2}\|\ho(t)\|_{\E_{1}} + \frac1\e \Big(\|\Q_{\e,t}(\ho(t),\ho(t))\|_{\E} + \|\Q_{\e,t}(\ho(t),\hu(t))\|_{\E}\\
+\|\Q_{\e,t}(\hu(t),\ho(t))\|_{\E}+\big\|\bm{\pi}_{0}\Q_{\e,t}(\hu(t),\hu(t))\big\|_{\E}\Big)\\
+ {\dfrac{1}{\e^{2}}\big\|\LLe\hu(t) - \LL\hu(t)\big\|_{\E}} 
 +\mathrm{a}_{3}\,\xi(t)\|\m_{2}\left(\mathbf{Id-P}_{0}\right)\hu(t)\|_{\E} + \|\mathbf{P}_{0}S_{\e}(t)\|_{\E}\,
\end{multline*}
holds true for any $\e \in (0,\e_0)$ and any $t\geq0.$
Using classical estimates for $\Q_{\e,t}$ and $\Q_1$ (see Proposition~\ref{prop:BIL}), we have
\begin{multline*}
\|\Q_{\e,t}(\ho(t),\ho(t))\|_{\E} + \|\Q_{\e,t}(\ho(t),\hu(t))\|_{\E}\\
+\|\Q_{\e,t}(\hu(t),\ho(t))\|_{\E} \lesssim \Big(\|\ho(t)\|_{\E}+\|\hu(t)\|_{\E_{1}}\Big)
\|\ho(t)\|_{\E_{1}}\,.
\end{multline*}
Using now Corollary~\ref{lem:LeL1} with $\l=\l_{\e}(t)$, we deduce
$$
 \dfrac{1}{\e^{2}}\big\|\LLe\hu(t) - \LL\hu(t)\big\|_{\E} 
\lesssim \e^{-2}\l_{\e}^{\g}(t) \|\hu(t)\|_{\E_{2}} \lesssim z^\gamma(t)\|\hu(t)\|_{\E_{2}}\,,$$
where we used the definition of $\E_{2}$ in~\eqref{def:E1E2}, that of $\l_{\e}(t)$ and $z(t)$ in~\eqref{def:zt}.  From the conservation of mass and momentum, one deduces from \eqref{def:pi0} that
\begin{multline*}\bm{\pi}_{0}\Q_{\e,t}(\hu(t),\hu(t))=\frac{1}{\sqrt{6}\en_{\star}}\varphi_{5}\,\M\,\int_{\R^{3}}\Q_{\e,t}(\hu(t),\hu(t))(w)|w|^{2}\, \d w\\
=-\frac{1}{\sqrt{6}\en_{\star}}\varphi_{5}\,\M\,\int_{\R^{3} \times \R^{3}}\hu(t,x,w)\hu(t,x,\vb)\bm{\Psi}_{\e,t}(|w-\vb|^{2})\, \d w\, \d \vb\end{multline*}
where we recall that $\bm{\Psi}_{\e,t} = \bm{\Psi}_{\re_{\e,t}}$, see~\eqref{def:Psire} for the definition of $\bm{\Psi}_\re$. 
Using Lemma~\ref{lem:Psiet}, 
we deduce easily that 
{
\begin{equation}\label{eq:pi0Q}\begin{split}
\|\bm{\pi}_{0}\Q_{\e,t}(\hu(t),\hu(t))\|_{\E} &\lesssim \e^{2}z^{\g}(t)\left(\int_{\R^{3}}\|\hu(t,\cdot,\vb)\|_{\W^{m,2}_{x}}\langle \vb\rangle^{3+\g}\, \d \vb\right)^{2} \\
&\lesssim \e^{2}z^{\g}(t)\|\hu(t)\|^2_{\E_1}
\end{split}\end{equation}}
where we used that $\W^{m,2}_{x}(\T^{3})$ is a Banach algebra since $m>\frac{3}{2}$. 
Since moreover, we clearly have
$$
\|\m_{2}\left(\mathbf{Id-P}_{0}\right)\hu(t)\|_{\E} \lesssim \|\hu(t)\|_{\E_{2}}\,,
$$ 
using that $\H \hookrightarrow \E_2 \hookrightarrow \E_1$, we are able to conclude that there exists $C>0$ such that 
\begin{equation*}
\dfrac{\d}{\d t}\|\ho(t)\|_{\E} \leq 
-\frac{1}{\e^2}\Big(\nu_{q}-\e\,C\big(\|\ho(t)\|_{\E}+\|\hu(t)\|_{\H}\big)\Big)\|\ho(t)\|_{\E_{1}}
+ C   z^\gamma(t) { \|\hu(t)\|_{\H}}+ \|\mathbf{P}_{0}S_{\e}(t)\|_{\E}
\end{equation*}
{where we used that from~\eqref{def:xi}-\eqref{def:zt}, $\xi(t) =\mathfrak{a}_{1} z^\gamma(t)$ as well as the fact that $\|\hu(t)\|_\H \leq \Delta_0 \lesssim 1$}. 
For $\mu_0\in(0,\nu_{q})$, we pick {$\e_1 \in (0,\e_0)$ such that $\nu_{q}-\e_1\,C\,\Mo \geq \mu_0$}.  
Consequently, 
we obtain that for any $\e \in (0,\e_1)$ and any $t \geq 0$,
\begin{equation}
\begin{split}\label{imp:h0}
\dfrac{\d}{\d t}\|\ho(t)\|_{\E} 
&\leq -\frac{\mu_0}{\e^2}\,\|\ho(t)\|_{\E_{1}} + Cz^\gamma(t) { \|\hu(t)\|_{\H}} +\|\mathbf{P}_{0}S_{\e}(t)\|_{\E}\,.
\end{split}
\end{equation}
Now, we estimate $\mathbf{P}_{0}S_{\e}$ thanks to Lemma \ref{lem:Ener}. Namely, recalling  \eqref{eq:conservSeps} and the definition of $\Pi_0$ in~\eqref{eq:PP0}, one has
$$\|\mathbf{P}_{0}S_{\e}(t)\|_{\E}=\|\Pi_{0}S_{\e}(t)\|_{\E}=\frac{1}{\sqrt{6}\en_{\star}}\|\varphi_{5}\M\|_{\E}\,\left|\int_{\R^{3}}S_{\e}(t,v)|v|^{2}\,\d v\right|$$
and, with the definition of $S_{\e}(t,v)$ in~\eqref{def:LetSe} and $\xi(t)$ in~\eqref{def:xi}, we deduce from Lemma \ref{lem:Ener} applied with $\l=\l_{\e}(t)=\e^{\frac{2}{\g}}z(t)$,
$$\left|\int_{\R^{3}}S_{\e}(t,v) |v|^{2}\,\d v\right|\lesssim\max\left(\e^{2\frac{\overline{\g}}{\g}-3}z^{\overline{\g}}(t),\e\,z^{2\g}(t),\e^{\frac{2\overline{\g}}{\g}-1}z^{\overline{\g}+\g}(t)\right).$$
Since  
$z(\cdot)$ is decreasing in time (see~\eqref{def:zt}), we deduce that
\begin{equation}\label{eq:Pi0Seps}
{\|\mathbf{P}_{0}S_{\e}(t)\|_{\E} 
\lesssim\e^{\theta} z^{\overline \theta}(t)\,, \qquad \theta:=\min\left(2\frac{\overline{\g}}{\g}-3,1\right) \geq 0\,, \quad \overline \theta := \min(\overline\g, 2\g) \geq \frac{3\g}{2}\,.}
\end{equation}
{Inserting this into \eqref{imp:h0} and using that $\gamma \leq \frac{3\g}{2} \leq \overline\theta$  {with $\|\hu(t)\|_{\H} \leq \Delta_0 \lesssim 1$,} we deduce that
\begin{equation} \label{eq:diffh0}
\dfrac{\d}{\d t}\|\ho(t)\|_{\E} 
\leq -\frac{\mu_0}{\e^2}\,\|\ho(t)\|_{\E_{1}} + Cz^\gamma(t)\,.
\end{equation}
On the one hand,~\eqref{eq:diffh0} implies that 
$$\|\ho(t)\|_{\E} \lesssim \|h_{\rm in}^\e\|_{\E}\exp\left(-\frac{\mu_{0}}{\e^{2}}t\right) + \int_{0}^{t}\exp\left(-\frac{\mu_{0}}{\e^{2}}(t-s)\right)z^{ \g}(s)\,\d s\,.$$
Using the explicit expression of $z(\cdot)$ given in~\eqref{def:zt} together with Lemma~\ref{lem:IntAlg} gives~\eqref{eq:estimah00}.
On the other hand, multiplying~\eqref{eq:diffh0} by $\eta(t)$ and using that $\eta'(t) \leq 0$ for any $t \geq 0$, we obtain that 
	$$
	\frac{\d}{\d t} \left(\eta(t) \|h^0(t)\|_\E\right) \leq -\eta(t) \frac{\mu_0}{\e^2}\,\|\ho(t)\|_{\E_{1}} + C\eta(t) z^\gamma(t)\,.
	$$
Integrating this differential inequality in time and using the assumptions made on $\eta$ yields the wanted inequality~\eqref{eq:estimah01}.}
\end{proof}

\subsection{Estimating $\mathbf{P}_{0}\hu$}
{
Let us point out that getting estimates on $\hu$ is trickier than in~\cite{bmam} since the energy is no longer preserved which induces additional difficulties. We use a strategy similar (though more involved) to the one used in \cite{ALT} and start with a basic estimate on ~$\mathbb{P}_{0} h^1$ where $\mathbb{P}_0$ is defined in~\eqref{eq:PP0}. One has
\begin{equation} \label{eq:estimP0h1}
	\|\mathbb{P}_0 \hu(t)\|_\H \lesssim \| \mathbb{P}_{0} \hu(t)\|_{\E} \lesssim \|\ho(t)\|_\E\,.
	\end{equation}
 The proof of this estimate is straightforward using the explicit expression of $\mathbb{P}_0$ given in~\eqref{eq:PP0} which implies that $\mathbb{P}_0 \in \mathscr{B}(\E,\H)$. We also use the additional property that $\mathbb{P}_{0}S_{\e}=0$ (see \eqref{eq:conservSeps}) from which $\mathbb{P}_{0}h(t)=0$. The second key observation regards the action of $\Pi_{0}$ on the linearized operator $\G_{\e,t}$ (defined in~\eqref{def:Get}). 
\begin{lem} \label{lem:Pi0}
There exists $\e_2 \in (0,\e_1)$ (where $\e_1$ has been defined in Proposition~\ref{prop:h0}) such that for any $\e \in (0,\e_{2})$ , there is a \emph{non-increasing} mapping $\lambda_\e\::\:\R^{+} \to \R$ and another one $r_\e\::\:\R^{+}\to \R^{+}$ such that
\begin{equation} \label{pi0Ge}
\Pi_0 [\G_{\e,t} h(t)] = - \lambda_\e(t) \Pi_0 h(t) + r_\e(t)\varphi_{5} \mathcal M \end{equation}
where we recall that $h=h(t)$ is a solution to~\eqref{eq:heps} and $\varphi_{5}$ is defined in \eqref{def:Psii}. Moreover, $\lambda_\e(t) >0$ and  $r_\e(t) \in \R$ are such that for any~$t \geq 0$ and any $\e \in (0,\e_2)$,
\begin{equation}\label{eq:stimerese}
z^\g(t) \lesssim \lambda_\e(t) \lesssim z^{\g}(t)\,, \qquad \left|{r}_\e(t)\right| \lesssim \lae(t)\left\|\left(\mathbf{Id-P}_{0}\right)h(t)\right\|_{L^{1}_{x,v}({\m_{3+\g}})}\,.
\end{equation}
\end{lem}
\begin{proof} 
\noindent {\it Step 1.} One first recalls that, due to conservation of mass,  if $h=h(t)$ is a solution to \eqref{eq:heps} then
\begin{equation}\label{eq:Pi0H}
\Pi_{0}h(t)=\frac{1}{\en_\star \sqrt{6}}\left(\int_{\T^{3}\times\R^{3}}h\,|v|^{2}\, \d v\, \d x\right)\varphi_{5} \M\,.\end{equation}
We check then easily that
\begin{equation*}\begin{split}
\Pi_{0}\left[\G_{\e,t}h(t)\right]&=\frac{1}{\sqrt{6}\en_{\star}}\left[ \int_{\T^{3}\times\R^{3}}\left[\e^{-2}\LLe h(t,x,v)-\xi(t)\mathrm{div}_{v}(vh(t,x,v)\right]\,|v|^{2}\,\d v\,\d x\right]\varphi_{5}\M\\
&=\frac{1}{\sqrt{6}\en_{\star}}\left[\e^{-2}\int_{\T^{3}\times\R^{3}} \LLe h(t,x,v)|v|^{2}\,\d v\,\d x\right.\\
&\phantom{+++++} \left.+2\xi(t)\int_{\T^{3}\times \R^{3}}h(t,x,v)\,|v|^{2}\,\d v \,\d x\right]\varphi_{5}\M\,.\end{split}\end{equation*}
Concerning the first term, we have that
\begin{equation*}\begin{split}\int_{\T^{3}\times\R^{3}} \LLe h(t,x,v)&|v|^{2}\,\d v\,\d x=-2\int_{\T^{3}}\d x\int_{\R^{3}\times \R^{3}}h(t,x,v)\M(\vb)\bm{\Psi}_{\e,t}(|v-\vet|^{2})\,\d \vb\,\d v\\
&=-2\int_{\T^{3}}\d x\int_{\R^{3}\times \R^{3}}\Pi_{0}h(t,x,v)\M(\vb)\bm{\Psi}_{\e,t}(|v-\vet|^{2})\,\d \vb\,\d v\\
&\phantom{++} -2\int_{\T^{3}}\d x\int_{\R^{3}\times \R^{3}}\left(\mathbf{Id}-\Pi_{0}\right)h(t,x,v)\M(\vb)\bm{\Psi}_{\e,t}(|v-\vet|^{2})\,\d \vb\,\d v\,.
\end{split}\end{equation*}
Observing that $(\mathbf{Id}-\Pi_{0})h=\mathbf{(Id-P_{0})}h$ we set
$$r_{\e}(t):=-\frac{2}{\e^{2}\en_{\star}\sqrt{6}}\int_{\R^{3}\times \R^{3}}\left(\mathbf{Id-P}_{0}\right)h(t,x,v)\M(\vb)\bm{\Psi}_{\e,t}(|v-\vet|^{2})\,\d \vb\,\d v$$
and, using \eqref{eq:Pi0H}, we obtain easily that
$$\Pi_{0}\left[\G_{\e,t}h(t)\right]=-\lambda_{\e}(t)\Pi_{0}h(t) + r_{\e}(t)\varphi_{5}\M\,$$
where we define
\begin{equation}\label{eq:LambdaE}
\lambda_{\e}(t):=\frac{2}{\e^{2}}\left[\frac{1}{\en_\star \sqrt{6}}\int_{\R^{3}\times \R^{3}}\varphi_{5}(v)\M(v)\M(\vb)\bm{\Psi}_{\e,t}(|v-\vet|^{2})\,\d \vb \, \d v  
-\e^{2}\xi(t)\right]\,,\end{equation}
which is indeed independent of $h$. 

\noindent {\it Step 2.} We are now going to prove that, for $\e$ small enough, 
\begin{equation} \label{eq:zlambda}
z^\g(t) \lesssim \lambda_{\e}(t) \lesssim z^{\g}(t)\,.
\end{equation} We first prove the bound from above. 
As in the proof of Lemma \ref{lem:Ener}, our strategy consists in estimating, for $\l >0$ and $\re_{\l}(r)=\re(\l r)$, $r >0$ the quantity
\begin{equation*}
I_{\l}:=\int_{\R^{3}\times \R^{3}}\left(|v|^{2}-3\en_{\star}\right)\M(v)\M(\vet)\bm{\Psi}_{\re_{\l}}(|v-\vet|^{2})\,\d \vb\,\d v - 6\en_{\star}^{2}\mathfrak{a}_{1}\l^{\g}=:\l^{\g}J_{\l}-6\en_{\star}^{2}\mathfrak{a}_{1}\l^{\g}\end{equation*}
where we observe that, for $\l=\l_{\e}(t)$ and since $\e^{2}\xi(t)=\mathfrak{a}_{1}\l^{\g}_{\e}(t)$, one has
$$\lambda_{\e}(t)=\frac{1}{3\vartheta^{2}_{\star}\e^{2}}J_{\l_{\e}(t)}.$$
Using again $\bm{\Psi}_{\re_{\l}}(r^{2})=\l^{-3}\bm{\Psi}_{\re}(\l^{2}r^{2})$, one has
$$J_{\l}=\frac{1}{\l^{3+\g}}\int_{\R^{3}\times \R^{3}}\left(|v|^{2}-3\en_{\star}\right)\M(v)\M(\vet)\bm{\Psi}_{\re_{\l}}(|v-\vet|^{2})\,\d \vb\,\d v=:J^{1}_{\l}-J_{\l}^{2}\,.$$
We compute $J_{\l}^{2}$ as in the proof of Lemma \ref{lem:Ener} and recalling that $\M(v)=\en_{\star}^{-\frac{3}{2}}M(\en_{\star}^{-\frac{1}{2}}v)$ where~$M$ is defined in~\eqref{def:M}, we get
$$J_{\l}^{2}=\frac{3\en_{\star}}{\l^{3+\g}}\int_{\R^{3}\times \R^{3}}\M(v)\M(\vet)\bm{\Psi}_{\re_{\l}}(|v-\vet|^{2})\,\d \vet \,\d v=\frac{3\en_{\star}}{\l^{3+\g}}\int_{\R^{3}}M(y)\bm{\Psi}_{\re}(2\l^{2}\en_{\star}^{2}|y|^{2})\,\d y\,.$$
In the same way,
$$J_{\l}^{1}=\frac{\en_{\star}}{\l^{3+\g}}\int_{\R^{3}\times \R^{3}} |v|^{2} M(v)M(\vet)\bm{\Psi}_{\re}(\l^{2}\en_{\star}|v-\vet|^{2})\,\d \vb \, \d v\,.$$
As in \cite[Lemma A.1]{MiMo3}, we can then write
$$J_{\l}^{1}=\frac{\en_{\star}}{4\l^{3+\g}}\int_{\R^{3}\times \R^{3}}\left( |v+\vet|^{2}+|v-\vet|^{2}\right) M(v)M(\vet)\bm{\Psi}_{\re}(\l^{2}\en_{\star}|v-\vet|^{2})\,\d \vb\,\d v$$
and, using the unitary change of variable $y=\frac{v+\vet}{\sqrt{2}},$ $y_*=\frac{v-\vet}{\sqrt{2}}$, since $M(v)M(\vet)=M(y)M(y_*)$, we deduce that
$$J_{\l}^{1}=\frac{\en_{\star}}{2\l^{3+\g}}\int_{\R^{3}\times \R^{3}}\left(|y|^{2}+|y_*|^{2}\right)M(y)M(y_*)\bm{\Psi}_{\re}(2\l^{2}\en_{\star}|y_*|^{2})\,\d y\,\d y_*\,.$$
Thus, since $\int_{\R^{3}}|y|^{2}M(y)\,\d y=3$ and $\int_{\R^{3}}M(y)\,\d x=y$, we deduce that
$$J_{\l}^{1}=\frac{\en_{\star}}{2\l^{3+\g}}\left(3\int_{\R^{3}} M(y)\bm{\Psi}_{\re}(2\l^{2}\en_{\star}|y|^{2})\,\d y+\int_{\R^{3}}|y|^{2}M(y)\bm{\Psi}_{\re}(2\l^{2}\en_{\star}|y|^{2})\,\d y\right).$$
Combining these identities, we get
$$J_{\l}=\frac{\en_{\star}}{2\l^{3+\g}}\int_{\R^{3}}\left(|y|^{2}-3\right)M(y)\bm{\Psi}_{\re}(2\l^{2}\en_{\star}|y|^{2})\,\d y\,.$$
As in the proof of Lemma \ref{lem:Ener}, we use that
$$\left|\frac{1}{\l^{3+\g}}\bm{\Psi}_{\re}(\l^{2}r^{2})-\frac{\mathfrak{a}_{0}}{4+\g}r^{3+\g}\right| \lesssim \l^{\overline{\g}-\g}r^{ {3+\overline\g}}+\l^{\g}r^{3+2\g}+\l^{\overline{\g}}r^{3+\g+\overline{\g}}$$
to deduce that
$$\left|J_{\l}-\frac{\mathfrak{a}_{0}\en_{\star}}{2(4+\g)}(2\en_{\star})^{\frac{3+\g}{2}}(K_{2+\gamma}-3K_\gamma)\right| \lesssim \max\left(\l^{\overline{\g}-\g},\l^{\g},\l^{\overline{\g}}\right)
$$
where $K_s$ has been defined in~\eqref{def:Kgamma}. 
We thus deduce that there exists $C>0$ such that
$$J_{\l} \geq \frac{\mathfrak{a}_{0}\en_{\star}}{2(4+\g)}(2\en_{\star})^{\frac{3+\g}{2}}(K_{2+\gamma}-3K_\gamma)  -C\max\left(\l^{\overline{\g}-\g},\l^{\g},\l^{\overline{\g}}\right)\,.$$
After multiplying it with $\l^{\g}$ we deduce that
\begin{align*}
&I_{\l}=\l^{\g}J_{\l}-6\en_{\star}^{2}\mathfrak{a}_{1}\l^{\g}  \\
&\quad
\geq \left(\frac{\mathfrak{a}_{0}\en_{\star}}{2(4+\g)}(2\en_{\star})^{\frac{3+\g}{2}}(K_{2+\gamma}-3K_\gamma)-6\en_{\star}^{2}\mathfrak{a}_{1}\right)\l^{\g}-C\max\left(\l^{\overline{\g}}\,,\l^{2\g}\,,\l^{\overline{\g}+\g}\right)\\
&\quad
=: K \ell^\gamma - C\max\left(\l^{\overline{\g}}\,,\l^{2\g}\,,\l^{\overline{\g}+\g}\right)\,.
\end{align*}
We observe that, with \eqref{def:a1}, 
$$
K=\frac{\mathfrak{a}_{0}\en_{\star}}{2(4+\g)}(2\en_{\star})^{\frac{3+\g}{2}}\left(K_{2+\g}-5K_{\g}\right)\,,
$$
where $K_{\g+2}-5K_{\g} >0$ according to Lemma \ref{lem:eqK}. Therefore, using that $\overline{\g} >\g$, there is $\l_{2} >0$ small enough such that
$$J_{\l} \geq \frac{K}{2}\l^{\g}, \qquad \forall \, \l \in (0,\l_{2})\,.$$
Applying this to $\l=\l_{\e}(t)$ and choosing $\e_{2}$ small enough so that $\l_{\e}(t)<\l_{2}$ for any $\e \in (0,\e_{2})$ and any $t\geq0$, we deduce  
$$\lambda_{\e}(t)=\frac{1}{3\en_{\star}^{2}\e^{2}}J_{\l_{\e}(t)} \geq \frac{\mathfrak{a}_2}{2} \frac{1}{\e^2} \l_{\e}^{\g}(t)\geq \frac{\mathfrak{a}_2}{2} z^{\g}(t)$$
since $\l_{\e}(t)=\e^{\frac{2}{\g}}z(t)$ and where 
\begin{equation} \label{def:a2}
\mathfrak{a}_2 := \frac{1}{3\en_{\star}^{2}} K = \frac{1}{3\en_{\star}^{2}} \frac{\mathfrak{a}_{0}\en_{\star}}{2(4+\g)}(2\en_{\star})^{\frac{3+\g}{2}}\left(K_{2+\g}-5K_{\g}\right)\,.
\end{equation} 
Using again \eqref{eq:PsiE}, we also deduce the upper bound $\lambda_\e(t) \lesssim z^\gamma(t)$ simply neglecting the nonpositive term $-\e^{2}\xi(t)$ in~\eqref{eq:LambdaE}. The bound on $r_{\e}(t)$ is also obvious using \eqref{eq:PsiE} and~\eqref{eq:zlambda}. 
\end{proof}

\begin{nb} \label{rem:a2}
We notice that we can be more precise in~\eqref{eq:zlambda}. Indeed, if $\delta \in (0,1)$ is given, one can choose $\e_2 = \e_2(\delta)$ small enough so that for any $\e \in (0,\e_2)$ and any $t \geq 0$, 
\begin{equation} \label{eq:stimerese2}
(1-\delta) \mathfrak{a}_2 \leq \lae(t) z(t)^{-\g} \leq (1+\delta) \mathfrak{a}_2
\end{equation}
with $\mathfrak{a}_2$ defined in~\eqref{def:a2}.
Recalling that $\mathfrak{a}_1$ is defined in~\eqref{def:a1}, it is also worth noticing that using the proof of Lemma~\ref{lem:eqK}, we have
\begin{equation} \label{eq:a2a1}
\frac{\mathfrak{a}_2}{\mathfrak{a}_1} = \frac{K_{\g+2}-5K_\g}{K_\g} = \g+1\,. 
\end{equation}
\end{nb}
The quantity $\lae(t)$ defined in~\eqref{eq:LambdaE} shall play an important role in what follows. We thus introduce a notation for its supremum in $\e$ and in time: 
\begin{equation} \label{def:lambda*}
\lambda_\star:= \sup_{\e \in (0,1)} \sup_{t \geq 0} \lae(t)\,.
\end{equation} 
We are now able to derive a nice estimate on $\mathbf{P}_0 h^1$. 
\begin{lem}\label{lem:energy}
Consider $\delta \in \left(0,\frac{1}{\g+1}\right)$. For any $\e \in (0,\e_{2})$ (where $\e_2=\e_2(\delta)$ has been defined in Remark~\ref{rem:a2}) and any $t \geq 0$,
\begin{equation}\label{eq:lemenergy} 
\begin{aligned}
&\left\|\mathbf{P}_{0}\hu(t)\right\|_{\H} \\
&\quad
\lesssim \|\ho(t)\|_{\E} + \|h_{\rm in}\|_{ \E}\,\exp\left(-\int_{0}^{t}\lae(\tau)\,\d\tau\right) 
+\e^\theta \exp\left(- (1-\delta) \frac{\overline \theta-\g}{\g+1} \int_0^t \lae(\tau) \, \d\tau \right) \\
&\qquad +\int_{0}^{t}\lae(s)\exp\left(-\int_{s}^{t}\lae(\tau)\,\d \tau\right)\,\left(\|\ho(s)\|_{\E}   + \left\|\left(\mathbf{Id-P}_{0}\right)\hu(s)\right\|_{\H}\right)\d s\,
\end{aligned}
\end{equation} 
where $\theta$ and $\overline \theta$ have been defined in~\eqref{eq:Pi0Seps}.  
\end{lem}
\begin{nb}
We remark that in the case of true viscoelastic hard spheres, we have $\overline \g = 2\g = 2/5$ so that the rate given by the second integral is 
$\exp \left(-\frac16 (1-\delta) \int_0^t \lae(\tau) \, \d\tau\right)$.
\end{nb}
\begin{proof} Since $\mathbf{P}_{0}\hu=\mathbb{P}_0\hu+\Pi_{0}\hu$ and 
$\Pi_0 \hu = \Pi_0 h- \mathbf{P}_0 \ho - \mathbb{P}_0 \hu$, using \eqref{eq:estimP0h1} and the fact that $\mathbf{P}_{0} \in \mathscr{B}(\E,\H)$, one sees that
\begin{equation}\label{eq:PoPi0}
\|\mathbf{P}_{0}\hu(t)\|_{\H} \lesssim \|\Pi_{0}h(t)\|_{\H}+\|\ho(t)\|_{\E}\,.\end{equation}
We thus only need to estimate $\|\Pi_{0}h\|_{\H}$. Recalling that $h$ solves
$\partial_{t} h= \G_{\e,t}h + \frac1\e\Q_{\e,t}(h,h)+S_{\e}$,
we deduce from \eqref{pi0Ge} that
{
\begin{equation*}
\partial_{t}\big(\Pi_{0} h(t)\big)=-\lae(t)\Pi_{0}h(t) + r_{\e}(t)\varphi_{5}\M+ \frac1\e\Pi_{0} \Q_{\e,t}(h(t),h(t))+\Pi_{0}S_{\e}(t)\,,\end{equation*}}
so that
\begin{multline}\label{e1-pert-en}
\Pi_{0}  h(t)= \Pi_{0}  h_{\rm in}\,\exp\left(-\int_{0}^{t}\lae(\tau)\,\d\tau\right) \\
+ \int^{t}_{0}\exp\left(-\int_{s}^{t}\lae(\tau)\,\d \tau\right)\left(\frac1\e\Pi_{0} \Q_{\e,s}(h(s),h(s)) + r_{\e}(s)\varphi_{5}\M+\Pi_{0}S_{\e}(s)\right) \, \d s\,.
\end{multline}
For any $s \geq0,$ one has
\begin{multline*}
\left\|\Pi_{0}\Q_{\e,s}(h(s),h(s))\right\|_{\H}
=\frac{1}{\sqrt{6}\en_{\star}}\|\varphi_{5}\M\|_{\H}\\\left|\int_{\T^{3}}\d x\int_{\R^{3}\times \R^{3}}h(s,x,v)h(s,x,v)\bm{\Psi}_{\e,s}(|v-\vet|^{2})\,\d \vb\,\d v \right|\,.
\end{multline*}
Therefore, arguing as in \eqref{eq:pi0Q},
\begin{multline*}\left\|\Pi_{0}\Q_{\e,s}(h(s),h(s))\right\|_{\H}
\lesssim \e^{2}z^{\g}(s)\int_{\T^{3}}\d x\left[\int_{\R^{3}}|h(s,x,v)|\m_{3+\g}(v)\, \d v\right]^{2}\\
\lesssim \e^{2}z^{\g}(s)\left[\int_{\R^{3}}\m_{3+\g}(v)\left(\int_{\T^{3}}|h(s,x,v)|^{2}\, \d x\right)^{\frac{1}{2}}\, \d v\right]^{2}
\end{multline*}
where we used  Minkowski's integral inequality in the last estimate. We deduce then that \begin{equation}\label{eq:Pi0Qe}
\left\|\Pi_{0}\Q_{\e,s}(h(s),h(s))\right\|_{\H} \lesssim
\e^{2}z^{\g}(s)\|h(s)\|_{L^{1}_{v}L^{2}_{x}(\m_{3+\g})}^{2} \lesssim \e^{2} z^{\g}(s)\|h(s)\|_{\E}^{2}\end{equation}
where we used  that $q \geq 3+\g$. \smallskip
\noindent 
Now,~\eqref{eq:Pi0Seps} directly implies that
\begin{equation}\label{eq:Pi0Seps2}\|\Pi_{0}S_{\e}(s)\|_{\H} \lesssim \e^{\theta}z^{\overline\theta}(s)\,. 
\end{equation}
We then take the $\|\cdot\|_{\H}$-norm of \eqref{e1-pert-en}, we use \eqref{eq:Pi0Qe} - \eqref{eq:Pi0Seps2} as well as the bound on $r_{\e}$ in~\eqref{eq:stimerese}. Using also the first estimate of~\eqref{eq:stimerese}, we  obtain
\begin{multline}\label{eq:Pi0h}
\big\| \Pi_{0}  h(t) \big\|_{ \H} \lesssim \| \Pi_{0}  h_{\rm in}\|_{ \H}\,\exp\left(-\int_{0}^{t}\lae(\tau)\,\d\tau\right) + \\
+\int_{0}^{t}\exp\left(-\int_{s}^{t}\lae(\tau)\,\d \tau\right)\lae(s)\,\left(\e\|h(s)\|^{2}_{\E} + \e^{\theta} \lae^{\frac{\overline\theta}{\gamma}-1}(s)+ \left\|\left(\mathbf{Id-P}_{0}\right)h(s)\right\|_{{L^{1}_{x,v}(\m_{3})}}\right)\,\d s\,,
\end{multline}
where we notice that $\frac{\overline\theta}{\g} - 1 \geq \frac{1}{2}$. {We are now going to prove that  
\begin{equation} \label{eq:estimint}
\int_0^t \exp\left(-\int_s^t \lae(\tau) \, \d\tau \right) \lae^{\frac{\overline \theta}{\gamma}}(s) \, \d s
\lesssim \exp\left(-(1-\delta) \frac{\overline \theta - \g}{\g+1} \int_0^t \lae(\tau) \, \d\tau \right) \,, \quad \forall \, t \geq 0 \,.
\end{equation}
For this, we recall that from~\eqref{eq:stimerese2} and~\eqref{def:xi}
$$
\lae(\tau) \geq \frac{\mathfrak{a}_2(1-\delta)}{\mathfrak{b}_{1}^{\frac{\g}{\g+1}}+\g\mathfrak{a}_{1}\tau}\,, \quad \forall \, \tau \geq 0\,,
$$
from which we deduce 
$$
\int_s^t \lae(\tau) \, \d \tau \geq (1-\delta) \frac{\g+1}{\g} \left[\log\left(\mathfrak{b}_1^\frac{\g}{\g+1} + \g \mathfrak{a}_1 t \right) - \log\left(\mathfrak{b}_1^\frac{\g}{\g+1} + \g \mathfrak{a}_1 s \right) \right]\,, \quad \forall \, 0 \leq s \leq t \,,
$$
where we used the fact that $\mathfrak{a}_2=\mathfrak{a}_1(\g+1)$ from~\eqref{eq:a2a1}. 
Using again~\eqref{eq:stimerese2} and~\eqref{def:xi}, we obtain 
\begin{align*}
&\int_0^t \exp\left(-\int_s^t \lae(\tau) \, \d\tau \right) \lae^{\frac{\overline \theta}{\gamma}}(s) \, \d s \\
&\quad 
\lesssim \left(\mathfrak{b}_1^\frac{\g}{\g+1} + \g \mathfrak{a}_1 t \right)^{-(1-\delta)\frac{\g+1}{\g}} 
\int_0^t \left(\mathfrak{b}_1^\frac{\g}{\g+1} + \g \mathfrak{a}_1 s \right)^{-\frac{\overline \theta}{\gamma}+(1-\delta)\frac{\g+1}{\g}} \, \d s\,.
\end{align*}
Noticing now that $\overline \theta = \min(\overline\g,2\g) \leq 2\g$, we deduce that
$$
-\frac{\overline \theta}{\g} +1+ (1-\delta) \frac{\g+1}{\g}
\geq -1 + (1-\delta) \frac{\g+1}{\g} >0
$$
since $\delta<\frac{1}{\g+1}$. As a consequence, we obtain that 
$$
\int_0^t \left(\mathfrak{b}_1^\frac{\g}{\g+1} + \g \mathfrak{a}_1 s \right)^{-\frac{\overline \theta}{\gamma}+(1-\delta)\frac{\g+1}{\g}} \, \d s
\lesssim 
\left(\mathfrak{b}_1^\frac{\g}{\g+1} + \g \mathfrak{a}_1 t \right)^{-\frac{\overline \theta}{\gamma}+1+(1-\delta)\frac{\g+1}{\g}}
$$
so that 
\begin{equation*}
\int_0^t \exp\left(-\int_s^t \lae(\tau) \, \d\tau \right) \lae(s)^{\frac{\overline \theta}{\gamma}} \, \d s  
\lesssim \left(\mathfrak{b}_1^\frac{\g}{\g+1} + \g \mathfrak{a}_1 t \right)^{-\frac{\overline \theta}{\gamma}+1} \,.
\end{equation*}
Finally, using ~\eqref{eq:stimerese2} and~\eqref{def:xi} together with the fact that $\overline \theta -\g \geq \frac{\g}{2} >0$, we remark that
\begin{equation*}\begin{split}
\bigg(\mathfrak{b}_1^\frac{\g}{\g+1}  &+ \g \mathfrak{a}_1 t \bigg)^{-\frac{\overline \theta}{\gamma}+1}
 = \exp \left( -\left(\overline \theta - \gamma\right) \mathfrak{a}_1 \int_0^t \left(\mathfrak{b}_1^\frac{\g}{\g+1} + \g \mathfrak{a}_1 s\right)^{-1} \, \d s \right) \\
 &\lesssim \exp \left( -\left(\overline \theta - \gamma\right)\frac{\mathfrak{a}_1}{(1+\delta) \mathfrak{a}_2} \int_0^t \lae(s) \, \d s\right)=\exp \left( - \frac{\overline \theta - \g}{(1+\delta) (\g+1)} \int_0^t \lae(s) \, \d s\right) \\
 &\lesssim \exp \left(-\left(1-\delta\right)  \frac{\overline \theta - \g}{\g+1}  \int_0^t \lae(s) \, \d s\right)
\end{split}\end{equation*}
where we used  the fact that $(1+\delta)^{-1} \geq 1-\delta$ for $\delta \in (0,1)$ in the last estimate. 
}
From this, we deduce that 
\begin{multline}\label{eq:lemenergy1} 
\big\| \Pi_{0}  h(t) \big\|_{ \H} \lesssim \|h_{\rm in}\|_{ \E}\,\exp\left(-\int_{0}^{t}\lae(\tau)\,\d\tau\right) +\e^\theta\exp\left(-\left(1-\delta\right)  \frac{\overline \theta - \g}{\g+1}\int_{0}^{t}\lae(\tau)\,\d\tau\right)\\
+\int_{0}^{t}\lae(s)\exp\left(-\int_{s}^{t}\lae(\tau)\d \tau\right)\,\left(\e\|h(s)\|_{\E}^{2}  + \left\|\left(\mathbf{Id-P}_{0}\right)h(s)\right\|_{{L^{1}_{x,v}(\m_{3})}}\right)\d s\,.
\end{multline}
Observing  that 
$$\left\|\left(\mathbf{Id-P}_{0}\right)h(s)\right\|_{{L^{1}_{x,v}(\m_{3})}} \lesssim \left\|\left(\mathbf{Id-P}_{0}\right)\hu(s)\right\|_{\H}+\|\ho(s)\|_{\E}$$
and using~\eqref{def:Mo},
$$\|h(s)\|_{\E}^{2} \lesssim  \|\hu(s)\|_{\H}^2 + \|\ho(s)\|_{\E}^2 \lesssim \Delta_0^2 \,,$$ we deduce then \eqref{eq:lemenergy} from \eqref{eq:lemenergy1} and 
\eqref{eq:PoPi0} and the fact that $\e \leq \e^\theta$ since $\theta \leq 1$.  \end{proof}

Inserting the estimate on $\|\ho(t)\|_{\E}$ provided by \eqref{eq:estimah00} allows us to derive a more precise estimate of~$\mathbf{P}_{0}\hu(t)$ in the following proposition:
\begin{prop}
\label{prop:precisePo}
Consider $\delta \in \left(0,\frac{1}{\g+1}\right)$. There exists an explicit $\e_{3} \in (0,\e_{2})$ (where $\e_2$ has been defined in Remark~\ref{rem:a2}) such that for any $\e \in (0,\e_{3})$ and any $t \geq 0$, there holds
\begin{multline}\label{eq:P0h11}
\|\mathbf{P}_{0}\hu(t)\|_{\H} \lesssim \|h_{\rm in}\|_{\E}\exp\left(-\int_{0}^{t}\lae(\tau)\,\d \tau\right) + \e^\theta \exp\left(-(1-\delta)\frac{\overline \theta - \g}{\g+1}\int_{0}^{t}\lae(\tau)\,\d\tau\right)\\ + \int_{0}^{t}\lae(s)\exp\left(-\int_{s}^{t}\lae(\tau)\,\d \tau\right)\left\|\left(\mathbf{Id-P}_{0}\right)\hu(s)\right\|_{\H} \d s
\end{multline}
where $\theta$ and $\overline \theta$ are defined in~\eqref{eq:Pi0Seps}. 
\end{prop}
\begin{proof} 
We insert in \eqref{eq:lemenergy} the bound for $\|\ho(t)\|_{\E}$ obtained in \eqref{eq:estimah00}. First, we remark that using~\eqref{eq:stimerese2} and~\eqref{def:xi}, we have
	\begin{equation} \label{eq:zgexp}
	z^\gamma(t) \lesssim \exp\left(- (1-\delta) \frac{\g}{\g+1} \int_0^t \lae(\tau) \, \d\tau\right)
	\end{equation}
so that
	\[
	\|\ho(t)\|_\E \lesssim \|h_{\rm in}\|_\E \exp\left(-\frac{\mu_0}{\e^2}t\right) + \e^2 \exp\left(- (1-\delta) \frac{\g}{\g+1} \int_0^t \lae(\tau) \, \d\tau\right)\,. 
	\]
We deduce that
\begin{align*} 
&\| \mathbf{P}_{0} \hu(t)  \|_{\H}  
\lesssim \e^\theta \exp\left(- (1-\delta) \frac{\overline \theta-\g}{\g+1} \int_0^t \lae(\tau) \, \d\tau\right)\\
&+  
\|h_{\rm in}\|_\E \Bigg[\exp\left(-\frac{\mu_0}{\e^2}t\right) + \exp\left(-\!\int_{0}^{t}\lae(\tau)\,\d\tau\right) 
+\int_{0}^{t}\lae(s)\exp\left(-\!\int_{s}^{t}\lae(\tau)\,\d \tau -\frac{\mu_0}{\e^{2}}s\right)\d s\Bigg] \\
&\quad
+\int_{0}^{t}\lae(s)\exp\left(-\int_{s}^{t}\lae(\tau)\,\d \tau\right)\left\|\left(\mathbf{Id-P}_{0}\right)\hu(s)\right\|_{\H} \d s
\end{align*}
where we used the fact that $\e^2 \leq \e^\theta$ since $\theta \leq 1$ as well as $\g\geq \overline\theta-\g$ by definition of $\overline \theta$,
and noticed that 
\begin{align*}
&\int_0^t \lae(s) \exp \left( - \int_s^t \lae(\tau) \, \d\tau \right) \exp\left(- (1-\delta) \frac{\g}{\g+1} \int_0^s\lae(\tau) \, \d\tau\right) \, \d s \\
&\quad \leq \exp\left(-(1-\delta) \frac{\g}{\g+1}  \int_0^t\lae(\tau) \, \d\tau\right) \\
&\qquad \qquad \qquad \qquad 
\times \int_0^t \lae(s) \exp \left( - \left(1-(1-\delta) \frac{\g}{\g+1} \right) \int_s^t \lae(\tau) \, \d\tau \right) \, \d s \\
&\quad \lesssim \exp\left(-(1-\delta) \frac{\g}{\g+1}  \int_0^t\lae(\tau) \, \d\tau\right)
\lesssim \exp\left(-(1-\delta) \frac{\overline \theta - \g}{\g+1}  \int_0^t\lae(\tau) \, \d\tau\right)
\end{align*}
since $\frac{(1-\delta) \g}{\g+1} <1$. 
We can then choose~$\e_{3} \in (0,\e_{2})$ such that $\mu_0 \geq 2\e^{2}\lambda_\star$ for $\e \in (0,\e_{3})$. With this, one has
\begin{multline*}\exp\left(-\frac{\mu_0}{\e^2}t\right) + \exp\left(-\int_{0}^{t}\lae(\tau)\,\d\tau\right) \\
+\int_{0}^{t}\lae(s)\exp\left(-\int_{s}^{t}\lae(\tau)\,\d \tau -\frac{\mu_0}{\e^{2}}s\right)\,\d s  \lesssim \exp\left(-\int_{0}^{t}\lae(\tau)\,\d \tau\right)
\end{multline*}
which gives \eqref{eq:P0h11}. \end{proof} 
\color{black}
We then prove some technical  estimates that shall be useful to prove the following Corollary. \begin{lem}\label{lem:integdouble} Consider $\eta: \R^+ \to \R^+$. For any $\e \in (0,\e_3)$ (where $\e_3$ is defined in Proposition~\ref{prop:precisePo}) and any $t \geq 0$, 
\begin{multline*}
\int_{0}^{t}\lae(s)\exp\left(-\int_{s}^{t}\lae(r)\,\d r\right)\int_{0}^{s}\exp\left(-\frac{\mu_0}{\e^{2}}(s-\tau)\right)\eta(\tau)\,\d \tau \, \d s \\
\lesssim \e^2 \int_{0}^{t}{\lae(s)}\exp\left(-\int_{s}^{t}\lae(\tau)\, \d \tau\right)\eta(s)\,\d s\,.
\end{multline*} 
Moreover, for any $\delta \in (0,1)$, 
\begin{multline}\label{eq:integsquare2}
\left(\int_{0}^{t}\lae(s)\exp\left(-\int_{s}^{t}\lae(\tau)\,\d \tau\right)\eta(s)\,\d s\right)^{2} \\
\leq \frac{1}{2\delta}\int_{0}^{t}\lae(s)\,\eta^{2}(s)\exp\left(-2(1-\delta)\int_{s}^{t}\lae(\tau)\,\d \tau\right)\,\d s\,.\end{multline}
\end{lem}
\begin{proof} We have
\begin{multline*}
\int_{0}^{t}\lae(s)\exp\left(-\int_{s}^{t}\lae(r)\,\d r\right)\int_{0}^{s}\exp\left(-\frac{\mu_0}{\e^{2}}(s-\tau)\right)\eta(\tau)\,\d \tau\, \d s\\
=\exp\left(-\int_{0}^{t}\lae(r)\,\d r\right)\int_{0}^{t}\exp\left(\frac{\mu_0}{\e^{2}}\tau\right)\eta(\tau)\int_{\tau}^{t}\lae(s)\exp\left(\int_{0}^{s}\lae(r)\,\d r-\frac{\mu_0}{\e^{2}}s\right)\,\d s\, \d\tau\,.
\end{multline*} Let us introduce
\[
I(\tau, t) := \int_{\tau}^t \lae(s) \exp\left(\int_0^s \lae(r) \d r - \frac{\mu_0 s}{\varepsilon^2} \right) \d s.
\]
We factor out the dominant exponential term at $s = \tau$ to write
\begin{equation*}\begin{split}
I(\tau, t)
&= \exp\left( \int_0^\tau \lae(r) \,\d r - \frac{\mu_0\tau}{\varepsilon^2} \right)
\int_{\tau}^t \lae(s)
\exp\left( \int_\tau^s \lae(r) \,\d r - \frac{\mu_0}{\varepsilon^2}(s-\tau) \right) \,\d s\\
&=\exp\left( \int_0^\tau \lae(r) \,\d r - \frac{\mu_0\tau}{\varepsilon^2} \right)
\int_0^{t - \tau} \lae(\tau + u)
\exp\left( \int_0^u \lambda_{\varepsilon}(\tau + r) \,\d r - \frac{\mu_0u}{\varepsilon^2} \right) \, \d u\,.\end{split}\end{equation*}
Since 
$
\int_0^u \lae(\tau + r) \, \d r \leq \lambda_{\star}u,$
we have
$$I(\tau, t)
\leq \exp\left( \int_0^\tau \lae(r) \, \d r - \frac{\mu_0 \tau}{\varepsilon^2} \right)\int_{0}^{t-\tau}\lae(\tau+u)\exp\left(\lambda_{\star}u-\frac{\mu_0}{\e^{2}}u\right)\d u\,.$$
{We then notice that $\lae(\cdot)$ is such that fo any $u$, $\tau \geq 0$,
$$\lae(\tau+u) \lesssim \lae(\tau)\,.$$
Indeed, from~\eqref{eq:zlambda} and the fact that $z(\cdot)$ is non-increasing, we have 
$$
\lae(\tau + u) \lesssim z^\gamma(\tau+u) \leq z^\gamma(\tau) \lesssim \lae(\tau) \,.
$$}
We deduce that, as soon as $\lae(t) \e^{2} \leq \frac{\mu_0}{2}$ for any $t \geq 0$, which is the case for $\e \in (0,\e_3)$ by definition of $\e_3$ (see the proof of Proposition~\ref{prop:precisePo}), we have
\begin{equation*}\begin{split}
I(\tau,t) &\lesssim  {\lae(\tau)}\exp\left( \int_0^\tau \lae(r) \,\d r - \frac{\mu_0\tau}{\varepsilon^2} \right)\int_{0}^{\infty}\exp\left(-\frac{\mu_0}{2\e^{2}}u\right)\d u\\
&\lesssim  {\lae(\tau)\e^{2}}\exp\left( \int_0^\tau \lae(r) \, \d r - \frac{\mu_0\tau}{\varepsilon^2} \right)\end{split}\end{equation*}
which gives the first part of the result. To prove \eqref{eq:integsquare2}, we adapt the proof of \cite[Eq. (4.12)]{ALT} (given for constant mapping $\lae$) and use Cauchy Schwarz inequality to deduce that,
\begin{align*}
&\quad \left(\int_{0}^{t}\lae(s)\exp\left(-\int_{s}^{t}\lae(r)\d r\right)\eta(s)\,\d s\right)^{2}\\
&\leq \left(\int_{0}^{t}\lae(s)\exp\left(-2\delta\int_{s}^{t}\lae(r)\d r\right)\d s\right)\left(\int_{0}^{t}  {\lae(s)}\exp\left(-2(1-\delta)\int_{s}^{t}\lae(r)\d r\right)\eta(s)^{2}\,\d s\right) \\
&\leq \frac{1}{2\delta} \left(\int_{0}^{t}  {\lae(s)}\exp\left(-2(1-\delta)\int_{s}^{t}\lae(r)\d r\right)\eta(s)^{2}\,\d s\right)\,,
\end{align*}
which gives the result. \end{proof} 

  \begin{cor}
Consider $\delta \in \left(0,\frac{1}{\g+1}\right)$. For any $\e \in (0,\e_{3})$ (where $\e_3$ is defined in Proposition~\ref{prop:precisePo}) and any~$t \geq 0$, there holds
\begin{multline} \label{eq:P0h12-b}
\|\mathbf{P}_{0} \hu(t)  \|^2_{\H} \lesssim \e^{2\theta} \exp\left(-2(1-\delta) \frac{\g}{\g+1}\int_{0}^{t}\lae(\tau)\,\d\tau\right)+  
\|h_{\rm in}\|_{\E}^{2}  \exp\left(-2\int_{0}^{t}\lae(\tau)\,\d\tau\right) \\
+\frac{1}{\delta} \int_{0}^{t}\lae(s)\exp\left(-2(1-\delta)\int_{s}^{t}\lae(\tau)\,\d \tau\right)\left\|\left(\mathbf{Id-P}_{0}\right)\hu(s)\right\|_{\H}^{2} \,\d s \end{multline}
where $\theta$ is defined in~\eqref{eq:Pi0Seps}.
\end{cor}
\begin{proof}
We first deduce from~\eqref{eq:P0h11} that
\begin{multline*}
\|\mathbf{P}_{0} \hu(t)  \|^2_{\H} \lesssim \e^{2\theta}  \exp\left(-2(1-\delta)\frac{\g}{\g+1}\int_0^t \lae(\tau)\,\d\tau\right)+
\|h_{\rm in}\|_{\E}^{2}  \exp\left(-2\int_{0}^{t}\lae(\tau)\d\tau\right)\\ 
+\left(\int_{0}^{t}\lae(s)\exp\left(-\int_{s}^{t}\lae(\tau)\d \tau\right)\left\|\left(\mathbf{Id-P}_{0}\right)\hu(s)\right\|_{\H} \d s\right)^{2}\,.
\end{multline*}
We then estimate the squared integral using  \eqref{eq:integsquare2} to deduce \eqref{eq:P0h12-b}. 
\end{proof}
\color{black}
\subsection{Estimating the complement $(\mathbf{Id-P}_0)\hu$} 

Let us focus on an estimate on $\mathbf{P}_0^\perp \hu(t)$ with $\mathbf{P}_0^\perp=\mathbf{Id-P}_{0}$, the orthogonal projection onto $\left(\mathrm{Ker}(\G_{\e})\right)^{\perp}$ in the Hilbert space $L^2_{x,v}(\M^{-\frac{1}{2}})$.  The same notation for the linearized operator $\G_{\e}$ in the spaces $\E$ and $\H$ is used. As in \cite{ALT}, it is crucial in our analysis that the nonlinear term that appears in \eqref{eq:h1} in the form $(\mathbf{Id}-\bm{\pi}_0) \Q_{\e,t}(\hu,\hu)$ since it ensures that 
$$\bm{\pi}_{0}(\mathbf{Id}-\bm{\pi}_0) \Q_{\e,t}(\hu,\hu)=0\,.$$
The approach in this section is similar in that of \cite[Section 4.4]{ALT} and resorts to the hypocoercivity results recalled in Section \ref{subsec:hypo}. Namely, using Proposition~\ref{prop:hypoco} together with Lemmas~\ref{lem:hypodrift} and~\ref{theo:briant2.4}, we are able to obtain some nice estimate on $\mathbf{P}_0^\perp\hu$.}
\begin{lem}\label{lem:h1perp} With the notations of  Proposition \ref{prop:hypoco}, assume that $\Mo$ defined in~\eqref{def:Mo} is small enough so that 
\begin{equation}\label{eq:mu1}
\widetilde{\mu}_1:=2\mathrm{a}_{1}-c_{1}\Mo^{2} >0\end{equation}
where  $c_{1} >0$ is a universal constant defined in~\eqref{def:c1}. Then, for any $\e \in (0,\e_{3})$ (where $\e_{3}$ is defined in Proposition~\ref{prop:precisePo}) and any $t \geq 0$, there holds 
\begin{multline}\label{eq:huut1-0}
\|(\mathbf{Id} - \mathbf{P}_0) h^1 (t)\|_{\H}^{2} 
\lesssim 
 \int_{0}^{t}e^{-\mu_1(t-s)}\left[\left( \lambda_\star+\Mo\right) \|\hu(s)\|_{\H}^{2}+ {\lae(s)}^2\right]\, \d s\\
+\frac{1}{\e^2}\int_{0}^{t}e^{-\mu_1(t-s)}\|\hu(s)\|_{\H}\,\|\ho(s)\|_{\E}\,\d s\,\end{multline} 
where $\mu_1:=C_{\H}^2\widetilde\mu_1$ with $C_{\H}$ is the constant appearing in \eqref{eq:CH-CH} related to the equivalence between to $\vertiii{\cdot}_{\H}$ and $\|\cdot\|_{\H}.$
\end{lem}
\begin{proof} Set $\Psi(t) := \mathbf{P}_0^\perp h^1 (t)$ for any $t \geq 0$. We start by recalling that $h^1(0)=0$ so that $\Psi(0)=0$. 
One checks from \eqref{eq:h1} (using in particular that $\mathbf{P_0} \bm{\pi}_0 = \mathbf{P}_0$ for the second term) that \begin{align*}
\dfrac{1}{2}\dfrac{\d}{\d t}\vertiii{\huu(t)}_{\H}^{2}&=\langle\!\langle \G_{\e}\huu(t),\huu(t)\rangle\!\rangle_{\H} +\e^{-1}\langle\!\langle \left(\mathbf{Id}-\bm{\pi}_{0}\right)\Q_{\e,t}(\hu(t),\hu(t)),\huu(t)\rangle\!\rangle_{\H}\\
&\quad
+\xi(t)\bigg(-\langle\!\langle \mathrm{div}_{v}(v\huu(t)),\huu(t)\rangle\!\rangle_{\H}-\mathrm{a}_3\langle\!\langle  \m_{2}\huu(t),\huu(t)\rangle\!\rangle_{\H}\bigg)\\
&\quad
+\xi(t)\langle\!\langle  \mathbf{P}_{0} \mathrm{div}_{v}\left(v \hu\right)- \mathrm{div}_{v}\left(v \mathbf{P}_{0}\hu\right) + \mathrm{a}_3 \mathbf{P}_0 (\m_2 \huu),\huu(t)\rangle\!\rangle_{\H}\\
&\quad
+\langle\!\langle\mathbf{P}_{0}^{\perp}\A_{\e}\ho(t),\huu(t)\rangle\!\rangle_{\H}  +\langle\!\langle\mathbf{P}_{0}^{\perp}S_{\e}(t),\huu(t)\rangle\!\rangle_{\H} \\
&
=:I_{1}+I_{2}+I_{3}+I_{4}+I_{5}+I_{6}\,.
\end{align*}
Each of the terms are estimated independently as in the proof of \cite[Lemma 4.6]{ALT} to which we refer for further details. First, according to Proposition \ref{prop:hypoco},
$$
I_{1}\leq -\frac{\mathrm{a}_{1}}{\e^{2}}\|\left(\mathbf{Id}-\bm{\pi}_{0}\right)\huu(t)\|_{\H_{1}}^{2}-\mathrm{a}_{1}\|\huu(t)\|_{\H_{1}}^{2}$$
 and 
$$I_{3} \lesssim \e\,\xi(t) \|\huu(t)\|_{\H}^{2} \lesssim \xi(t)\|\hu(t)\|_{\H}^{2}
$$
where we used \eqref{lem:hypodrift} for the estimate of $I_{3}$ and notice the use of the parameter $\mathrm{a}_{3}$ in the splitting \eqref{eq:h0}-\eqref{eq:h1}. Moreover, we deduce from  Lemma \ref{theo:briant2.4} and  Young's inequality that there is $C >0$ such that
\begin{equation}\label{eq:I2}
I_{2} \leq \frac{\eta}{\e^{2}}\left\|\huu^{\perp}(t)\right\|_{\H_{1}}^{2}+\frac{C}{\eta}\|\hu(t)\|_{\H}^{2}\|\hu(t)\|_{\H_{1}}^{2} +C\|\hu(t)\|_{\H}^{3}\,, \qquad  \eta >0\,.
\end{equation}
Using the regularization properties of $\mathbf{P}_0$ in velocity together with~\eqref{def:xi} and \eqref{eq:zlambda}, one has
$$I_{4} \lesssim \lambda_\e(t)\|\hu(t)\|_{\H}^{2}\,.$$
Concerning $I_5$, using the regularization properties in velocity of $\A_\e$ (see Lemma~\ref{lem:dissip}), we have 
$$
I_5 \lesssim \frac{1}{\e^2} \|h^1(t)\|_\H \|h^0(t)\|_\E\,. 
$$
Finally, to estimate $I_{6}$ we recall that {$S_\e=S_\e(t,v)$ is defined in~\eqref{def:LetSe}}. Moreover, on can check that 
$$
\operatorname{div}_v(v\M)=3\M-v\cdot \nabla_{v}\M = 3\M-\frac{|v|^{2}}{\en_{\star}}\M=-\sqrt{6}\varphi_{5}\M
$$ 
from which we deduce that 
$\mathrm{div}_{v}(v\M) \in \mathrm{Range}(\bm{\pi}_{0})$
and thus $\left(\mathbf{Id}-\bm{\pi}_{0}\right) \operatorname{div}_v(v\M)=0$.
In particular, noticing that $S_\e$ depends only on $v$ so that $\left(\mathbf{Id-P}_{0}\right)S_{\e}=\left(\mathbf{Id}-\bm{\pi}_{0}\right)S_{\e}$, one has
$$\left(\mathbf{Id-P}_{0}\right)S_{\e}(t)=\e^{-3}\left(\mathbf{Id}-\bm{\pi}_{0}\right)\Q_{\e,t}(\M,\M) \,.$$
Using the definition of $\langle\!\langle\cdot,\cdot\rangle\!\rangle_{\H}$, since $\Q_{\e,t}(\M,\M)=\Q_{\e,t}(\M,\M)-\Q_{1}(\M,\M)$ is \emph{independent of $x$}) 
\begin{equation*}
\langle\!\langle\mathbf{P}_{0}^{\perp}S_{\e}(t),\huu(t)\rangle\!\rangle_{\H}
=\e^{-3}\langle  \left(\mathbf{Id}-\bm{\pi}_{0}\right)\Q_{\e,t}(\M,\M), \huu^{\perp}(t)\rangle_{\H} 
\end{equation*} and
$$\langle\!\langle\mathbf{P}_{0}^{\perp}S_{\e}(t),\huu(t)\rangle\!\rangle_{\H} \lesssim \e^{-3}\|\Q_{\e,t}(\M,\M)\|_{\H}\,\|\huu^{\perp}(t)\|_{\H}
\,.$$
Now, using Lemma \ref{lem:QEMM} with $\ell = \ell_\e(t) = \e^{\frac2\g} z(t)$ together with estimate~\eqref{eq:zlambda}, $$\|\Q_{\e,t}(\M,\M)\|_{\H} \lesssim \e^{2}\lae(t)\,,$$
which, with Young's inequality, implies that
\begin{equation} \label{eq:I6}
I_{6} \lesssim \e^{-1}\lae(t)\|\huu^\perp(t)\|_{\H}
\lesssim \frac{\eta}{\e^{2}}\|\huu^{\perp}(t)\|_{\H}^{2}+\frac{C}{\eta}{\lae(t)^2}\,, \qquad \eta >0\,.
\end{equation}
Therefore, choosing  {$2\eta \leq  \mathrm{a}_{1}$} in \eqref{eq:I2}-\eqref{eq:I6}, one sees that  there is some positive constant $c_{0} >0$ such that
\begin{multline}\label{eq:bmh1}
\dfrac{\d}{\d t}\vertiii{\huu(t)}_{\H}^{2} \leq -2\mathrm{a}_{1}\|\huu(t)\|_{\H_{1}}^2 + c_0\|\hu(t)\|_{\H}^2\left(\|\hu(t)\|_{\H_{1}}^2 +\|\hu(t)\|_{\H}\right) +c_{0} {\lae(t)^2}\\
+c_0\lae(t) \|\hu(t)\|_{\H}^{2}+\frac{c_{0}}{\e^{2}}\|\ho(t)\|_{\E} \|\hu(t)\|_\H\,.
\end{multline}
Using then~\eqref{def:Mo}, we deduce first that there is an explicit $c_{1} >0$ such that
\begin{equation}\label{def:c1}
c_0 \|\hu(t)\|_{\H}^2\,\left(\|\hu(t)\|_{\H_{1}}^2 + \|\hu(t)\|_\H\right)
\leq c_1 \Mo^{2}\|\huu(t)\|_{\H_{1}}^{2}+c_1\Mo \|\hu(t)\|^2_\H
\end{equation} 
and then  that, for $\Mo$ small enough so that  $\widetilde \mu_1:={2\mathrm{a}_{1}}-c_{1}\Mo^{2} >0$,
  there exists $c_2 >0$ such that
\begin{multline*}
\dfrac{\d}{\d t}\vertiii{\huu(t)}_{\H}^{2} \leq -\widetilde\mu_1\,\|\huu(t)\|_{\H_{1}}^2   
+ c_2\left( \lae(t)+\Mo\right) \|\hu(t)\|_{\H}^{2} \\+\frac{c_2}{\e^2} \|\hu(t)\|_{\H}\,\|\ho(t)\|_{\E}+c_{2}{\lae(t)^2}  \,.\end{multline*}
Integrating this differential inequality yields the desired estimate \eqref{eq:huut1-0}, recalling that $\hu(0)=0$ and $\|\huu(t)\|_{\H_{1}} \geq \|\huu(t)\|_{\H} \geq C_{\H}\vertiii{\huu(t)}_{\H}.$ 
\end{proof} 
To complete the estimate of $\|(\mathbf{Id} - \mathbf{P}_0) \hu(t)\|_{\H}^{2}$ we need to estimate the  last integral term in~\eqref{eq:huut1-0}.
\begin{lem}\label{lem:huut1}
Consider  $\delta \in \left(0,\frac{1}{\g+1}\right)$. There exists {an explicit $\e_{4}\in (0,\e_{3})$} (where $\e_3$ has been defined in Proposition~\ref{prop:precisePo}) such that for any $\e \in (0,\e_{4})$ and any $t \geq 0$,
\begin{multline}\label{lem:huut1-b}
\frac{1}{\e^2}\int_{0}^{t}e^{-\mu_1(t-s)}\|\hu(s)\|_{\H}\,\|\ho(s)\|_{\E}\,\d s \\
\lesssim \|h_{\rm in}\|_{\E}e^{-\mu_1t} + \exp\left(-2(1-\delta) \frac{\g}{\g+1} \int_0^t \lae(s) \, \d s \right) +  \lambda_\star^{2\delta} \int_0^t e^{-\mu_1(t-s)} \|\hu(s)\|^2_\H \, \d s\,.\end{multline}
\end{lem}
\begin{proof} 
We use~\eqref{eq:estimah00},~\eqref{def:Mo} and~\eqref{eq:zlambda} to deduce that 
\begin{multline*}
\frac{1}{\e^2}\int_{0}^{t}e^{-\mu_1(t-s)}\|\hu(s)\|_{\H}\,\|\ho(s)\|_{\E}\,\d s \\
\lesssim \frac{1}{\e^2} \|h_{\rm in}\|_\E \Delta_0 e^{-\mu_1 t} \int_0^t e^{- \left(\frac{\mu_0}{\e^2} - \mu_1\right)s} \, \d s 
+ \int_0^t e^{-\mu_1(t-s)} \lae(s) \|h^1(s)\|_\H \, \d s\,.
\end{multline*}
Choosing then $\e_4$ small enough so that $2 \e_4^2 \mu_1 \leq \mu_0$ and using Young's inequality for the second term, we deduce 
\begin{multline*}
\frac{1}{\e^2}\int_{0}^{t}e^{-\mu_1(t-s)}\|\hu(s)\|_{\H}\,\|\ho(s)\|_{\E}\,\d s \\
\lesssim  \|h_{\rm in}\|_\E e^{-\mu_1 t}
+ \int_0^t e^{-\mu_1(t-s)} \lae(s)^{2(1-\frac{\delta}{2})} \, \d s  
+ \int_0^t e^{-\mu_1(t-s)} \lae(s)^{\delta}  \|h^1(s)\|^2_\H \, \d s \,. 
\end{multline*}
Arguing as in Lemma \ref{lem:IntAlg} it is not difficult to check that
\begin{equation}\label{eq:intlambdae}
\int_0^t e^{-\mu_1(t-s)} \lae (s)^{\alpha} \, \d s \lesssim \lae(t)^{\alpha} \qquad \forall \,\alpha >0\,,
\end{equation} and applying this with $\alpha= 2(1-\frac{\delta}{2})$ gives 
\begin{multline*}
\frac{1}{\e^2}\int_{0}^{t}e^{-\mu_1(t-s)}\|\hu(s)\|_{\H}\,\|\ho(s)\|_{\E}\,\d s \\
\lesssim \|h_{\rm in}\|_{\E}e^{-\mu_1t} + \lae(t)^{2-\delta} +  \lambda_\star^{2\delta} \int_0^t e^{-\mu_1(t-s)} \|\hu(s)\|^2_\H \, \d s\,.
\end{multline*}
To conclude, we use~\eqref{eq:zgexp} (with $\delta/2$ instead of $\delta$) to deduce that 
$$
\lae(t)^{2-\delta}\lesssim \exp\left(- (2-\delta)\left(1-\frac{\delta}{2}\right) \frac{\g}{\g+1}\int_0^t \lae(s) \, \d s \right)\,. 
$$
This yields the wanted result since $(2-\delta) (1-\frac{\delta}{2}) = 2 (1-\frac{\delta}{2})^2 \geq 2(1-\delta)$. 
\end{proof}
 \color{black}

We deduce from the previous the following main estimate for $\|(\mathbf{Id} - \mathbf{P}_0) \hu(t)\|_{\H}$:
\begin{prop} \label{prop:Psi}
Consider $\delta \in \left(0,\frac{1}{\g+1}\right)$. For any $\e \in (0,\e_{4})$ (where $\e_{4}$ is defined in Lemma~\ref{lem:huut1}), it holds
\begin{multline}\label{eq:Psi}
\|(\mathbf{Id} - \mathbf{P}_0) h^1 (t)\|_{\H}^{2} 
\lesssim  {\|h_{\rm in}\|_{\E} e^{-\mu_1t}
+ \exp\left(-2(1-\delta) \frac{\g}{\g+1} \int_0^t \lae(s) \, \d s \right) } \\+
\left(\lambda_\star^{\min(1,2\delta)}+\Mo\right) \int_{0}^{t}e^{-\mu_1(t-s)}\|\hu(s)\|_{\H}^{2}\, \d s\,.
\end{multline} 
In particular, if one assumes 
$2\lambda_{\star} < \mu_1$
(where $\mu_1$ is defined in Lemma~\ref{lem:h1perp}),
then, for any $\e \in (0,\e_{4})$ and any $\delta \in (0,1)$, it holds 
\begin{multline}\label{eq:intPsib}
\int_{0}^{t}\lae(s)\exp\left(-2(1-\delta)\int_{s}^{t}\lae(\tau)\,\d \tau\right)\left\|\left(\mathbf{Id-P}_{0}\right)\hu(s)\right\|_{\H}^{2} \d s \\
\lesssim\|h_{\rm in}\|_{\E} \exp\left(-2(1-\delta)\int_{0}^{t}\lae(\tau)\,\d\tau\right) 
+ \exp\left(-2(1-\delta) \frac{\g}{\g+1} \int_0^t \lae(\tau)\, \d\tau\right) \\
+ \left(\lambda_\star^{\min(1,2\delta)} +\Mo\right) \int_{0}^{t}\lae(s) \exp\left(-2(1-\delta)\int_{s}^{t}\lae(\tau)\,\d\tau\right)\|\hu(s)\|_{\H}^{2}\, \d s\,.
\end{multline}
\end{prop}

\begin{proof}  We insert the estimates obtained in Lemma \ref{lem:huut1}  into \eqref{eq:huut1-0}. We then use~\eqref{eq:intlambdae} with $\alpha=2$ to get \eqref{eq:Psi} and then use the fact that 
$$
\lae(t)^{2} \lesssim \lae(t)^{2(1-\delta)} \lesssim \exp \left(- 2(1-\delta) \frac{\g}{\g+1} \int_0^t \lae(s) \, \d s \right)\,.
$$
To get~\eqref{eq:intPsib},  we use the facts that 
$$\int_{0}^{t}\lae(s)\exp\left(-2(1-\delta)\int_{s}^{t}\lae(\tau)\,\d \tau-\mu_1s\right)\d s \lesssim \exp\left(-2(1-\delta) \int_{0}^{t}\lae(\tau)\,\d \tau\right)\,,$$
and 
\begin{align*}
&\int_0^t \lae(s) \exp \left(-2(1-\delta) \int_s^t \lae(\tau) \, \d\tau \right) \exp \left(-2(1-\delta) \frac{\g}{\g+1} \int_0^s \lae(\tau) \, \d\tau \right) \\
&\quad 
= \exp \left(-2(1-\delta) \frac{\g}{\g+1} \int_0^t \lae(\tau) \, \d\tau \right) 
\int_0^t \lae(s) \exp \left(-2(1-\delta) \frac{1}{\g+1} \int_s^t \lae(\tau) \, \d\tau \right) \, \d s \\
&\quad
\lesssim  \exp \left(-2(1-\delta) \frac{\g}{\g+1} \int_0^t \lae(\tau) \, \d\tau \right) 
\end{align*}
and, as in Lemma \ref{lem:integdouble},
\begin{multline*}
\int_{0}^{t}\lae(s)\exp\left(-2(1-\delta)\int_{s}^{t}\lae(r)\,\d r\right)\int_{0}^{s}e^{-\mu_1(s-\tau)}\eta(\tau)\,\d\tau \,\d s\\
\lesssim \int_{0}^{t}\exp\left(-2(1-\delta)\int_{\tau}^{t}\lae(r)\,\d r\right) {\lae(\tau)}\eta(\tau)\,\d\tau\end{multline*}
for any nonnegative $\eta(\cdot)$  provided $\mu_1 >2\lambda_{\star}$ so that $\mu_1 > 2(1-\delta)\lae(\tau)$ for any $\tau$. We thus obtain the wanted result.
\end{proof}

\subsection{Final \emph{a priori} estimates}
We deduce from the above the following estimate on $h^1$:
\begin{prop}\label{prop:h1fin}
Consider $\delta \in \left(0,\frac{1}{\g+1}\right)$. Under the assumptions of Lemma \ref{lem:h1perp}, we have for any~$t \geq0$, any~$\e \in (0,\e_{4})$ (where $\e_{4}$ is defined in Lemma~\ref{lem:huut1}) and for $\lambda_{\star} < 2\mu_{1}$,
\begin{equation}\label{eq:hufin}
\begin{aligned}
\|h^1(t)\|^2_\H 
&\lesssim \left(1 + \|h_{\rm in}\|_\E +\|h_{\rm in}\|_\E^2\right) \exp\left(-2(1-\delta) \frac{\g}{\g+1} \int_0^t \lae(\tau) \, \d \tau\right) \\
&\quad
+\left(\lambda_\star^{\min(1,2\delta)} + \Delta_0\right) \int_{0}^{t}\lae(s) \exp\left(-2(1-\delta)\int_{s}^{t}\lae(\tau)\,\d \tau\right)\|\hu(s)\|_{\H}^{2}\,\d s \\
&\quad
+  \left(\lambda_\star^{\min(1,2\delta)}+\Mo\right)  \int_0^t e^{-\mu_1(t-s)}\|h^1(s)\|^2_\H \, \d s \,.
\end{aligned}
\end{equation}
 
\end{prop} 
\begin{proof} Since $\|\hu(t)\|_{\H}^{2}=\|\mathbf{P}_{0}\hu(t)\|_{\H}^{2}+\|\left(\mathbf{Id-P}_0\right)\hu(t)\|_{\H}^{2}$, we gather estimates~\eqref{eq:Psi} and~\eqref{eq:P0h12-b} (combined with~\eqref{eq:intPsib}) to obtain that 
\begin{equation*} 
\begin{aligned}
\|h^1(t)\|^2_\H 
&\lesssim \|h_{\rm in}\|_\E \,e^{-\mu_1 t} +
 \left(\|h_{\rm in}\|^2_\E + \frac{1}{\delta}\|h_{\rm in}\|_\E\right) \, \exp\left(-2(1-\delta)\int_{0}^{t}\lae(\tau)\,\d \tau\right) 
 \\
 &\quad + \frac1\delta\exp\left(-2(1-\delta) \frac{\g}{\g+1} \int_0^t \lae(\tau) \, \d \tau\right)\\
&\quad
+ \frac1\delta (\lambda_\star^{\min(1,2\delta)} + \Delta_0) \int_{0}^{t}\lae(s) \exp\left(-2(1-\delta)\int_{s}^{t}\lae(\tau)\,\d \tau\right)\|\hu(s)\|_{\H}^{2}\,\d s \\
&\quad
+  \left(\lambda_\star^{\min(1,2\delta)}+\Mo\right)  \int_0^t e^{-\mu_1(t-s)}\|h^1(s)\|^2_\H \, \d s \,,
\end{aligned}
\end{equation*}
which allows us to conclude to~\eqref{eq:hufin}.  
\end{proof}

By a refined Gronwall type argument, we are able to derive from this the following decay rate for~$\|\hu(t)\|_{\H}$:
\begin{cor}\label{cor:h1fin}
Let $\delta \in \left(0,\frac{1}{\g+1}\right)$. There exist $\e_\delta$,  $\lambda_\delta, \Delta_{0,\delta}>0$ and $C_\delta>0$ depending on $\delta$ such that for any $\e \in (0,\e_\delta)$, $\lambda_\star \in [0,\lambda_\delta)$ , $\Mo \in (0,\Delta_{0,\delta})$ (where $\Mo$ is defined in~\eqref{def:Mo}),
\begin{equation}\label{eq:h1fin}
\|\hu(t)\|_{\H}^{2}\leq C_\delta  \left(1+\|h_{\rm in}\|_{\E}+ \|h_{\rm in}\|_{\E}^{2}\right)\lae(t)^{2(1-\delta)}
\end{equation}
holds true for any $t \geq0.$
\end{cor}
\begin{proof}
We set $\delta':=\frac{\delta}{4}$. First, we use the estimate~\eqref{eq:hufin} with $\delta'$ instead of $\delta$ which gives 
\begin{align*}
\|h^1(t)\|^2_\H
&\lesssim \mathcal K_0 \exp\left(-2(1-\delta') \frac{\g}{\g+1} \int_0^t \lae(\tau) \, \d \tau\right) \\
&\quad
+ \mathcal K_1 \int_{0}^{t}\lae(s) \exp\left(-2(1-\delta')\int_{s}^{t}\lae(\tau)\,\d \tau\right)\|\hu(s)\|_{\H}^{2}\,\d s \\
&\quad
+ \mathcal K_1  \int_0^t e^{-\mu_1(t-s)}\|h^1(s)\|^2_\H \, \d s \,,
\end{align*}
where $\mathcal K_0 := 1 + \|h_{\rm in}\|_\E +\|h_{\rm in}\|_\E^2$ and $\mathcal K_1 := \lambda_\star^{\min(1,2\delta)}+\Mo$. Note that in the previous inequality, the multiplicative constant depends on $\delta$. We also have to keep in mind that in the end, we can choose $\mathcal K_1$ as small as necessary. Introducing then the new function 
\begin{equation} \label{def:bmx}
\bm{x}(t) := \|h^1(t)\|^2_\H \exp \left(2(1-\delta') \frac{\g}{\g+1} \int_0^t \lae(\tau) \, \d\tau\right)\,, 
\end{equation}
we deduce that 
\begin{align*}
\bm{x}(t) &\lesssim \mathcal K_0 +
\mathcal K_1 \int_0^t \lae(s) \exp\left(-2(1-\delta')\frac{1}{\g+1}\int_s^t \lae(\tau)\,\d\tau\right) \bm{x}(s) \, \d s \\
&\quad
+ \mathcal K_1 \int_0^t \exp\left(-\int_s^t \left(\mu_1-2(1-\delta') \frac{\g}{\g+1}\lae(\tau)\right)\,\d\tau\right) \bm{x}(s) \, \d s 
=: \bm{y}(t)\,.
\end{align*}
We then compute
\begin{align*}
\bm{y}'(t) 
&=
-2 \mathcal K_1(1-\delta')\frac{1}{\g+1} \lae(t) \int_0^t \lae(s) \exp\left(-2(1-\delta')\frac{1}{\g+1}\int_s^t \lae(\tau)\,\d\tau\right) \bm{x}(s) \, \d s \\
&\quad 
- \mathcal K_1\left(\mu_1 - 2(1-\delta') \frac{\g}{\g+1} \lae(t)\right) \\
&\qquad \qquad \qquad 
\times \int_0^t \exp\left(-\int_s^t \left(\mu_1-2(1-\delta') \frac{\g}{\g+1}\lae(\tau)\right)\,\d\tau\right) \bm{x}(s) \, \d s \\
&\quad
+ \mathcal K_1 \left( \lae(t)+1\right) \bm{x}(t) 
\end{align*}
so that 
\begin{align*}
\bm{y}'(t)
&\leq
\left(\mu_1 - 2(1-\delta') \frac{\g}{\g+1} \lae(t)\right) \\
&\qquad \qquad 
\times \left[\mathcal K_0 +
\mathcal K_1 \int_0^t \lae(s) \exp\left(-2(1-\delta')\frac{1}{\g+1}\int_s^t \lae(\tau)\,\d\tau\right) \bm{x}(s) \, \d s\right] \\
&\quad
- \left(\mu_1 - 2(1-\delta') \frac{\g}{\g+1} \lae(t)\right)\bm{y}(t) + \mathcal K_1 \left( \lae(t)+1\right) \bm{x}(t)\,.
\end{align*}
We then choose $\lambda_\star = \sup_{\e,t} \lae(t)$ small enough so that $\mu_1 - 2(1-\delta') \frac{\g}{\g+1} \lae(t) \geq \frac{\mu_1}{2}$ for any~$\e$ and any~$t \geq 0$. We also use the fact that there exists $C>0$ such that for any $\e$ and $t$, 
$$
\mathcal K_1 \left( \lae(t)+1\right) \bm{x}(t) \leq C \mathcal K_1 \bm{y}(t)\,.
$$
We can then choose $\mathcal{K}_1$ small enough so that $\mu_1 - 2(1-\delta') \frac{\g}{\g+1} \lae(t) - C \mathcal K_1 \geq \frac{\mu_1}{4}$. The previous inequality then implies 
\begin{align*}
\bm{y}'(t) &\leq - \frac{\mu_1}{4} \bm{y}(t) + 
C \left[\mathcal K_0 +
\mathcal K_1 \int_0^t \lae(s) \exp\left(-2(1-\delta')\frac{1}{\g+1}\int_s^t \lae(\tau)\,\d\tau\right) \bm{y}(s) \, \d s\right] \,.
\end{align*}
Multiplying this inequality by $\exp(\mu_1 t/4)$ and integrating in time gives 
$$
\bm{y}(t) \leq \bm{y}(0) + C\mathcal{K}_0 \int_0^t \exp\left(-\frac{\mu_1}{4}(t-s)\right) \, \d s + C \mathcal{K}_1 \int_0^t \lae(s) \bm{y}(s) \, \d s\,.
$$
Recalling that $\bm{y}(0) = \mathcal{K}_0$, we deduce that 
$$
\bm{y}(t) \lesssim \mathcal{K}_0 + \mathcal{K}_1 \int_0^t \lae(s) \bm{y}(s) \, \d s\,.
$$
The classical Gr\"onwall lemma then gives that there exists $C>0$ such that for any $t \geq 0$,
$$
\bm{x}(t) \lesssim \bm{y}(t) \lesssim \mathcal{K}_0 \exp\left( C \mathcal{K}_1 \int_0^t \lae(\tau) \, \d\tau \right) \,.
$$
Recalling the definition of $\bm{x}(t)$ in~\eqref{def:bmx}, we deduce that 
$$
\|h^1(t)\|^2_\H \lesssim \mathcal K_0 \exp\left( - \left(2(1-\delta') \frac{\g}{\g+1} - C \mathcal{K}_1\right) \int_0^t \lae(\tau) \, \d\tau \right) \,.
$$
We reduce again the value of $\mathcal{K}_1$ (which amounts to reduce the value of $\lambda_\star$ and $\Delta_0$) in order to have $C \mathcal K_1 \leq \frac{\delta}{2} \frac{\g}{\g+1}$ so that $2(1-\delta') \frac{\g}{\g+1} - C \mathcal{K}_1 \geq  {(2-\delta)\frac{\g}{1+\g}}$ and thus 
$$
\|h^1(t)\|^2_\H \lesssim \mathcal K_0 \exp\left( - 2\left(1-\frac{\delta}{2}\right) \frac{\g}{\g+1} \int_0^t \lae(\tau) \, \d\tau \right) \,.
$$
Using again~\eqref{eq:stimerese2} and~\eqref{def:xi}, we obtain 
$$
\|h^1(t)\|^2_\H \lesssim \mathcal K_0 \, \lae(t)^{2\left(1-\frac{\delta}{2}\right)^2} \,,
$$
which yields the conclusion since $(1-\frac{\delta}{2})^2 \geq 1-\delta$. 
\end{proof}

We are now able to state the main result of this section which provides a result of decay for the solution $h$ to~\eqref{eq:heps}:
\begin{theo}\label{theo:h-relaxation}  Let $\delta \in \left(0,\frac{1}{\g+1}\right)$. There exist $\e_\delta$,  $\lambda_\delta, \Delta_{0,\delta}>0$ and ${C}_\delta>0$ depending on $\delta$ such that for any $\e \in (0,\e_\delta)$, $\lambda_\star \in [0,\lambda_\delta)$ , $\Mo \in (0,\Delta_{0,\delta})$ (where $\Mo$ is defined in~\eqref{def:Mo}), the solution $h=h_{\e}$ to \eqref{eq:heps} satisfy
\begin{equation}\label{eq:finH}
\| h(t)\|_{\E} \leq {C}_\delta \,\Big(1 +  \| h_{\rm in}\|_{\E} \Big)\,\lambda_{\e}(t)^{1-\delta} + \e^{2}\lambda_{\e}(t)\,, \quad \forall \, t \geq 0\,. 
\end{equation}

\end{theo}
\begin{proof}
The result is simply obtained by gathering the estimate on $\|h^{0}(t)\|_{\E}$ in Proposition \ref{prop:h0} together with \eqref{eq:h1fin} to get
$$\|h(t)\|_{\E} \lesssim \|h_{\rm in}\|_{\E}\exp\left(-\frac{\mu_{0}}{\e^{2}}t\right)+\e^{2}\lambda_{\e}(t) + (1+\|h_{\rm in}\|_{\E})\lambda_{\e}(t)^{1-\delta}$$
where we used again \eqref{eq:zlambda}.   Keeping only the dominant terms yields \eqref{eq:finH}. 

\end{proof}

\section{Cauchy theory and hydrodynamic limit}\label{sec:hydro}

In this section, we deduce from the \emph{a priori} estimates in Section \ref{sec:nonlinear} the well-posedness of \eqref{eq:heps}. We also deduce some compactness of the family of solutions $(h_{\e})_{\e}$ and prove its convergence as~$\e \to 0$ (in some suitable sense) yielding the incompressible Navier-Stokes system as a limiting hydrodynamic description of the solution.  
\subsection{Cauchy theory}\label{sec:Cauchy}

Recall first that the functional spaces $\E$ and $\E_1$ are defined in~\eqref{eq:defE}-\eqref{eq:defE1}. By resuming the ideas of \cite{ALT},  the \emph{a priori} estimates in Section \ref{sec:nonlinear} are enough to prove the following well-posedness result. We leave the details to the reader.

\begin{theo}\label{theo:main-cauchy1}
Let $\delta \in \left(0,\frac{1}{\gamma+1}\right)$. There exists a pair $(\varepsilon_\delta,{\alpha_{0}})$ of positive constants that depend on the mass and energy of $F_{\mathrm{in}}^\e$,~$m$ and $q$ that appear in the definition of $\E$, ($\varepsilon_\delta$  additionally depending on $\delta$ whereas $\alpha_{0}$ is not) such that, for $\e\in(0,\e_\delta)$, if
$$\|h_{\mathrm{in}}^\e\|_{ \E} \leq {\alpha_{0}}$$
then the equation \eqref{eq:heps} has a unique solution $h_{\e} \in \mathcal{C}\big([0,\infty); \E\big) \cap L^1_{\rm loc}\big((0,\infty);\E_1\big)$. Moreover, it satisfies the following decay estimate 
$$
\left\|h_\e(t)\right\|_{ \E}\leq C (\delta,{\alpha_{0}})\,\left[\mathfrak{b}_{1}^{\frac{\g}{\g+1}}+\g\mathfrak{a}_{1}t\right]^{-(1-\delta)},
\qquad \forall \, t \in (0,T)\,,
$$
for some positive constant $C(\delta,{ \alpha_{0}}) >0$ independent of $\e$.\end{theo}    

\begin{nb} \label{nb:main-cauchy1}
Under the same assumptions, we can actually prove the following estimates (which will be useful in what follows) on $h^0_\e$ and $h^1_\e$ that are respectively solutions to~\eqref{eq:h0} and~\eqref{eq:h1}. Let~$T>0$,  according to Proposition~\ref{prop:h0}, 
	\begin{equation} \label{eq:estimh0}
	\|h^0_\e\|_{L^\infty((0,\infty)\,;\,\E)} \lesssim 1
	\quad \text{and} \quad 
	\|h^0_\e\|_{L^1((0,T)\,;\,\E_1)} \lesssim \e^2 
	\end{equation}
as well as
	\begin{equation} \label{eq:estimh1}
	\|h^1_\e\|_{L^\infty((0,\infty)\,;\,\H)} \lesssim 1
	\quad \text{and} \quad 
	\|h^1_\e\|_{L^2((0,T)\,;\,\H_1)} \lesssim 1
	\end{equation}
where we recall that the spaces $\E_1$, $\H$ and $\H_1$ are defined in~\eqref{eq:defE1} and~\eqref{eq:defH}. 
Notice that in the previous inequalities, the multiplicative constants only involve quantities related to the initial data of the problem {and $T>0$ for the $L^1$ and $L^2$ norms in time} and are independent of~$\e$. Such estimates are actually consequence of the construction of the solutions through an iterative scheme build upon the splitting \eqref{eq:h0}--\eqref{eq:h1} such that \eqref{eq:estimh0} and \eqref{eq:estimh1} hold true.
\end{nb}
 One can prove the following estimate for time-averages of the microscopic part of $h_\e$, namely on $(\mathbf{Id}-\bm{\pi}_{0})h_\e$, which in particular tells that this microscopic part vanishes in the limit $\e \to 0$:

\begin{lem} \label{lem:hperpholder}
For any $0 \leq t_1 \leq t_2 \leq T$, there holds:
	\begin{equation} \label{eq:hperpholder}
	\int_{t_1}^{t_2} \|(\mathbf{Id} - \bm{\pi}_0)h_\e(\tau)\|_{\E} \, \d \tau \lesssim \e \sqrt{t_2-t_1}\,,
	\end{equation}

where {the multiplicative constant depends on $T>0$} and we recall that $\bm{\pi}_0$ is the projection onto the kernel of $\mathscr{L}_1$ defined in~\eqref{def:pi0}. 
\end{lem}
\begin{proof} 
We first remark that 
	\begin{multline*}
	\int_{t_1}^{t_2} \|(\mathbf{Id} - \bm{\pi}_0)h_\e(\tau)\|_{\E} \, \d\tau
	\lesssim 
	\left(\int_{t_1}^{t_2} \|(\mathbf{Id} - \bm{\pi}_0)h^0_\e(\tau)\|^2_{\E} \, \d \tau\right)^{\frac{1}{2}} \sqrt{t_2-t_1} \\
	+ \left(\int_{t_1}^{t_2} \|(\mathbf{Id} - \bm{\pi}_0)h^1_\e(\tau)\|^2_{\H_1} \, \d\tau\right)^{\frac{1}{2}} \sqrt{t_2-t_1}\,.
	\end{multline*}
The first term is estimated thanks to~\eqref{eq:estimh0}, which gives:
	$$
	\int_{t_1}^{t_2} \|(\mathbf{Id} - \bm{\pi}_0)h^0_\e(\tau)\|^2_{\E} \, \d \tau
	\lesssim \|(\mathbf{Id} - \bm{\pi}_0)h^0_\e\|_{L^\infty((0,T) \, ; \, \E)} 
	\|(\mathbf{Id} - \bm{\pi}_0)h^0_\e\|_{L^1((0,T) \, ; \, \E_1)} \lesssim \e^2.  
	$$
Concerning the second one, we perform similar computations as in the proof of Lemma~\ref{lem:h1perp}. We recall that $h^1_\e$ solves~\eqref{eq:h1}
and consider $\vertiii{\cdot}_{\H}$ an hypocoercive norm on $\H$ (see~Proposition~\ref{prop:hypoco}). We then have, with $\huu(t)=\left(\mathbf{Id-P}_{0}\right)\hu_\e(t)$,	
{
	\begin{multline*}
	\frac12 \frac{\d}{\d t} \vertiii{\huu(t)}_{\H}^2 
	\leq - \frac{\mathrm{a}_1}{2\e^2} \|(\mathbf{Id}- \bm{\pi}_0)\huu(t)\|_{\H_1}^2 
	 - \mathrm{a}_1 \| \huu(t)\|_{\H_1}^2 +\lambda_{\e}(t)^{2}\\
	 +C \|h^1_\e(t)\|^2_\H \left(\|h^1_\e(t)\|^2_{\H_1} + \|h^1_\e(t)\|_\H + \lambda_{\e}(t)  \right) + \frac{C}{\e^2} \|h^0_\e(t)\|_\E \|h^1_\e(t)\|_\H\,
	\end{multline*}}
from which we deduce that 
	\begin{multline*}
	\frac{1}{\e^2} \int_{t_1}^{t_2} \|(\mathbf{Id} - \bm{\pi}_0)\huu(\tau)\|^2_{\H_1} \, \d\tau 
	\lesssim \|h^1_\e(t_1)\|^2_{\H} +\int_{t_{1}}^{t_{2}}\lambda_{\e}(\tau)^{2}\,\d\tau\\
	+ \int_{t_1}^{t_2}  \|h^1_\e(\tau)\|^2_\H \left(\|h^1_\e(\tau)\|^2_{\H_1} + \|h^1_\e(\tau)\|_\H +  {\lae(\tau)}\right) \, \d\tau  
	+ \frac{1}{\e^2}\int_{t_1}^{t_2} \|h^0_\e(\tau)\|_\E \|h^1_\e(\tau)\|_\H \, \d \tau \\
	\lesssim  \|h^1_\e(t_1)\|^2_{\H} +\int_{t_{1}}^{t_{2}}\lambda_{\e}(\tau)^{2}\,\d\tau +  \int_{t_1}^{t_2}  \|h^1_\e(\tau)\|^2_\H \left(\|h^1_\e(\tau)\|^2_{\H_1} + \|h^1_\e(\tau)\|_\H \right) \, \d\tau  \\
	+ \frac{1}{\e^2}\int_{t_1}^{t_2} \|h^0_\e(\tau)\|_\E \|h^1_\e(\tau)\|_\H \, \d \tau
	\lesssim 1
	\end{multline*}
where we used the fact that the norm $\vertiii{\cdot}_{\H}$ is equivalent to the usual one $\|\cdot\|_\H$ uniformly in $\e$ as well as~\eqref{eq:estimh0} and~\eqref{eq:estimh1} to get the last estimate.  Now, since $\bm{\pi}_{0}\mathbf{P}_{0}=\mathbf{P}_{0}$, one has
$$(\mathbf{Id} - \bm{\pi}_0)h^1_\e=(\mathbf{Id} - \bm{\pi}_0)\huu(t)$$
and one deduces that
$$\int_{t_1}^{t_2}\|(\mathbf{Id} - \bm{\pi}_0)h^1_\e(\tau)\|^2_{\H_1} \, \d\tau \lesssim \e^2 $$
and this allows to conclude to the wanted estimate. 
\end{proof}

\subsection{Hydrodynamic limit}\label{sec:hyd} We explain here the main steps in the derivation of the hydrodynamic for solutions to \eqref{eq:heps}. Our approach is a direct adaptation of the one in \cite{ALT} and we only enlighten the main changes necessary to achieve the hydrodynamic limit. 

We deduce the following convergence result (whose proof is immediate using estimates~\eqref{eq:estimh0},~\eqref{eq:estimh1} together with~\eqref{eq:hperpholder}):
\begin{prop}\label{theo:strong-conv}
There exists $\bm{h}=\bm{\pi}_{0}(\bm{h}) \in L^{2}\left((0,T);\H\right)$ such that up to extraction of a subsequence, one has
	\begin{equation}\begin{cases}\label{eq:mode-conv}
	\left\{\ho_{\e}\right\}_{\e} \text{converges to $0$ strongly  in } L^{1}((0,T)\,;{\E_1})\,, \\[0.2cm]
	\left\{\hu_{\e}\right\}_{\e} \text{converges to $\bm{h}$ weakly in } L^{2}\left((0,T)\,;\H\right)\,.
	\end{cases}\end{equation}
In particular, there exist 
	\begin{multline*}
	\varrho \in L^{2}\left((0,T)\,;\W^{m,2}_{x}(\T^{d})\right)\,, \qquad
	u \in L^{2}\left((0,T)\,;\left(\W^{m,2}_{x}(\T^{d})\right)^{d}\right)\,,  \\
	\vE \in L^{2}\left((0,T)\,;\W^{m,2}_x(\T^{d})\right)\,, 
	\end{multline*}
such that
	\begin{equation}\label{def:bmh}
	\bm{h}(t,x,v)=\left(\varrho(t,x)+u(t,x)\cdot v + \frac{1}{2}\vE(t,x)(|v|^{2}-3\en_{\star})\right)\M(v)
	\end{equation}
where $\M$ is the Maxwellian distribution introduced in \eqref{def:M}.
\end{prop}

\begin{nb}\label{nb:mode}
The convergence \eqref{eq:mode-conv} can be made even more precise since from Lemma~\ref{lem:hperpholder}, we also have
\begin{equation*}
\left\{\left(\mathbf{Id}-\bm{\pi}_{0}\right)\hu_{\e}\right\}_\e \text{ converges strongly to $0$ in } L^{2}\left((0,T)\,;\H\right).\end{equation*}
This means somehow that the only part of $h_{\e}$ which prevents the strong convergence towards $\bm{h}$ is~$\left\{\bm{\pi}_{0}\hu_{\e}\right\}_{\e}$. 
\end{nb}

\begin{nb}
Notice that the hydrodynamic quantites $(\varrho,u,\theta)$ in~\eqref{def:bmh} can be expressed in terms of~$\bm{h}$ through the following equalities:
	\begin{multline} \label{def:rho-u-theta} 
	\varrho(t,x) = \int_{\R^3} \bm{h}(t,x,v) \, \d v\,, 
	\quad
	u(t,x) = \frac{1}{\en_{\star}} \int_{\R^3} \bm{h}(t,x,v) v \, \d v\,, 
	\\
	\vE(t,x) =  \int_{\R^3} \bm{h}(t,x,v) \frac{|v|^2-3\en_{\star}}{3\en_{\star}^2}\, \d v\,. 
	\end{multline}
\end{nb}

Because of Theorem~\ref{theo:strong-conv} and for simplicity sake, from here on, we will write that our sequences converge even if it is true up to an extraction. We now aim to fully characterise the limit $\bm{h}$ obtained in Theorem \ref{theo:strong-conv}. To do so, we are going to identify the limit equations satisfied by the macroscopic quantities
$(\varrho,u,\vE)$  in \eqref{def:bmh} following the same lines as in the elastic case and more precisely the same path of {\cite{BaGoLe2,golseSR}} exploiting the fact that {the mode of convergence in Theorem \ref{theo:strong-conv} is stronger than the one of {\cite{BaGoLe2,golseSR}}}.  The regime of weak inelasticity is central in the analysis.  {The main idea is to write equations satisfied by averages in velocity of~$h_\e$ and to study the convergence of each term. To this end, we begin by a result about convergence of velocity averages of $h_{\e}$ and in what follows, we will use the following notation: for~$g=g(x,v)$, $$\langle g \rangle := \int_{\R^3} g(\cdot,v) \, \d v\,$$
which is now a function of the spatial variable only.
}

\begin{lem}\label{lem:average}
Let consider $\{h_{\e}\}_\e$ that converges to $\bm{h}$ in the sense of Theorem \ref{theo:strong-conv}. Then, for any function $\psi=\psi(v)$ such that
$|\psi(v)| \lesssim \m_{q}(v)\,$,
one has
\begin{equation}\label{eq:distr}
\la \psi\,h_{\e}\ra \xrightarrow[\e \to 0]{} \la \psi\,\bm{h}\ra \qquad \text{in} \qquad \mathscr{D}'_{t,x}\,.\end{equation}
\end{lem}

\begin{lem}\label{lem:boussi} With the notations of Theorem \ref{theo:strong-conv}, the limit $\bm{h}$ given by \eqref{def:bmh} satisfies on~$(0,T)\times \T^3$ the \emph{incompressibility condition}
\begin{equation}\label{eq:incomp}
\mathrm{div}_{x} u=0\,, 
\end{equation}
as well as \emph{Boussinesq relation}
\begin{equation}\label{eq:boussi}
\nabla_{x}\left(\varrho+\en_{\star}\vE\right)\,=0\,.
\end{equation}
As a consequence, introducing for almost every $t \in (0,T)$, 
\begin{equation} \label{def:Et}
E(t):= \int_{\T^{3}}\,\vE(t,x)\, \d x\,,
\end{equation}
one has \emph{strengthened Boussinesq relation}: for almost every $(t,x) \in (0,T)\times \T^{3}$,
\begin{equation}\label{eq:boussi2}
\varrho(t,x) + \en_{\star}\left(\vE(t,x)- E(t)\right)=0\,.
\end{equation}
\end{lem}
The proof is similar to \cite[Lemma 6.7]{ALT} and is omitted here. 
We introduce 
	\begin{equation} \label{def:bmA}
	\bm{A} (v) := v \otimes v - \frac1 3 |v|^2 \mathbf{Id} \quad \text{and} \quad p_\e := \la \frac1 3 |v|^2 h_\e \ra
	\end{equation}
so that  $
	\la v \otimes v \, h_\e \ra = \la \bm{A} \, h_\e \ra + p_\e \, \mathbf{Id}.
	$
We integrate in velocity equation~\eqref{eq:heps} that we multiply in turn by $1$, $v_{i}$, $\frac{1}{2}\,|v|^{2}$, to obtain
	\begin{subequations}\label{moments}
	\begin{equation}\label{eq:mass}
	\partial_{t}\la h_{\e}\ra +\frac{1}{\e }\mathrm{div}_{x}\la v\,h_{\e}\ra=0\,,
	\end{equation}

	\begin{equation}\label{eq:bulk}
	\partial_{t} \la v\,h_{\e}\ra +\frac{1}{\e }
	\mathrm{Div}_{x}\la \bm{A} \, h_{\e} \ra  + \frac{1}{\e } \nabla_{x} p_\e =\xi (t)\la v\,h_{\e}\ra\,,
	\end{equation}

	\begin{equation}\label{eq:energy} 
	\partial_{t}\la \tfrac{1}{2}|v|^{2}h_{\e}\ra+\frac{1}{\e}\mathrm{div}_{x}\,\la \tfrac{1}{2}|v|^{2}v\,h_{\e}\ra\,
	= \mathscr{J}_{\e}(t)+2\xi (t)\la \tfrac{1}{2}|v|^{2}h_{\e}\ra  
	\end{equation}
where
\begin{equation}\label{eq:JJE}
\mathscr{J}_{\e}(t)=\frac{1}{2}\int_{\R^{3}}|v|^{2}\left(\frac{1}{\e^2} \LLe h_\e + \frac1\e \Q_{\e,t} (h_\e,h_\e) + S_\e(t,v)
\right)\,|v|^{2}\,\d v\,.\end{equation}
	\end{subequations} 
Notice that, using \eqref{def:rho-u-theta} as well as Lemma~\ref{lem:average} and Lemma \ref{lem:boussi}, 
\begin{multline*}
\mathrm{div}_{x}\la v\,h_{\e}\ra \xrightarrow[\e \to 0]{} {\en_{\star}}\mathrm{div}_{x}u=0\,,\qquad 
\la \frac{1}{2}|v|^{2}h_{\e}\ra \xrightarrow[\e \to 0]{} \frac{3\en_{\star}}{2}\left(\varrho+\en_{\star}\vE\right)\,,\\
\nabla_{x}p_{\e} \xrightarrow[\e \to 0]{} \frac{1}{3}\nabla_{x}\la |v|^{2}\bm{h}\ra =\en_{\star}\nabla_{x}(\varrho+\en_{\star}\vE)=0\,,\\
\la \bm{A}\,h_{\e}\ra \xrightarrow[\e \to 0]{} \la \bm{A}\,\bm{h}\ra=0\,, \\
\la \frac{1}{2}|v|^{2}v_{j}\,h_{\e}\ra\,\xrightarrow[\e \to 0]{} \la \frac{1}{2}|v|^{2}v_{j}\,\bm{h}\ra=\frac{1}{2}u_{j}\la |v|^{2}v_{j}^{2}\M\ra=\frac{5}{2}\en_{\star}^{2}u_{j}\,, \qquad \forall \, j=1,2,3\,,
\end{multline*}
where all the limits hold in {$\mathscr{D}'_{t,x}$} and where $\la \bm{A}\bm{h}\ra=0$ because $\bm{h} \in \mathrm{Ker}(\mathscr{L}_{1})$ and $\bm{A}\M \in \mathrm{Ker}(\mathscr{L}_{1})^\perp$ (see Lemma~\ref{lem:phipsi}). Moreover,
$$\xi(t)\la v\,h_{\e}\ra \xrightarrow[\e \to 0]{} {\en_{\star}}\xi(t) u \qquad \text{in} \qquad {\mathscr{D}'_{t,x}}\,.$$
The limit of $ \mathscr{J}_{\re}(t)$ is handled in the following lemma.
\begin{lem}\label{lem:Jre}
It holds that  
$$ \mathscr{J}_{\e}(t) \xrightarrow[\e \to 0]{}  \mathscr{J}_{0}(t)\qquad \text{in} \qquad \mathscr{D}'_{t,x}\,,$$
where
$$
\mathscr{J}_{0}(t,x):=-6\en_{\star}\xi(t)\left( \varrho(t,x)+\en_{\star}\frac{\g+3}{4}\,\vE(t,x)\right)\,, 
\qquad (t,x) \in (0,T)\times \T^{3}\,.
$$
\end{lem}
\begin{proof} We recall from \eqref{eq:JJE} that
\begin{multline*}
\mathscr{J}_{\e}(t,x)=\frac{1}{2}\int_{\R^{3}}|v|^{2}\left(\frac{1}{\e^2} \LLe h_\e + \frac1\e \Q_{\e,t} (h_\e,h_\e) + S_\e(t,v)
\right)|v|^{2}\d v\\
=:\mathscr{J}_{\e}^{1}(t,x)+\mathscr{J}_{\e}^{2}(t,x)+\mathscr{J}_{\e}^{3}(t)\end{multline*}
where, using \eqref{eq:Dre}, we have that 
$$\mathscr{J}_{\e}^{1}(t,x)=\frac{1}{2\e^{2}}\int_{\R^{3}}|v|^{2}\LLe h_{\e}(t,x,v)\d v=- \e^{-2}\int_{\R^{3}\times\R^{3}}h_{\e}(t,x,v)\M(\vb)\bm{\Psi}_{\e,t}(|v-\vet|^{2})\,\d\vet\,\d v\,,$$
\begin{equation*}\begin{split}
\mathscr{J}_{\e}^{2}(t,x)&=\frac{1}{2\e}\int_{\R^{3}}|v|^{2}\Q_{\e,t}(h_{\e},h_{\e})(t,x,v)\,\d v\\
&=-\frac{1}{2\e}\int_{\R^{3}\times\R^{3}}h_{\e}(t,x,v)h_{\e}(t,x,\vb)\bm{\Psi}_{\e,t}(|v-\vet|^{2})\,\d\vet\,\d v\,,\end{split}\end{equation*}
and
$$\mathscr{J}_{\e}^{3}(t)=\frac{1}{2}\int_{\R^{3}}|v|^{2}{S}_{\e}(t,v)\,\d v\,,$$
is actually independent of $x$. 
One deduces from \eqref{eq:PsiE} that
$$\left|\mathscr{J}_{\e}^{2}(t,x)\right| \lesssim \e\,z^\g(t)\int_{\R^{3}\times\R^{3}}h_{\e}(t,x,v)h_{\e}(t,x,\vb)|v-\vet|^{3+\g} \,\d\vet\,\d v$$
and it is clear {from Minkowski's integral inequality} that the $\W^{m,2}_{x}(\T^{3})$ norm of the last term in the right-side is controlled by $\|h_{\e}\|_{\E}^{2}$. Thus
$$\mathscr{J}_{\e}^{2}\xrightarrow[\e \to 0]{} 0\qquad \text{in} \qquad L^{1}((0,T);\W^{m,2}_{x}(\T^{3}))\,.$$
One already saw (see \eqref{eq:Pi0Seps}) that
$$|\mathscr{J}_{\e}^{3}(t)| \lesssim \e^{\theta}z^{\overline{\theta}}(t)$$
so that 
$$\lim_{\e\to0}\mathscr{J}_{\e}^{3}(t)=0 \qquad \text{in} \qquad  L^{\infty}((0,T))\,.$$
We only have to investigate the limit of $\mathscr{J}_{\e}^{1}(t,x)$. Recalling from \eqref{eq:PsiEE} and~\eqref{def:xi} that, for any $r >0$,
$$\lim_{\e\to0}\e^{-2}\bm{\Psi}_{\e,t}(r^{2})=\frac{\mathfrak{a}_{0}}{4+\g}r^{3+\g}\lim_{\e\to0}\frac{\l_{\e}(t)^{\g}}{\e^{2}}=\frac{\mathfrak{a}_{0}}{\mathfrak{a}_{1}}\frac{r^{3+\g}}{4+\g}\xi(t)\,,$$
we deduce without major difficulty from  Theorem \ref{theo:strong-conv} that
$$\mathscr{J}_{\e}^{1}(t,x) \xrightarrow[\e \to 0]{} -\frac{\mathfrak{a}_{0}}{\mathfrak{a}_{1}}\frac{\xi(t)}{4+\g}\int_{\R^{3}\times\R^{3}}\bm{h}(t,x,v)\M(\vet)|v-\vet|^{3+\g}\,\d\vet\,\d v=:\mathscr{J}_{0}(t,x)
\quad \text{in} \quad \mathscr{D}'_{t,x}\,.$$
Now, by direct inspection from \eqref{def:bmh} 
$$
\mathscr{J}_{0}(t,x)=-\frac{\mathfrak{a}_{0}}{\mathfrak{a}_{1}}\frac{(2\en_{\star})^{\frac{3+\g}{2}}}{4+\g}K_{\g}\xi(t)\left(\varrho(t,x)+\frac{\en_{\star}}{4}\left(\frac{K_{2+\gamma}}{K_{\g}}-3\right)\vE(t,x)\right)$$
where we used, as in \textit{Step 2} of the proof of Lemma  \ref{lem:Pi0}, that  \begin{multline*}
\int_{\R^3 \times \R^3}\M(v)\M(v_{\ast})|v-v_{\ast}|^{3+\g}\, \d v_{\ast}\, \d v=(2\en_{\star})^{\frac{3+\g}{2}}K_{\g}\,,\\
\int_{\R^3 \times \R^3}\M(v)\M(v_{\ast})(|v|^{2}-3\en_{\star})|v-v_{\ast}|^{3+\g}\, \d v_{\ast}\, \d v=\frac{\en_{\star}}{2}(2\en_{\star})^{\frac{3+\g}{2}}(K_{2+\gamma}-3K_\gamma)\,.\end{multline*}
To conclude, we observe that from~\eqref{def:a1}, we have
$$\frac{\mathfrak{a}_{0}}{\mathfrak{a}_{1}}\frac{(2\en_{\star})^{\frac{3+\g}{2}}}{4+\g}K_{\g}=6\en_{\star}$$
while, from Lemma \ref{lem:eqK} $K_{2+\g}=K_{\g}(\g+6)$.   \end{proof}
Using a series of technical Lemmas as in \cite{ALT} (see Appendix \ref{app:hydro} for the main steps of the proof and the slight changes necessary with respect to \cite{ALT}), we can prove the following
\begin{prop}\label{prop:limit1}
The limit velocity $u$ in \eqref{def:bmh} satisfies
\begin{equation}\label{eq:bulk2}
\partial_{t}u-\frac{\nu_{0}}{\en_{\star}}\,\Delta_{x}u + {\en_{\star}}\mathrm{Div}_{x}\left(u\otimes u\right)+\nabla_{x}p= \xi(t)\,u
\end{equation}
while the limit temperature $\vE$ in \eqref{def:bmh} satisfies 
\begin{equation} \label{eq:temperature}
\partial_{t}\vE 
-\frac{\nu_{1}}{\en_{\star}^{2}}\,\Delta_{x}\vE + \en_{\star}\,u\cdot \nabla_{x}\vE
=\frac{2}{5\en_{\star}^{2}}\mathscr{J}_{0}(t,\cdot)+\frac{6}{5}\xi(t)\,E(t)+\frac{2}{5}\frac{\d}{\d t}E
\end{equation}
where we recall that $\mathscr{J}_0$ is defined in Lemma~\ref{lem:Jre} and $E$ is defined in~\eqref{def:Et}. 
\end{prop} 
\begin{nb}
The viscosity and heat conductivity coefficients~$\nu_{0}$ and~$\nu_{1}$ are explicit and fully determined by the elastic linearized collision operator~$\mathscr{L}_1$ (see Lemma~\ref{lem:phipsi}).  
\end{nb}
\begin{nb}
Notice also that, due to~\eqref{eq:incomp}, $\mathrm{Div}_{x}(u\otimes u)=\left(u \cdot \nabla_{x}\right)u$ and \eqref{eq:bulk2} is nothing but a {\emph{reinforced}} Navier-Stokes equation associated to a divergence-free source time-dependent term given by~$\xi(t)u$ which can be interpreted as  an energy supply/self-consistent force acting on the hydrodynamical system because of the self-similar rescaling.
\end{nb}
 \begin{prop}\label{prop:limit2} It holds that
	$$
	E(t)=0\,, \qquad \forall \, t \in [0,T]
	$$
where $E=E(t)$ is defined in~\eqref{def:Et}. Consequently, the limiting temperature $\vE$ in \eqref{def:bmh} satisfies
	\begin{equation}\label{eq:temperature2}
	\partial_{t}\,\vE(t,x)-\frac{\nu_{1}}{\en_{\star}^{2}}\,\Delta_{x}\vE(t,x) + \en_{\star}\,u(t,x)\cdot \nabla_{x}\vE(t,x)
	={3\frac{1-\g}{2}}\en_{\star}^{2}\xi(t)\,\vE(t,x)\,. 
	\end{equation}
where $\nu_{1}$ is defined in Lemma~\ref{lem:phipsi}. Moreover, the strong Boussinesq relation holds true:
	\begin{equation}\label{eq:strongBoussinesq2}
	\varrho(t,x)+\en_{\star}\vE(t,x)=0\,, \qquad \forall \, (t,x) \in [0,T] \times \T^3\,. 
	\end{equation}
\end{prop}
\begin{proof} The proof is similar to that of \cite[Prop. 6.19]{ALT} where one sees, as in the \emph{op. cit.} that
$$\dfrac{\d}{\d t}E(t)=-c_{0}\, \xi(t) \,E(t)$$
for some positive $c_{0} >0$. One proves as in \cite[Prop. 6.19]{ALT} that $E(0)=0$ and this shows the result. The limiting equation \eqref{eq:temperature2} then follows from   \eqref{eq:temperature} together with \eqref{eq:strongBoussinesq2}.\end{proof}
\section{About Haff's law}  \label{sec:orig} 
We show here how the well-posedness of \eqref{eq:heps} as well as the decay of the solution obtained in \eqref{eq:finH} provide a quite precise description of the asymptotic behaviour for the original physical problem \eqref{Bol-e} (see Theorem~\ref{theo:haff}) by carefully estimating the error between the solution $F_{\e}(t,x,v)$ to \eqref{Bol-e} and the limit $\bm{h}(t,x,v)$ defined in \eqref{def:bmh}. We keep the presentation informal, explaining only the main idea allowing to derive both \emph{global and local} versions of Haff's law for granular gases. We recall here that Haff's law as predicted in the seminal paper \cite{haff} asserts that the temperature of a freely cooling granular gases of hard-spheres decays like $(1+t)^{-\frac{2}{\g+1}}$ as $t \to \infty.$ It has been proven rigorously in the spatially homogeneous case in  \cite{SIAM}. We derive here a version of Haff's law valid for small values of $\e$ in the spatially inhomogeneous framework we adopted. 
 
In this section, we work under the conditions of Theorem~\ref{theo:main-cauchy1}. We recall that the solution $F_\e$ to \eqref{Bol-e} reads
\begin{equation}\label{eq:ScalING}
{F}_{\e}(t,x,v)= V_\e(t)^{3}f_{\e}\big(\tau_\e (t),x,V_\e(t)v\big)\end{equation}
where $f_{\e}=\M + \e h_\e$ with $h_\e$ solution to \eqref{eq:heps} and $V_\e(\cdot), \tau_\e(\cdot)$ are defined in~\eqref{eq:Ve}-\eqref{eq:taueps}. To capture the long term behaviour in Haff's law, we introduce 
 the  error term $\bm{e}_{\e}$ given by 
$$\bm{e}_{\e}(t,x,v):={V}_{\e}(t)^{3}\,\big( h_{\e}(\tau_{\e}(t),x,V_{\e}(t)v) - \bm{h}(\tau_{\e}(t),x,V_{\e}(t)v) \big)\,.$$ 
\begin{lem}\label{lemma-app-pp}
 Let $\delta \in \left(0,\frac{1}{\g+1}\right)$. For {any fixed $\e \in (0,\e_0)$} (where $\e_0$ is defined in Theorem \ref{theo:main-cauchy1}), the following estimate holds up to possibly extracting a subsequence,
\begin{equation}\label{error-term}
\big| \langle \bm{e}_{\e}(t),|v|^{\kappa}\varphi\rangle \big| \leq C(\varphi,\delta,\alpha_0)\, {\e^{-2(1-\delta)} V_\e(t)^{-\kappa-\g(1-\delta)}}\,,\qquad  \,\forall \, \varphi\in \W^{1,\infty}_vL^\infty_x(\T^3\times \R^3)\,,
\end{equation}
holds true for any $ 0\leq\kappa \leq q - 1$
where we  recall $q \geq 3+\g$ is the strength of the weight in the definition of $\E$ and $C(\varphi,\delta,\alpha_0)$ is a constant depending only on $\|\varphi\|_{\W^{1,\infty}_vL^\infty_x},\delta$ and $\alpha_0$ as defined in Theorem \ref{theo:main-cauchy1}.\end{lem}
 
\begin{proof} 
Let $\delta \in \left(0,\frac{1}{\g+1}\right)$. After a change of variables it follows that, for any suitable test-function~$\varphi$,
\begin{align*}
&\langle \bm{e}_{\e}(t),|v|^{\kappa}\varphi\rangle \\
&\quad=V_\e(t)^{-\kappa}\int_{\T^3\times \R^3} \Big( h_{\e}(\tau_\e (t),x,v) - \bm{h}(\tau_\e (t),x,v) \Big)|v|^{\kappa}\big(\varphi(x,V_\e (t)^{-1}v)-\varphi(x,0)\big)\,\d v\,\d x\\
&\qquad+V_\e(t)^{-\kappa}\int_{\T^3\times \R^3} \Big( h_{\e}(\tau_\e (t),x,v) - \bm{h}(\tau_\e(t),x,v) \Big)|v|^{\kappa}\varphi(x,0)\,\d v\,\d x\\
&\quad=:\mathcal{I}_{1}(t) + \mathcal{I}_{2}(t)\,.
\end{align*}
Note that, up to a subsequence, $\bm{h}$ is the weak$-\star$ limit of $\{h_{\e}\}_{\e}$ in $L^{\infty}\big((0,\infty);\E\big)$.  Thus, {for any~$t >0$,
$\|\bm{h}\|_{L^{\infty}((t,\infty)\,;\,\E)} \leq \liminf_{\e\searrow0}\|h_{\e}\|_{L^{\infty}((t,\infty)\,;\,\E)}.$
Consequently, for $\delta \in \left(0,\frac{1}{\g+1}\right)$   thanks to Theorem~\ref{theo:h-relaxation}, it holds that
\begin{equation}\label{tiempito}
\|\bm{h}\|_{L^{\infty}((t,\infty)\,;\,\E)} \leq C_{\delta}\left(z^\g(t)\right)^{1-\delta}\,,\qquad \forall \, t>0\,.
\end{equation}}
In regard of $\mathcal{I}_{1}(t)$, note that
$$\big| \varphi(x,V_\e(t)^{-1}v) - \varphi(x,0)\big | \leq \|\varphi\|_{\W^{1,\infty}_vL^\infty_x}\,V_\e (t)^{-1}|v|\,,$$
so that the following holds:
\begin{equation*}\begin{split}
\big|\mathcal{I}_{1}(t)\big| &\leq \|\varphi\|_{\W^{1,\infty}_vL^\infty_x}\,V_\e(t)^{-\kappa-1}\| h_{\e}(\tau_\e(t)) - \bm{h}(\tau_\e (t)) \|_{ L^{1}_{v,x}(\m_{\kappa+1}) } \\
&\leq \|\varphi\|_{\W^{1,\infty}_vL^\infty_x}\,V_\e (t)^{-\kappa-1}\| h_{\e}(\tau_\e(t)) - \bm{h}(\tau_\e (t)) \|_\E \leq C(\varphi,\delta,\alpha_{0}) \, V_\e (t)^{-\kappa-1} z(\tau_\e(t))^{\g(1-\delta)} \,.
\end{split}\end{equation*}
Observing that from~\eqref{eq:taueps}-\eqref{def:zt} and~\eqref{eq:Ve},
\begin{equation*} 
z(\tau_\e(t))= \left(\mathfrak{b}_{1}+ (1+\g)\e^{\frac{2}{\g}} \mathfrak{a}_{1}t\right)^{-\frac{1}{\g+1}}= {\e^{-\frac{2}{\g}}V_\e(t)^{-1}}\,,
\end{equation*}
we deduce that
$$
\Big|\mathcal{I}_{1}(t)\Big| \leq C(\varphi,\delta,\alpha_0) {\e^{-2(1-\delta)} V_\e(t)^{-\kappa-1-\g(1-\delta)}}\,. 
$$
The term $\mathcal{I}_{2}(t)$ is treated similarly and we obtain easily that 
\begin{align*}
\Big|\mathcal{I}_{2}(t) \Big| \leq  C(\varphi,\delta,\alpha_0) {\e^{-2(1-\delta)} V_\e(t)^{-\kappa-\g(1-\delta)}}
\end{align*}
which ends the proof.
\end{proof}

Let us explain how Lemma \ref{lemma-app-pp} allows to deduce the large time behaviour of ${F}_{\e}$ in the weak sense defined through \eqref{error-term}.
  Indeed, recalling the relations \eqref{eq:ScalING} together with Theorem \ref{theo:main} one has\begin{align*}
{F}_{\e}(t,x,v)&={V}_\e(t)^{3}f_{\e}\big(\tau_\e(t),x,V_\e(t)v\big)\\
&={V}_\e (t)^{3}\Big(\M(V_\e (t)v) + \e\,h_{\e}(\tau_\e (t),x,V_\e (t)v) \Big)\\
&={V}_\e (t)^{3}\Big( \M(V_\e(t)v) + \e\,\bm{h}(\tau_\e (t),x,V_\e (t)v)\Big) + \e\,\bm{e}_{\e}(t,x,v) \,,
\end{align*}
with
\begin{multline*}
\bm{h}(\tau_\e (t),x,V_\e (t)v)\\
=\Big(\varrho(\tau_\e (t),x)+u(\tau_\e(t),x)\cdot (V_\e(t)v) 
+ \frac{1}{2}\vE(\tau_\e(t),x)(|V_\e (t)v|^{2}-3\en_{\star})\Big)\M(V_\e (t)v)\,.
\end{multline*}
This allows to compute the temperature as defined in~\eqref{def:Teps} by
\begin{equation*}\begin{split}
{\bm{T}}_{\e}(t)&= \int_{\mathbb{T}^{3}\times\mathbb{R}^{3}}{F}_{\e}(t,x,v)|v|^{2}\, \d v\,\d x \\
&=V_\e(t)^{-2}\bigg(\int_{\R^{3}}\M(v)|v|^{2}\, \d v + \frac{\e}{2 }\int_{\mathbb{R}^{3}}\big( |v|^{2} - 3\en_{\star} \big)|v|^{2}\mathcal{M}(v)\,\d v \int_{\mathbb{T}^{3}}\vE(\tau_\e (t),x)\,\d x \\
&\phantom{++}+\e\int_{\R^{3}}|v|^{2}\M(v)\, \d v\int_{\T^{3}}\varrho(\tau_\e (t),x)\, \d x\bigg)+ { \e} \int_{\mathbb{T}^{3}\times\mathbb{R}^{3}}\bm{e}_\e(t,x,v)|v|^2\d v\,\d x\\
&= \frac{3\en_{\star}}{V_\e(t)^2}\bigg(1 + \e\left(3\en_{\star}\int_{\T^{3}}\vE(\tau (t),x)\, \d x + 2\int_{\T^{3}}\varrho(\tau(t),x)\, \d x\right)\bigg) \\
&\phantom{++++} + { \e} \int_{\mathbb{T}^{3}\times\mathbb{R}^{3}}\bm{e}_\e(t,x,v)|v|^2\,\d v\,\d x\,.
\end{split}\end{equation*}
Recalling from Lemma \ref{lem:boussi} and Proposition~\ref{prop:limit2} that
$$
\ds\int_{\T^3}\varrho(\tau_\e(t),x)\,\d x=-\en_\star\int_{\T^3}\vartheta(\tau_\e(t),x)\,\d x=0\,,
$$ 
we deduce
$${\bm{T}}_{\e}(t)= \frac{3\en_{\star}}{V_\e(t)^2} + { \e} \int_{\mathbb{T}^{3}\times\mathbb{R}^{3}}\bm{e}_\e(t,x,v)|v|^2\,\d v\,\d x\,.$$
Applying Lemma \ref{lemma-app-pp} with $\varphi=1,$ $\kappa=2$ and \eqref{eq:Ve}, we get
$$\left|\bm{T}_\e(t)-\frac{3\en_{\star}}{V_\e(t)^2}\right| \lesssim \e^{1-2(1-\delta)}V_\e(t)^{-2-\g(1-\delta)}=\e^{1+\frac{4}{\g}}\left(\mathfrak{b}_{1}+(1+\g)\mathfrak{a}_{1}\e^{\frac{2}{\g}}t\right)^{-\frac{2+\g(1-\delta)}{\g+1}}$$

which leads to an explicit expression for \textit{Haff's law} as  {for $\e\in(0,\e_0)$ fixed}, 
$$\bm{T}_{\e}(t) \simeq \frac{3\en_{\star}}{ V_\e(t)^{2}}\,, \qquad  t \to \infty\,.$$
We can actually show that the Haff's law holds uniformly \textit{locally} in space due to the boundedness of the solutions that we treat here.  This is not expected in a general context where more general solutions are considered.   Consider  $0 \leq \kappa \leq q-1$}. Note that
$$
\int_{\mathbb{R}^{3}}f_{\e}( \tau_\e(t),x,w)|w|^{\kappa}\, \d w 
= \int_{\mathbb{R}^{3}}\M(w)|w|^{\kappa}\, \d w + \e\,\int_{\mathbb{R}^{3}}h_{\e}(\tau_\e(t),x,w)|w|^{\kappa}\, \d w\,.
$$
Thanks to Sobolev embedding it holds that for any $t \geq 0$,
\begin{equation*}
\bigg| \sup_{x\in\mathbb{T}^{3}}\int_{\mathbb{R}^{3}}h_{\e}(\tau_\e (t),x,w)|w|^{\kappa}\, \d v \bigg| \leq { C} \| h_{\e}( \tau_\e (t) ) \|_{\E} \leq  { C(\alpha_{0})}\,.
\end{equation*}
Therefore, for sufficiently small $\e>0$, there exist two positive constants ${ c(\alpha_{0})}$ and  {${ \bar c(\alpha_{0})}$} such that
\begin{equation*}
{c(\alpha_{0})}\leq\int_{\mathbb{R}^{3}}f_{\e}(\tau_\e (t),x,w)|w|^{\kappa}\, \d w \leq { \bar c(\alpha_{0})}\,,\qquad \forall \, t\geq0\,,
\end{equation*}
which leads, for the physical problem, to
\begin{equation*}
V_\e(t)^{-\kappa}{ c(\alpha_{0})}\leq\int_{\mathbb{R}^{3}} {F}_{\e}(t,x,v)|v|^{\kappa}\, \d v \leq V_\e (t)^{-\kappa}{ \bar c(\alpha_{0})}\,,\qquad \forall \, t\geq0\,.
\end{equation*}
In particular, this estimate renders a local version of Haff's law
\begin{equation*}
\int_{\mathbb{R}^{d}}{F}_{\e}(t,x,v)|v|^{2}\, \d v \simeq { V_\e(t)^{-2}} \simeq \e^{\frac{4}{\g}}\left(\mathfrak{b}_{1}+(1+\g)\mathfrak{a}_{1}\e^{\frac{2}{\g}}t\right)^{-\frac{2}{\g+1}} \,, \qquad \forall \, t \geq0\,, \qquad \forall \,x \in \T^{3}\,.
\end{equation*}
This proves Theorem \ref{theo:haff}.
\begin{nb} Notice that  {Haff's law} is valid for  for any \emph{fixed} $\e\in (0,\e_0)$ but it is not uniform with respect to $\e \to 0$. This is an easy consequence of the scaling   \eqref{eq:Fescal} which implies that for $0 \leq k \leq q$,
$$\int_{\T^3\times\R^3}|v|^k F_\e(t,x,v)\,\d v\,\d x=V_\e(t)^{-k}\int_{\T^3\times\R^3}f_\e(\tau_\e(t),x,w)|w|^k\,\d w\,\d x\,,$$
so that, for any $t >0,$
$$\lim_{\e \to 0}\int_{\T^3\times\R^3}|v|^k F_\e(t,x,v)\,\d v\,\d x=0$$
where we recall that $f_\e \to \M$ and $V_\e(t)^{-1} \to 0$ as $\e \to 0$ for any fixed $t >0.$ To observe a version of Haff's law which is \emph{uniform} with respect to $\e$, we need to resort to an additional scaling of the form 
$$\widetilde{F}_\e(t,x,v)=V_1(t)^3 f_\e\left(\tau_1(t),x,vV_1(t)\right)\,, \qquad \forall\,(x,v,t) \in \T^3 \times \R^3\times \R^+\,.$$
The above argument would apply leading an uniform version of Haff's law for the rescaled~$\widetilde{F}_\e$.

\end{nb}

\begin{figure}[!h]\label{fig1}
\centering
\begin{tikzpicture}[node distance=3cm]
\node (F) [rectangle, rounded corners, minimum width=3cm, minimum height=1cm, text centered, draw=black] {$F$ solution to free cooling BE};
\node (Fe) [rectangle, rounded corners, minimum width=3cm, minimum height=1cm, , text width=4cm, text centered, draw=black, below of=F] {$F_\e$ : order $1$ time and space, highly inelastic};
\node (fe) [rectangle, rounded corners, minimum width=3cm, minimum height=1cm, text width=4cm, text centered, draw=black, below of=Fe, xshift=-4.3cm] {$f_\e$ : quasi-conservative,  close to a Maxwellian state, hydrodynamic limit observed, Haff's law observed for any \textbf{fixed} $\e >0$};
\node (tFe) [rectangle, rounded corners, minimum width=3cm, minimum height=1cm, , text width=4cm, text centered, draw=black, below of=Fe, xshift=4.3cm] {$\widetilde{F}_\e$ : order 1 inelasticity, Haff's law observed uniformly w.r.t. $\e$};
\draw [thick,->,>=stealth] (F) -- node [anchor=west] {{Diffusive} scaling: time $\sim\frac{1}{\varepsilon^2}$, space $\sim\frac{1}{\varepsilon}$} (Fe);
\draw [thick,<->,>=stealth] (Fe) -- node [anchor=east] {Self-similar scaling: $\tau_\e$, $V_\e\quad$} (fe);
\draw [thick,<->,>=stealth] (fe) -- node [anchor=north] {Self-similar scaling : $\tau_1$, $V_1$}(tFe);
\draw [thick,<->,>=stealth] (tFe) -- node [anchor=west] {$\quad\e$-singular time-velocity scaling}(Fe);
\end{tikzpicture}
\caption{Different scales and the hydrodynamic limit.}\label{strategy}
\end{figure}
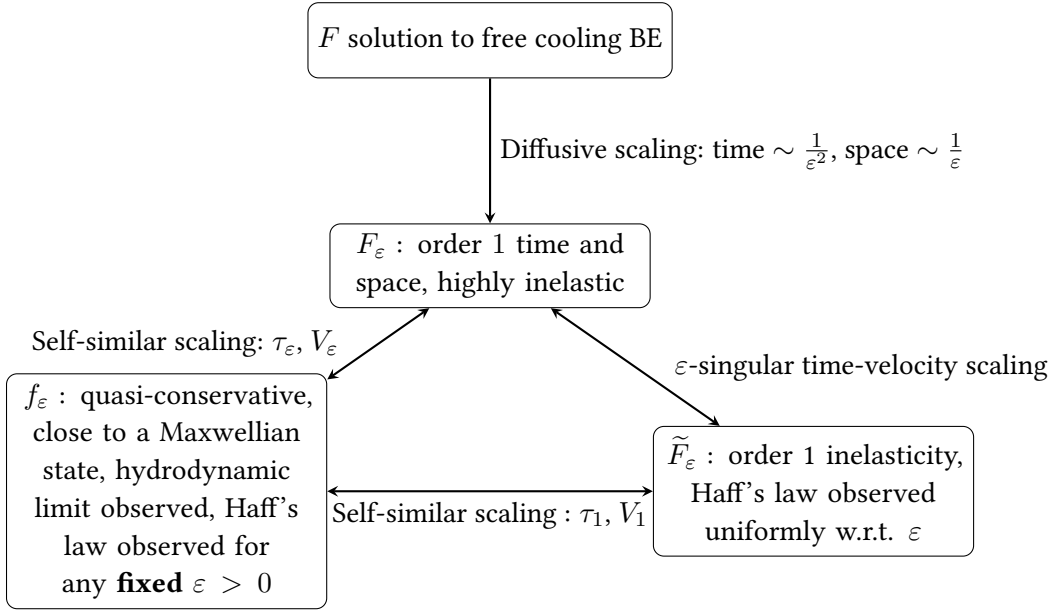

\appendix

\section{Tools for the Hydrodynamic limit}\label{app:hydro}

We establish here a general result about Gaussian distribution (in dimension $d \geq1$) \begin{lem}\label{lem:eqK} Introduce, for any $\alpha \in \R$,
$$I(\alpha):=\int_{\mathbb{R}^d} M(x)\, |x|^{3+\alpha} \left( |x|^2 - (d+2) \right) \, \d x$$
where  $M$ denotes the standard  $d$-dimensional Gaussian distribution
\[
M(x) = (2\pi)^{-\frac{d}{2}} \exp\left( -\frac{|x|^2}{2} \right), \quad x \in \mathbb{R}^d\,.
\]
Then, 
\[
I(\alpha) > 0 \quad \Longleftrightarrow \quad \alpha > -1\,.
\]
\end{lem}
\begin{proof} Computing $I(\alpha)$ using spherical coordinates, we obtain
\[
\int_{\mathbb{R}^d} M(x)\, |x|^{3+\alpha} \left( |x|^2 - (d+2) \right) \d x = |\S^{d-1}|(2\pi)^{-\frac{d}{2}} \int_0^\infty e^{-\frac{r^2}{2}} r^{d + 2 + \alpha} \left( r^2 - (d+2) \right) \d r\,.
\]  The result is then reduced to determine the sign of the radial  integral:
\begin{equation*} 
J(\alpha):= \int_0^\infty e^{-\frac{r^2}{2}} r^{d + 2 + \alpha} (r^2 - (d+2)) \, \d r\,.\end{equation*}
These integral can be evaluated using the substitution $u = \frac{r^2}{2}$ and the Gamma function identity:
\[
\int_0^\infty r^k e^{-\frac{r^2}{2}} \, dr = 2^{\frac{k-1}{2}} \Gamma\left( \frac{k+1}{2} \right), \quad \forall \, k > -1\,.
\]
With this identity,
\begin{equation*}\begin{split}
J(\alpha) &= 2^{\frac{d+3+\alpha}{2}} \Gamma\left( \frac{d+5+\alpha}{2} \right) 
- (d+2) \cdot 2^{\frac{d+1+\alpha}{2}} \Gamma\left( \frac{d+3+\alpha}{2} \right)\,\\
&= 2^{\frac{d+1+\alpha}{2}} \left[ 2 \cdot \Gamma\left( \frac{d+5+\alpha}{2} \right) 
- (d+2) \cdot \Gamma\left( \frac{d+3+\alpha}{2} \right) \right]\,,
\end{split}\end{equation*}
and  using the Gamma function identity $\Gamma(z+1) = z\Gamma(z)$ we obtain
$$
J(\alpha) = 2^{\frac{d+1+\alpha}{2}} \cdot \Gamma\left( \frac{d+3+\alpha}{2} \right) \cdot (\alpha + 1)\,,$$ and get the conclusion. 
\end{proof}
We now recall some well-known properties of the linearized elastic operator $\mathscr{L}_{1}$ useful for the  hydrodynamic limit as well as the main steps used for the proof of the results in Section \ref{sec:hyd}. 
 \begin{lem}\label{lem:phipsi}  Let $\bm{A}$  be the traceless tensor defined in~\eqref{def:bmA} and let $\bm{b}$ be the vector defined by
 \begin{equation}\label{def:bmb}
 \bm{b}(v)=\frac{1}{2}\left(|v|^{2}-5\en_{\star}\right) v\,, \qquad v \in\R^{3}\,.\end{equation}
{One has that $\bm{A}\M$, $\bm{b}\M \in \left(\mathrm{Ker}(\mathscr{L}_{1})\right)^\perp$ in $L^2_v(\M^{-\frac12})$} and there exists two radial functions $\chi_{i}=\chi_{i}(|v|)$, $i=1,2$, such that
$$\phi(v)=\chi_{1}(|v|)\bm{A}(v) \in \mathscr{M}_{3}(\R) \qquad \text{and} \qquad \psi(v)=\chi_{2}(|v|)\bm{b}(v) \in \R^{3}\,,$$
satisfy
\begin{equation}\label{eq:L1chi}
\mathscr{L}_{1}(\phi\,\M)=-\bm{A}\,\M\,,\qquad \mathscr{L}_{1}(\psi\,\M)=-\bm{b}\,\M\,.
\end{equation}
Moreover,
\begin{multline}\label{eq:visco}
\la \phi^{i,j}\mathscr{L}_{1}(\phi^{k,\ell}\M)\ra=-{\nu}_{0}\left(\delta_{ik}\delta_{j\ell}+\delta_{i\ell}\delta_{jk}-\frac{2}{3}\delta_{ij}\delta_{kl}\right)\\
\la \psi_{i}\mathscr{L}_{1}(\psi_{j}\M)\ra=-\frac{5}{2}\nu_{1}\,\delta_{ij}\,, \qquad \forall \, i,j,k,\ell \in \{1,2,3\}\,,\end{multline}
with
$${\nu}_{0}:=-\frac{1}{10}\la \phi\,:\,\mathscr{L}_{1}(\phi\M)\ra \geq 0\,, \qquad \nu_{1}:=-\frac{2}{15}\la \psi \cdot \mathscr{L}_{1}(\psi\M)\ra \geq 0\,.$$
{Finally, 
$$|\phi^{i,j}(v)| \lesssim \m_{3}(v)\,, \qquad |\psi_{i}(v)| \lesssim \m_{4}(v)\,, \qquad \forall \, i,j \in \{1,\ldots,d\}\,.$$}
\end{lem}
For the proof of convergence in Theorem \ref{theo:main}, we follow the approach of \cite[Section 7]{ALT}. We only here recall some of main steps and necessary changes in order to adapt the proof of \cite{ALT}.  The key aspect of the proof of Proposition \ref{prop:limit1}  is the careful study of the equations of motion 
$$ {u}_{\e}(t,x):={\frac{1}{\en_{\star}}}\int_{\R^{3}}v\,h_{\e}(t,x,v)\, \d v\,,$$
and temperature
\begin{equation} \label{def:bmthetaeps}
\bm{\theta}_{\e}(t,x):=\frac{1}{2}\int_{\R^{3}}\big( |v|^{2}-5\en_{\star} \big)h_{\e}(t,x,v)\,\d v\,.
\end{equation}
Notice that \eqref{eq:bulk} reads
$$\partial_{t} u_{\e}(t,x)  + \frac{1}{\e \en_{\star} }
	\mathrm{Div}_{x}\la \bm{A} \, h_{\e} \ra  + \frac{1}{\e \en_{\star}} \nabla_{x} p_\e =\xi (t)\,u_{\e}(t,x)$$
where the gradient of pressure can be removed by applying the Leray projector $\mathcal{P}$, 
	whereas
\begin{equation}\label{eq:trick}\partial_{t}\bm{\theta}_{\e}(t,x)+\frac{1}{\e}\mathrm{div}_{x}\la \bm{b}\,h_{\e}\ra=\mathscr{J}_{\e}(t)+2\xi(t)\la \frac{1}{2}|v|^{2}h_{\e}\ra\,.\end{equation}
 With this, one has the following whose proof follows that of \cite[Lemma 6.10]{ALT} (noticing  a change in the definition of $\bm{u}_{\e}$ with respect to \cite[Lemma 6.10]{ALT}):

\begin{lem}\label{prop:strongue} Introduce for $(t,x) \in (0,T) \times \T^3$:
\begin{equation} \label{def:bmueps-bmthetaeps}
{\bm{u}_\e(t,x):=\exp\left(-\int_{0}^{t}\xi(s)\,\d s\right)\mathcal{P}{u}_{\e}(t,x)}
\end{equation}
and recall that $\bm{\theta}_{\e}$ is defined in~\eqref{def:bmthetaeps}.
Then, {$\{\partial_{t}\bm{u}_{\e}\}_{\e}$ and $\{\partial_{t}\bm{\theta}_{\e}\}_{\e}$ are bounded in the space $L^{1}\left((0,T)\,; {\W^{m,2}_{x}(\T^{d})}\right)$.}
Consequently, up to the extraction of a subsequence,
\begin{equation}\label{eq:Puest}
\int_{0}^{T}\left\|\mathcal{P}{u}_{\e}(t)-u(t)\right\|_{ {\W^{m-1,2}_{x}(\T^{d})}}\, \d t\xrightarrow[\e \to 0]{}0\end{equation}
and
\begin{equation}\label{eq:varthstr}
\int_{0}^{T}\left\|\bm{\theta}_{\e}(t,\cdot)-\bm{\theta}_{0}(t,\cdot)\right\|_{ {\W^{m-1,2}_{x}(\T^{d})}}\, \d t\xrightarrow[\e \to 0]{}0\end{equation}
where 
\begin{equation} \label{def:theta0}
\bm{\theta}_{0}(t,x):=\la\bm{b}\,\bm{h}\ra=\frac{3\en_{\star}}{2}\left(\varrho(t,x)+\en_{\star}\vE(t,x)\right)-\frac{5}{2}\en_{\star}\varrho(t,x)\,.
\end{equation} 
In other words, $\{\mathcal{P}{u}_{\e}\}_{\e}$ (resp. $\{\bm{\theta}_{\e}\}_{\e}$) converges strongly to $u=\mathcal{P} u$ (resp. $\bm{\theta}_{0}$) in the space $L^{1}\left((0,T)\,;{\W^{m-1,2}_{x}(\T^{d})}\right)$.
\end{lem}

We can now establish the following whose proof is a simple adaptation of \cite[Lemma 6.12 \& Lemma 6.13]{ALT}.
\begin{lem}\label{lem:limitA} In the distributional sense, 
\begin{equation}\label{eq:Ah}
\mathcal{P}\mathrm{Div}_{x}\left(\la \frac1\e\bm{A}\,h_{\e}\ra-\la \phi \,\Q_{1}\left(\bm{\pi}_{0}h_{\e},\bm{\pi}_{0}h_{\e}\right)\ra\right)\xrightarrow[\e \to 0]{} -{\nu_{0}}\,\Delta_{x}u\end{equation}
and
\begin{equation}\label{eq:limpsiL1}
\frac1\e\mathrm{div}_{x}\la \bm{b}\,h_{\e}\ra + \mathrm{div}_{x} \la \psi\,\Q_{1}(\bm{\pi}_{0}h_{\e},\bm{\pi}_{0}h_{\e}\ra
\xrightarrow[\e \to 0]{}-\frac{5}{2}\,\nu_{1}\,\Delta_{x}\vE
\end{equation}
where ${\nu}_{0},\nu_{1}$ and the Burnet functions $\phi,\psi$ are defined in Lemma~\ref{lem:phipsi}. 
 \end{lem} 
 Moreover, the following
\begin{lem}\label{lem:convect}
We have
$$\mathcal{P}\mathrm{Div}_{x}\la \phi\,\Q_{1}\left(\bm{\pi}_{0}h_{\e},\bm{\pi}_{0}h_{\e}\right)\ra\xrightarrow[\e \to 0]{}\en_{\star}^{2} \mathcal{P}\mathrm{Div}_{x}(u\otimes u)
 \qquad \text{in} \qquad \mathscr{D}'_{t,x}$$
and
$$\mathrm{div}_{x}\la \psi\,\Q_{1}(\bm{\pi}_{0}h_{\e},\bm{\pi}_{0}h_{\e})\ra\xrightarrow[\e \to 0]{}\frac{5}{2}\en_{\star}^{3}u \cdot \nabla_{x}\vE
 \qquad \text{in} \qquad \mathscr{D}'_{t,x}\,.$$
In particular,
\begin{equation}\label{eq:limphiL1-2}
\mathcal{P}\mathrm{Div}_{x}\la \frac1\e\bm{A}\,h_{\e}\ra
\xrightarrow[\e \to 0]{}
-\nu_{0}\Delta_{x}u + \en_{\star}^{2}\mathcal{P}\mathrm{Div}_{x}(u \otimes u) \qquad \text{in} \qquad \mathscr{D}'_{t,x}\end{equation}
while
\begin{equation}\label{eq:limpsiL1-2}
\mathrm{div}_{x}\la \frac1\e\bm{b}\,h_{\e}\ra
\xrightarrow[\e \to 0]{}
-\frac{5}{2}\left(\nu_{1}\,\Delta_{x}\vE - \en_{\star}^{3}u \cdot \nabla_{x}\vE\right) \qquad \text{in} \qquad \mathscr{D}'_{t,x}
\end{equation}
where $\nu_{0}$ and $\nu_{1}$ are defined in Lemma~\ref{lem:phipsi}.
\end{lem}
\begin{nb} The proof is exactly the same as that of Lemma \cite[Lemma 6.14]{ALT} inspired by \cite[Corollary~5.7]{golseSR}. We only point out here that, with notations of the proof of \cite{ALT}, we need to replace Eqs. (6.40) and (6.41) respectively with
\begin{equation*}\label{eq:bulkP}
\e\,\partial_{t}{u}_{\e} + \nabla_{x}\bm{\beta}_{\e}=\xi(t){u}_{\e} - \frac{1}{\en_{\star}}\mathrm{Div}_{x}\la \bm{A}\,h_{\e}\ra\end{equation*}
and
\begin{equation*}\label{energy21}
\e\partial_{t}\bm{\beta}_{\e}+\mathrm{div}_{x}\,\la \frac{1}{3\en_{\star}}|v|^{2}v\,h_{\e}\ra\,=\frac{2\e}{3\en_{\star}}\mathscr{J}_{\e}(t)+2\e\xi(t)\bm{\beta}_{\e}\,
\end{equation*}
We notice that the time-dependent factor $\xi(t)$ plays no role in the estimate here since $\xi \in L^{\infty}([0,\infty))$.
We also replace the definition of $\bm{F}_{\e}$ and $\bm{G}_{\e}$ as
\begin{multline*}
\bm{F}_{\e}:=\e\xi(t)\nabla_{x}\bm{U}_{\e}-\frac{1}{\en_{\star}}(\mathbf{Id}-\mathcal{P})\mathrm{Div}_{x}\la \bm{A}\,h_{\e}\ra\\
\bm{G}_{\e}:=-\frac{2}{3\en_{\star}}\mathrm{div}_{x}\la \bm{b}\,h_{\e}\ra+\frac{2\e}{3\en_{\star}}\mathscr{J}_{\e}(t) +2\e\xi(t)\bm{\beta}_{\e}\,,
\end{multline*}
and one checks again that 
$$
\|\bm{F}_{\e}\|_{L^{1}((0,T)\,;\, {\W^{m-1,2}_{x}(\T^{3})})} \lesssim \e \qquad \text{and} \qquad \|\bm{G}_{\e}\|_{L^{1}((0,T)\,;\, {\W^{m-1,2}_{x}(\T^{3})})} \lesssim \e\,.
$$
\end{nb}
\begin{proof}[Proof of Proposition~\ref{prop:limit1}] The proof of \eqref{eq:bulk2} is a straightforward consequence of the previous lemma. To investigate the evolution of $\vE$, we recall that $\bm{\theta}_{\e}$ satisfies \eqref{eq:trick}. 
We notice that
$$\frac{1}{\e^{3}}\mathscr{J}_{\e}(t)+2\xi(t)\la \frac{1}{2}|v|^{2}h_{\e}\ra \xrightarrow[\e \to 0]{} \mathscr{J}_{0}+3\en_{\star}\xi(t)\left(\varrho+\en_{\star}\vE\right)\,,$$
whereas
$$\bm{\theta}_{\e}\xrightarrow[\e \to 0]{} \la \bm{b}\bm{h}\ra=\frac{3\en_{\star}}{2}\left(\varrho+\en_{\star}\vE\right)-\frac{5}{2}\en_{\star}\varrho
\qquad \text{in} \qquad \mathscr{D}'_{t,x}\,.$$
We deduce from \eqref{eq:limpsiL1-2}, performing the distributional limit of \eqref{eq:trick}, that
\begin{multline}\label{eq:energy2}
\frac{3\en_{\star}}{2}\partial_{t}\left(\varrho+\en_{\star}\vE\right)-\frac{5}{2}\en_{\star}\partial_{t}\varrho
-\frac{5}{2}\nu_{1}\,\Delta_{x}\vE{+}\frac{5}{2}\en_{\star}^{3}\,u\cdot \nabla_{x}\vE
=\mathscr{J}_{0}+3\en_{\star}\xi(t)\left(\varrho+\en_{\star}\vE\right).
\end{multline}
Using the strengthened Boussinesq relation \eqref{eq:boussi2}, we see that 
$$\partial_{t}\left(\varrho+\en_{\star}\vE\right)=\en_{\star}\frac{\d}{\d t}E \qquad \text{and} \qquad \partial_{t}\varrho=-\en_{\star}\left(\partial_{t}\vE-\frac{\d}{\d t}E\right)\,,$$ 
and get the result.
\end{proof}
As in \cite[Lemma 6.5]{ALT} an important consequence of Ascoli-Arzela Theorem is the following
{\begin{lem} \label{lem:cont}
Consider the sequences $\{\bm{u}_\e\}_{\e}$ and~$\{\bm{\theta}_\e\}_\e$ defined in Lemma~\ref{prop:strongue}. 
The time-depending mappings 
$$
t \in [0,T] \longmapsto   \left\|\bm{u}_{\e}(t)\right\|_{ {\W^{m-1,2}_{x}(\T^{3})}}
\quad \text{and} \quad 
t \in [0,T] \longmapsto  \left\|\bm{\theta}_{\e}(t)\right\|_{ {\W^{m-1,2}_{x}(\T^{3})}}
$$
are H\"older continuous uniformly in $\e$. 
As a consequence, the limiting quantities $u$ and $\bm{\theta}_0$ belong to~$\mathcal{C}([0,T]\,;\W^{m-1,2}_x(\T^3))$.
\end{lem}} 

\section{Estimates on the collision operator} \label{app:collL}
We collect here some nonlinear estimates on the collision operator $\Q_{\re}$ associated to non constant restitution coefficient. 
We start with a technical lemma.
\begin{lem}\label{lem:Lam} For a restitution coefficient $\re(\cdot)$ satisfying Definition \ref{defiC}, the mapping defined as
$$\Lambda_{\re}(u):=\int_{\S^{2}}\frac{1}{\re^{2}(|u\cdot n|)J_{\re}(|u\cdot n|)}\exp\left(-\frac{1-\re^{2}(|u\cdot n|)}{2\en_{\star}}\left(u\cdot n\right)^{2}\right)\d n\,, \qquad u \in \R^{3}$$
is bounded by a constant depending on $\re(\cdot)$ and $\en_\star$. 
\end{lem}
\begin{proof} Namely, as in \cite{A}, one defines, for any $\beta >0$ the mapping
$$\psi_{\re,\beta}(z)=\frac{1}{\re^{2}(z)J_{\re}(z)}\exp\left(-\frac{\beta}{2}(1-\re^{2}(z))z^{2}\right), \qquad z >0$$
and set
$${\varphi}_{\re,\beta}(y)=\int_{0}^{1}\psi_{\re,\beta}(yz)\,\d z\,, \qquad y >0\,.$$
As simple computation using polar coordinates shows that
$$\Lambda_{\re}(u)=2|\S^{1}|\varphi_{\re,\beta}(|u|)$$
for the choice $\beta=\frac{1}{\en_{\star}}.$ We are therefore going to compute $\|\varphi_{\re,\beta}\|_{\infty}$ for any $\beta >0.$ Recall that,  since $\re(\cdot)$ is non-increasing with $\re(0)=1$, there is a uniquely determined $z_{1} >0$ such that
$$\re(r) \leq \frac{1}{2} \qquad \forall \,r >z_{1}\,, \qquad \re(r) > \frac{1}{2} \qquad \forall \,r \leq z_{1}\,.$$
Recalling that from~\eqref{eq:boundJac}, there is $\alpha >0$ such that $J_{\re}(r) \gtrsim \re(r)^{\alpha}$, 
we thus get that
$$\psi_{\re,\beta}(z) \lesssim (\re(z))^{-2-\alpha} \lesssim 1\,, \qquad \forall \,z \leq z_{1}$$
whereas, for $z >z_{1}$, since $1-\re^{2}(z) \geq \frac{1}{2}$, one thus has
$$\psi_{\re,\beta}(z) \lesssim \re(z)^{-2-\alpha}\exp\left(-\frac{\beta}{4}z^{2}\right)\,, \qquad \forall \,z > z_{1}\,.$$
Since $\eta_{\re}(r)=r\re(r)$ is nondecreasing, it holds that $\re(r) \geq \frac{\eta_{\re}(z_{1})}{r}$ for $r >z_{1}$ and 
$$\psi_{\re,\beta}(z) \lesssim z^{2+\alpha}\exp\left(-\frac{\beta}{4}z^{2}\right)\,, \qquad \forall \,z >z_{1}\,.$$
Since 
$$\varphi_{\re,\beta}(y)=\int_{0}^{1}\psi_{\re,\beta}(yz)\,\d z=\frac{1}{y}\int_{0}^{y}\psi_{\re,\beta}(z)\,\d z\,, \qquad \forall \, y >0\,,$$
one sees that $\sup_{y \in [0,z_{1}]}\varphi_{\re,\beta}(y) \lesssim 1$ while
$$\varphi_{\re,\beta}(y) \lesssim 1 + \int_{z_{1}}^{\infty}z^{2+\alpha}\exp\left(-\frac{\beta}{4}z^{2}\right)\d z \qquad \forall \,y > z_{1}\,.$$
In particular, we obtain  
\begin{equation}\label{eq:varphire}
\|\varphi_{\re,\beta}\|_{\infty} \lesssim 1+\beta^{-\frac{3+\alpha}{2}}\,,\end{equation} 
which proves the Lemma.
\end{proof}
 
Lemma~\ref{theo:briant2.4} allows us to prove the following trilinear estimates on $\Q_\re$: 
\begin{lem}\label{lem:trilinear}
For any  $f,g,\varphi \in \H_1$, it holds
\begin{equation}\label{eq:trilinear}
\left\langle \Q^{\pm}_{\re}(g,f),\varphi\right\rangle_{\H} \lesssim \|\varphi\|_{\H_{1}}\,\left(\|f\|_{\H_{1}}\|g\|_{\H}
+\|f\|_{\H}\|g\|_{\H_{1}}\right)
\end{equation}
where the functional spaces $\H$ and $\H_1$ are defined in~\eqref{eq:defH}. 
\end{lem}
\begin{proof}  Since $\W^{\ell,2}_{x}(\T^{3})$ is a Banach algebra for $\ell > \frac{3}{2}$ and $\Q_{\re}^\pm$ is local in $x$, it is easy to see that it is enough to prove the result for function $f,g,\varphi$ not depending on $x$. We use the strong form of $\Q_{\re}^{\pm}(g,f)$ given in \eqref{Boltstrong} with 
$$B_{0}(u,n)=|u|b_{0}(\widehat{u}\cdot \n)$$
and prove the result only for $\Q^{+}_{\re}(g,f)$, the proof for $\Q^{-}_{\re}(g,f)=\Q_{1}^{-}(g,f)$ being simpler and well-known. 
We set
\begin{equation}\label{eq:GFPHI}
F:=\M^{-\frac{1}{2}}f\,,\qquad G:=\M^{-\frac{1}{2}}g\,, \qquad \Phi:=\M^{-\frac{1}{2}}\varphi\,.
\end{equation}
Given $v \in\R^{3}$, one has  
\begin{multline*}
\left|\Q_{\re}^{+}\big(g,f\big)(v)\right| \leq  \int_{\R^{3}\times\S^{2}}\dfrac{B_0(u,\n)}{\re\big(|'u\cdot n|\big)J_{\re}\big(|'u\cdot n|\big)}\,|f('\!v)|\,|g('\!\vb)|\, \d n\, \d\vb\,,\\
\leq  \int_{\R^{3}\times\S^{2}}\dfrac{B_0(u,\n)}{\re\big(|'u\cdot n|\big)J_{\re}\big(|'u\cdot n|\big)}\,\M^{\frac{1}{2}}('\!v)\M^{\frac{1}{2}}('\!\vb)\,|F('\!v)|\,|G('\!\vb)|\, \d n\, \d\vb\,,
\end{multline*}
where, according to \eqref{energ}
$$\M^{\frac{1}{2}}('\!v)\M^{\frac{1}{2}}('\!\vb)=\M^{\frac{1}{2}}(v)\M^{\frac{1}{2}}(\vb)\exp\left(-\frac{1-\re^{2}(|'\!u\cdot n|)}{4\en_{\star}}\left('u\!\cdot n\right)^{2}\right)$$
so that
$$\left|\Q_{\re}^{+}\big(g,f\big)(v)\right| \leq   \M^{\frac{1}{2}}(v)\int_{\R^{3}\times\S^{2}}\M^{\frac{1}{2}}(\vb)b_{0}(\widehat{u}\cdot \n)\mathcal{K}(v,\vb,n) \frac{|F('\!v)|\,|G('\!\vb)|}{\sqrt{J_{\re}\big(|'u\cdot n|\big)}}\,  \d n\, \d\vb\,$$
with 
$$\mathcal{K}(v,\vb,n)=\dfrac{|u|}{\re\big(|'\!u\cdot n|\big)\sqrt{J_{\re}\big(|'u\cdot n|\big)}}\,\exp\left(-\frac{1-\re^{2}\big(|'u\cdot n|\big)}{4\en_{\star}}\left('u\!\cdot n\right)^{2}\right)\,.$$
Thanks to Cauchy-Schwarz inequality, we deduce that
\begin{multline*}
\left|\Q_{\re}^{+}\big(g,f\big)(v)\right| \leq  \M^{\frac{1}{2}}(v)\left(\int_{\R^{3}\times\S^{2}}\,|F('\!v)|^{2}\,|G('\!\vb)|^{2}\, \frac{b_{0}^{2}(\widehat{u}\cdot n)}{J_{\re}\big(|'u\cdot n|\big)}\, \d \n\,\d\vb\right)^{\frac{1}{2}}\\
\left(\int_{\R^{3}\times\S^{2}}\left(\mathcal{K}(v,\vb,n)\right)^{2}\M(\vb)\, \d \n\, \d \vb\right)^{\frac{1}{2}}
\end{multline*}
Now, for fixed $(v,\vb) \in \R^{3}\times\R^{3}$, we compute
\begin{align*}
&\int_{\S^{2}}\left(\mathcal{K}(v,\vb,n)\right)^{2}\,\d n\\
&\quad = \int_{\S^{2}}\frac{|u|^{2}}{\re^{2}\big(|'u\cdot n|\big)J_{\re}\big(|'u\cdot n|\big)}\exp\left(-\frac{1-\re^{2}\big(|'u\cdot n|\big)}{2\en_{\star}}\left('u\!\cdot n\right)^{2}\right)\d n
=: |u|^2 \Lambda_\re('\!u)\,.
\end{align*}
One can prove (see subsequent Lemma \ref{lem:Lam} for a full proof) that
$\Lambda_{\re}(\cdot) \in L^{\infty}$.
Thus, 
$$\left(\int_{\R^{3}\times\S^{2}}\left(\mathcal{K}(v,\vb,n)\right)^{2}\M(\vb)\, \d \vb\, \d \n\right)^{\frac{1}{2}}
\lesssim \left(\int_{\R^{3}}|v-\vet|^{2}\M(\vb)\,\d\vb \right)^{\frac{1}{2}} \lesssim \langle v\rangle.$$
Therefore
$$\left|\Q_{\re}^{+}\big(g,f\big)(v)\right| \lesssim  \M^{\frac{1}{2}}(v)\langle v\rangle\,\left(\int_{\R^{3}\times\S^{2}}\,b_{0}^{2}(\widehat{u}\cdot n)\frac{|F('\!v)|^{2}\,|G('\!\vb)|^{2}}{J_{\re}\big(|'u\cdot n|\big)}\, \d \n\, \d\vb\right)^{\frac{1}{2}}.$$
Mutiplying by $\varphi(v)\M^{-1}(v)$ and integrating over $\R^{3}$, we get
\begin{align*}
&\left\langle \Q_{\re}^{+}(g,f),\varphi\right\rangle_{L^{2}_{v}(\M^{-\frac{1}{2}})} \\
&\qquad \lesssim  \int_{\R^{3}}\Phi(v)\langle v\rangle\left(\int_{\R^{3}\times\S^{2}}\,b_{0}^{2}(\widehat{u}\cdot n)\,|F('\!v)|^{2}\,|G('\!\vb)|^{2}\, \frac{1}{J_{\re}\big(|'u\cdot n|\big)}\,\d \n\, \d\vb\right)^{\frac{1}{2}}\, \d v\,.
\end{align*}
Using Cauchy-Schwarz inequality again, we obtain
\begin{equation*}\begin{split}
&\left\langle \Q_{\re}^{+}(g,f),\varphi\right\rangle_{L^{2}_{v}(\M^{-\frac{1}{2}})} \\
&\qquad \lesssim \|  \langle\cdot\rangle^{\frac{1}{2}}\Phi\|_{L^{2}_{v}}\left(\int_{\R^{3}}\langle v\rangle \, \d v\int_{\R^{3}\times\S^{2}}\,|F('\!v)|^{2}\,|G('\!\vb)|^{2}\, \frac{b_{0}^{2}(\widehat{u}\cdot n)}{J_{\re}\big(|'u\cdot n|\big)}\,\d \n\,\, \d\vb\right)^{\frac{1}{2}}\\
&\qquad \lesssim \| \langle\cdot\rangle^{\frac{1}{2}}\Phi\|_{L^{2}_{v}}\left(\int_{\R^{3}\times \R^{3}\times\S^{2}}\langle v'\rangle\,|F(v)|^{2}\,|G(\vb)|^{2}\,b_{0}^{2}(\widehat{u'}\cdot n)\, \d n\, \d\vb \,\d v\right)^{\frac{1}{2}}\,,
\end{split}\end{equation*}
where we used the pre-post collisional change of variable in the last estimate. Since
$$\langle v'\rangle \lesssim \langle v\rangle + \langle \vb\rangle$$
we deduce that
$$\left\langle \Q_{\re}^{+}(g,f),\varphi\right\rangle_{L^{2}_{v}(\M^{-\frac{1}{2}})} \lesssim \| \langle\cdot\rangle^{\frac{1}{2}}\Phi\|_{L^{2}_{v}}\left(\|F\|_{L^{2}_{v}}\|\langle\cdot\rangle^{\frac{1}{2}} G\|_{L^{2}_{v}}+\|\langle\cdot\rangle^{\frac{1}{2}}\,F\|_{L^{2}_{v}}\|G\|_{L^{2}_{v}}\right)\,,$$
or equivalently
\begin{multline*}\label{eq:tril}
\left\langle \Q^{\pm}_{\re}(g,f),\varphi\right\rangle_{L^{2}_{v}(\M^{-\frac{1}{2}})} \leq C\|\varphi\|_{L^{2}_{v}(\M^{-\frac{1}{2}}\langle\cdot\rangle^{\frac{1}{2}})}\,\bigg(\|f\|_{L^{2}_{v}(\M^{-\frac{1}{2}}\langle \cdot\rangle^{\frac{1}{2}})}\|g\|_{L^{2}_{v}(\M^{-\frac{1}{2}})}\\
+\|f\|_{L^{2}_{v}(\M^{-\frac{1}{2}})}\|g\|_{L^{2}_{v}(\M^{-\frac{1}{2}}\langle\cdot\rangle^{\frac{1}{2}})}\bigg)\,,\end{multline*}
which ends the proof.
\end{proof}

\begin{nb}\label{nb:Lambdae} Replacing $\re(\cdot)$ with $\re_{\l}(\cdot)$ with$\l \lesssim 1$, we see that the constants involved in~\eqref{eq:trilinear} will depend on $\l$ only through  $\|\Lambda_{\re_{\l}}\|_{\infty}$ which, according to Lemma~\ref{lem:Lam}, is proportional to $\|\varphi_{\re_{\l},\beta}\|_{\infty}.$ Observing that, for $\l \lesssim 1$, one has $J_{\re_{\l}}(z)=J_{\re}(\l z),$ so that
$$\psi_{\re,\beta \ell^{-2}}(\l z)=\psi_{\re_{\l},\beta}(z)\,.$$
In particular, $\|\varphi_{\re_{\l},\beta}\|_{\infty}= {\sup_{y}|\varphi_{\re,\beta \l^{-2}}(\l\,y)|}=\|\varphi_{\re,\beta{\l^{-2}}}\|_{\infty}$
and one sees from \eqref{eq:varphire} that
$$\|\varphi_{\re_{\l},\beta}\|_{\infty} \lesssim 1+ \beta^{-\frac{3+\alpha}{2}} \ell^{3+\alpha}
\lesssim 1+ \beta^{-\frac{3+\alpha}{2}} $$
since $\l \lesssim 1$. This shows that, replacing $\re(\cdot)$ with $\re_{\l}(\cdot)$ with $\l \lesssim 1$, the constants involved in Lemma~\ref{lem:trilinear} can be made \emph{independent} of $\l$.\end{nb}

We now collect several results derived in previous contributions regarding convolution inequalities for the collision operator, see \cite[Proposition A.2, Remark A.3]{CMS}

\begin{prop}\label{propB:CMS} Let $\re(\cdot)$ be a given restitution coefficient belonging to the class $\mathcal{R}_\g$ and let
$${B}(u,\widehat{u}\cdot \sigma)=\Theta(|u|)\,b(\widehat{u}\cdot\sigma)$$
be a given collision kernel with $\Theta(r) \geq 0$ and $b(s)=b(-s)$ with $\mathrm{Supp}\:b \in [-1+\delta,1-\delta]$ for some $\delta >0$. Let $\Q^+_{{B},\re}$ and $\Q^+_{B,1}$ denote  the positive part of the collision operator associated to ${B}$ with restitution coefficient $\re(\cdot)$ and elastic interactions respectively. Then, there exists an explicit constant $C_0(\re) >0$ such that
\begin{multline*}
\left|\IR \left[\Q^+_{{B},e}(f,g)-\Q^+_{{B},1}(f,g)\right]\psi\,\d v\right| \leq C_0(\re)   \IR \Q^+_{ {B}_\gamma,1}(f , g )\,|\psi(v)|\,\d v \\
 + 2^{\gamma+6}\sup_{r >0} \dfrac{1-\re(r)}{r^\gamma}\int_0^1\d s\int_{\R^3} \Q^+_{{\overline {B}_\gamma},\tilde{\re}_s}\left(f, h\right)|\psi(v)|\,\d v \end{multline*}
where $h(v)=g(v)+|\nabla g(v)|$, the kernels ${B}_\gamma$ and ${\overline{B}_\gamma}$ are given by
$${B}_\gamma(u,\widehat{u} \cdot \sigma)= {B}(u,\widehat{u} \cdot \sigma)|u|^\gamma, \qquad  {\overline{B}_\gamma}(u,\widehat{u} \cdot \sigma)=\max({B}(u,\widehat{u} \cdot \sigma),|\nabla_u  {B}(u,\widehat{u} \cdot \sigma)|)|u|^{\gamma+2},$$
and moreover, $\tilde{\re}_s(\cdot)$ is a given explicit restitution belonging to the class $\mathcal{R}_0$ for any $s \in [0,1]$. Finally, under the scaling $\re_{\l}(r)=\re(\l r),$ $r >0,\l >0,$ the constant $C_{0}(\re)$ is such that
$$C_{0}(\re_{\l}) \leq \mathfrak{a}_{0}(1-\gamma)\l^\gamma$$
where $\mathfrak{a}_{0}$ is defined in \eqref{eq:Je-e} for the restitution coefficient $\re(\cdot)$ belonging to the class $\mathcal{R}_{\g}.$
\end{prop}
We conclude this Section with the following bilinear estimate for $\Q_{\re}(f,g)$ in $L^{1}$-spaces with polynomial weights as derived in \cite[Lemma 4.1]{Tr}: 
\begin{prop}\label{prop:BIL}
Let $\re(\cdot)$ be a restitution coefficient in the class $\mathcal{R}_{0}$ and let $B(u,\sigma)=|u|b(\widehat{u}\cdot \sigma)$ with $b(\cdot) \in \W^{1,1}((-1,1))$. Let $\Q_{B,\re}$ be the associated collision operator. Then, for any $q >0,$ there is $C_{q}(\re) >0$ depending on $\re(\cdot),q$ and $b(\cdot)$ such that
$$\left\|\Q_{B,\re}(f,g)\right\|_{L^{1}_{v}(\m_{q})} \leq C_{q}(\re)\left(\|g\|_{L^{1}_{v}(\m_{q+1})}\|\|f\|_{L^{1}_{v}(\m_{q})}+\|g\|_{L^{1}_{v}(\m_{q})}\|\|f\|_{L^{1}_{v}(\m_{q+1})}\right)\,$$
for any $f,g \in L^{1}_{v}(\m_{q+1})$. Moreover, under the scaling $\re_{\l}(r)=\re(\l r),$ $r >0,\l >0,$ one has $\sup_{0 \leq \l \leq 1}C_{q}(\re_{\l}) < \infty.$
\end{prop}

\section{The restitution coefficient for viscoelastic hard spheres}\label{app:vis}

We illustrate our main assumptions on the restitution coefficient made in Section \ref{Sec21} (see Definition \ref{defiC})  to show that they are all met by the most physically relevant model in the literature. This model is the one of viscoelastic cases whose properties of the restitution coefficient have been derived in~\cite{PoSc,Poschel} by expliciting the dissipating force associated to viscoelastic spherical particles. We refer to \cite[Chapter 3, Eq. (3.22)]{Poschel}  for a derivation of the series  representation
\begin{equation}\label{viscE}
\re(r)=1+ \sum_{k=1}^\infty (-1)^k \, a_k \, r^{\frac{k}{5}}\,, \qquad \forall \,r > 0
\end{equation}
of the restitution coefficient $\re(r)$ in this case (see \eqref{visc}). We point out that the restitution coefficient actually admits the more tractable implicit representation
   \begin{equation}\label{visco}
   \re(r) + \mathfrak{a}_{0} r^{\frac{1}{5}} \re(r)^{\frac{3}{5}}=1\end{equation}
where $\mathfrak{a}_{0}>0$ is a suitable positive constant depending on the material viscosity. Notice that \eqref{visco} ensures that $\re(r) \in [0,1)$ for any $r >0$ with $\re(0)=1$.
The series representation \eqref{viscE} shows that the mapping $r \geq 0 \mapsto \re(r)$ is differentiable while its derivative $\re'$ is satisfying
$$\re'(r)+\frac{\mathfrak{a}_{0}}{5}\left(r^{-\frac{4}{5}}\re(r)^{\frac{3}{5}}+3 r^{\frac{1}{5}}\re'(r)\re(r)^{-\frac{2}{5}}\right)=0$$
thanks to \eqref{visco}. In particular
$$\re'(r)=-\frac{\mathfrak{a}_{0}}{5}r^{-\frac{4}{5}}\re(r)^{\frac35}\left(1+\frac{3}{5}\mathfrak{a}_{0}r^{\frac{1}{5}}\re(r)^{-\frac{2}{5}}\right)^{-1}$$
which shows that $r \geq 0 \mapsto \re(r)$ is nonincreasing and point \textit{(1)} of Definition \ref{defiC} is met. Multiplying the last inequality by $r$, one obtains
$$r\re'(r)=-\frac{\mathfrak{a}_{0}}{5}r^{\frac{1}{5}}\re(r)^{\frac{3}{5}}\left(1+\frac{3}{5}\mathfrak{a}_{0}r^{\frac{1}{5}}\re(r)^{-\frac{2}{5}}\right)^{-1},$$
so that the derivative $\eta_{\re}'(r)$ of the
 function $\eta_{\re}(r):=r\,\re(r)$ 
satisfies 
\begin{equation}\label{eq:35}\eta'_{\re}(r)=\frac{1 + \frac{2\mathfrak{a}_{0}}{5}r^{\frac{1}{5}}\re(r)^{-\frac{2}{5}}}{1 + \frac{3\mathfrak{a}_{0}}{5}r^{\frac{1}{5}}\re(r)^{-\frac{2}{5}}}\re(r) >0\end{equation}
showing \textit{(2)} of Definition \ref{defiC} is satisfied. From \eqref{visco}, it holds
\begin{equation*}
\re(r)^{\frac{3}{5}}=\left(\re(r)^{\frac{2}{5}}+\mathfrak{a}_{0} r^{\frac{1}{5}}\right)^{-1}\end{equation*}
and, since $\re(\cdot) \in [0,1)$, we deduce that \textit{(3)} of Definition \ref{defiC} holds true with $\re_{0}=0.$ Now, using again~\eqref{viscE}
$$
|\re(r)-1+{a}_{0}r^{\frac{1}{5}}|=\left|\sum_{k=2}^{\infty}(-1)^{k}a_{k}r^{\frac{k}{5}}\right|$$
from which one sees that~\eqref{gamma0bis} holds with $\mathfrak{a}_{0}=a_0$ and $\re(\cdot)$ belongs to the class $\mathcal{R}_{\g}$ with $\g=\frac{1}{5}$ and $\overline{\g}=\frac{2}{5}=2\g >\frac{3}{2}\g.$ Now, from \eqref{eq:35}, 
$$\eta'_{\re}(r) =\frac{\re(r) + \frac{2\mathfrak{a}_{0}}{5}r^{\frac{1}{5}}\re(r)^{\frac{3}{5}}}{\re(r) + \frac{3\mathfrak{a}_{0}}{5}r^{\frac{1}{5}}\re(r)^{\frac{3}{5}}}\re(r) =\frac{\re(r)+\frac{2}{5}(1-\re(r))}{\re(r)+\frac{3}{5}(1-\re(r))}\re(r) $$
and, since $0 \leq \re(r) \leq 1,$ this implies that 
$$\eta'_{\re}(r) \geq \frac{2}{5}\re(r)$$
proving the first part of \eqref{eq:Je-e} holds true with $\alpha=1$. For the second estimate in \eqref{eq:Je-e} one observes 
$$\eta'_{\re}(r)-\re(r)=\left(\frac{\re(r)+\frac{2}{5}(1-\re(r))}{\re(r)+\frac{3}{5}(1-\re(r))}-1\right)\re(r)=-\frac{1}{5}\frac{\re(r)(1-\re(r))}{\re(r)+\frac{3}{5}(1-\re(r))}$$
from which
$$|\eta'_{\re}(r)-\re(r)|=\frac{\re(r)}{3+2\re(r)}(1-\re(r)) \leq  1-\re(r)$$
proving the second part of \eqref{eq:Je-e}. 
To show \eqref{eq:Je-e1} one first observes that, from \eqref{visco} and the fact that $\lim_{r \to \infty}\re(r)=0$ one has
$$\re(r) \simeq \mathfrak{a}_{0}^{-\frac{5}{3}}r^{-\frac{1}{3}} \qquad \text{as} \qquad r \to \infty$$
so that $\eta_{\re}(r) \simeq \mathfrak{a}_{0}^{-\frac{5}{3}}r^{\frac{2}{3}}$ as $r \to \infty.$ Thus, the inverse bijection $\eta_{\re}^{-1}(z) \simeq \mathfrak{a}_{0}^{\frac{5}{2}}z^{\frac{3}{2}}$ as $z \to \infty$ showing that \eqref{eq:Je-e1} holds true with $m=\frac{3}{2}.$

 \bibliographystyle{plainnat-linked}

\end{document}